\documentclass[11pt,a4paper,twoside]{article}

\usepackage{microtype}
\usepackage[T1]{fontenc}
\usepackage{helvet}

\usepackage[margin=2cm, top=4.5cm, bottom=3cm, head=3cm]{geometry}
\usepackage{lastpage}

\usepackage{array}
\usepackage[british]{babel}
\usepackage[squaren]{SIunits}
\usepackage[numbers]{natbib}
\usepackage{bibentry}
\usepackage{url}
\usepackage{hyperref}
\usepackage{verbatim}
\usepackage{multirow}
\usepackage{subfigure}
\usepackage{algorithm,algorithmic}
\usepackage[bf]{caption}
\usepackage{setspace}
\usepackage{float}
\usepackage{placeins}
\usepackage{tikz}
\usetikzlibrary{shapes.misc}
\usetikzlibrary{shapes.arrows}
\usetikzlibrary{shapes.callouts}
\usetikzlibrary{backgrounds}
\usetikzlibrary{shadows}
\usetikzlibrary{calc}

\hypersetup{colorlinks=false}

\usepackage{tikz}
\usetikzlibrary{backgrounds}
\usetikzlibrary{calc}
\usetikzlibrary{shapes}

\usepackage[thmmarks]{ntheorem}
\usepackage{amsmath,amssymb}

\theoremsymbol{\ensuremath\blacksquare}
\theorembodyfont{\normalfont}

\theoremsymbol{\blitza}

\usepackage{parskip}
\allowdisplaybreaks
\usepackage{fancyhdr}
\pgfdeclareimage[height=1.5cm]{ucamlogo}{CUED_logo}

\usepackage{fancybox}
\title{\shadowbox{TECHNICAL REPORT: CUED/F-INFENG/TR.693}\\\vspace{\baselineskip}\bf Sequential Monte Carlo Optimisation for Air Traffic Management}
\author{A.J. Eele\thanks{\tt aje46@cam.ac.uk}{~}~and J. M. Maciejowski\thanks{\tt jmm1@cam.ac.uk}\\
~\\
\it(Cambridge University Engineering Department)
}
\begin{document}

\pagestyle{empty}
\begin{tikzpicture}[remember picture, overlay]

\path (current page.north west) -- ++(5.75,-9.5) coordinate (aaaa);

\draw[very thick, draw=blue!50!black] (aaaa) rectangle ++(9.5,-6) coordinate (bbbb);

\path (aaaa) -- node[pos=0.5,text width=8.5cm, text badly centered]{%
\fontfamily{pcr}%
\fontseries{b}
\fontsize{12}{15}
\selectfont%
{Sequential Monte Carlo Optimisation for Air Traffic Management}\\
\fontseries{m}%
\selectfont%
\vspace{\baselineskip}
A.J. Eele \& J.M. Maciejowski\\
\vspace{\baselineskip}
CUED/F-INFENG/TR.693\\
May 2015
} (bbbb);

\end{tikzpicture}

\cleardoublepage
\setcounter{page}{1}
\pagestyle{fancyplain}
\maketitle

\vfill
\begin{abstract}

This report shows that significant reduction in fuel use could be achieved by the adoption of `free flight' type of trajectories in the Terminal Manoeuvring Area (TMA) of an airport, under the control of an algorithm which optimises the trajectories of all the aircraft within the TMA simultaneously while maintaining safe separation. We propose the real-time use of Monte Carlo optimisation in the framework of Model Predictive Control (MPC) as the trajectory planning algorithm. Implementation on a Graphical Processor Unit (GPU) allows the exploitation of the parallelism inherent in Monte Carlo methods, which results in solution speeds high enough to allow real-time use. We demonstrate the solution of very complicated scenarios with both arrival and departure aircraft, in three dimensions, in the presence of a stochastic wind model and non-convex safe-separation constraints. We evaluate our algorithm on flight data obtained in the London Gatwick Airport TMA, and show that fuel saving of about 30\% can be obtained. We also demonstrate the flexibility of our approach by adding noise-reduction objectives to the problem and observing the resulting modifications to arrival and departure trajectories.


\vspace{3\baselineskip}
\it\scriptsize
This work was supported by EPSRC (Engineering and Physical Sciences Research Council - UK) Grant No. EP/G066477/1
\end{abstract}
\vfill

\cleardoublepage
\tableofcontents 

\cleardoublepage
\section{Introduction}
\subsection{Contribution of this report}

Civil aviation is considered to have been responsible for about 2\% of anthropogenic carbon dioxide $(CO_2)$ emissions in 1990 \cite{IPCC1999}, and this percentage is widely believed to have increased since then. This report shows that significant reduction in fuel use, and consequently in $CO_2$ emissions, could be achieved by the adoption of `free flight' type of trajectories in the Terminal Manoeuvring Area (TMA) of an airport, under the control of an algorithm which optimises the trajectories of all the aircraft within the TMA simultaneously while maintaining safe separation.

The report outlines the technical problem formulation and solution methodology, and describes the (ground-based) computing technology required to solve the problem quickly enough to make the approach feasible. It also discusses briefly changes in on-board equipment and procedures which would be required.

Typical approach and departure trajectories generated by the algorithm are presented and discussed, including special cases such as an aircraft landing with low fuel reserve.

A major part of the report is concerned with demonstrating that our proposed method would give substantial fuel savings in realistic scenarios. For this purpose we have taken traffic data from the vicinity of London Gatwick Airport over a 24-hour period, estimated the fuel used in actual approach and departure trajectories, and used our algorithm to generate alternative trajectories for the same aircraft. Examination of 4 hours out of the 24, including periods of both low and high traffic density, indicates that a fuel-use reduction of about 30\% can be achieved over this duration --- specifically, a reduction of over 14 tonnes of fuel (kerosene), equivalent to over 46 tonnes of $CO_2$. These results have been obtained with inclusion of a realistic representation of uncertainty, principally in the form of a random wind and turbulence field. Aircraft characteristics have been taken from the BADA database~\cite{BADA}.

One feature of our approach is that the objectives being optimised can be easily re-defined. For example, noise pollution over heavily-populated areas can be reduced by penalising trajectories that overfly such areas at low altitudes. We demonstrate this capability, with reference to population centres in the vicinity of London Gatwick Airport. More radical re-definition of objectives is also possible; for example, instead of reducing fuel use, the same approach could be used to increase traffic density and runway utilisation.

%

\subsection{Potential for fuel and emissions reductions in aircraft operations}

Work carried out by the Society of British Aerospace Companies and NATS in the UK has identified areas where savings in $CO_2$ emissions could be obtained~\cite{NATS09}. Savings of 10\% in $CO_2$ emissions were targeted from operational changes in air-traffic management. These savings were broken down between four flight modes:
\begin{itemize}
\item 3.25\% Climb: smooth continuous climb departures are more fuel efficient than stepped climb. Management of other aircraft to allow space for such climbing is needed.
\item 1.5\% Cruise: although generally well optimised, further savings can be found via improving matching of aircraft performance characteristics to route decisions.
\item 4.75\% Descent: continuous descent approaches, where aircraft approach the runway in a long glide, offer significant fuel savings when compared to descent in bursts with levelling out. Reduction in time spent holding also contributes significantly to this.
\item 0.5\% Ground Operations: reducing delays on stands and taxiways.
\end{itemize}
Of the planned savings, 80\% were to be gained from the climb and descent phases of each flight. In continuous climbs or descents the aim is to minimise the amount of time aircraft spend holding at fixed altitudes. In a congested airspace, particularly around the terminal manoeuvring area (TMA) of an airport, achieving this is extremely difficult for existing air-traffic management systems and procedures. This motivates our focus on operations within a TMA in this report.

\subsection{Trajectory Optimisation in a TMA}

Many existing proposals for multi-aircraft trajectory optimisation focus on cruise-like conditions where trajectories can be confined to two dimensions, or movement in the third dimension is highly discouraged for reasons of passenger comfort~\citep{HPS02,GT00,BP00}. Some work has been done using a novel conic approximation of the TMA to extend a 2-dimensional approach~\citep{PKMKZ08}. Methods such as mixed integer linear programming (MILP) require a linearisation of the problem and are too slow on centralised problems with large numbers of vehicles, due to the number of binary variables required~\citep{RH02}. 

Works more specific to the TMA often focus primarily on runway ordering optimisation through dynamic programming techniques~\citep{MPB10}. These techniques can include stochastic elements to account for uncertainty within the problem \citep{S12}. The work of \citep{K12} compares direct and indirect optimisation methods for a single aircraft's landing trajectory with the possible use of a dynamic programming technique to allow for real time operation. The work generated Shortest and Fastest Continuous Descent Approach (SF-CDA) paths able to reduce commercial aircraft annoyances from noise and fuel consumption.

This report tackles the TMA problem by combining the use of model predictive control (MPC) with sequential Monte Carlo (SMC) optimisation. In model predictive control an (open loop) optimal control problem, defined over a finite prediction horizon, is solved at each time step \citep{M02}. The solution is a trajectory of control actions, defined over the horizon. An initial segment of this solution is applied to the system being controlled (in our case, the aircraft in the TMA), then when new state measurements have been obtained a new solution is obtained and the process is repeated indefinitely.  More details of MPC are given in section \ref{sec:MPC}.

A paramount concern of air-traffic management is the maintenance of safe separation of aircraft. Separation constraints are inherently non-convex: there may exist two safe trajectories, each going either side of an exclusion region for example, but a convex (weighted-average) combination of these trajectories is not necessarily safe, because it may penetrate the exclusion region. Optimisation problems with non-convex constraints (and/or non-convex objective functions) are known as `non-convex' problems, and are inherently difficult to solve~\citep{BV04}, because they can have multiple local optima, and a `go downhill' search strategy may well result in becoming stuck at one of those local minima.

Stochastic optimisation algorithms address the problem of non-convex optimisation. They do so by `learning' from exploration where good solutions lie in the search space and concentrating further exploration there, but with occasional large jumps to probe whether a more promising region should be explored. Well-known examples of stochastic optimisation methods are `simulated annealing'~\citep{KirGelVec83} and `genetic algorithms'~\citep{Gol89}, while a less well-known one is `sequential Monte Carlo'~\citep{DFG01}, which is the one we employ in this study. A `nice to have' feature of most stochastic optimisation algorithms is that they guarantee to find a solution arbitrarily close to the global optimum (or set of global optima, in case there is more than one), \emph{if} one runs the algorithm for sufficiently many iterations. A further and important advantage of stochastic optimisation methods is that the solutions found have a degree of robustness in cases where the problem to be solved itself contains some randomness (in our case: random wind and turbulence, uncertain aircraft parameters, etc).

The use of stochastic optimisation, and specifically SMC algorithms, in the context of MPC and with application to planning trajectories for aircraft and UAV's has previously been proposed in~ \citep{KML08,EM11,VGS11}. In these and other studies it has been noted that conventional implementation of the SMC algorithm resulted in computation times which are too long for real-time use. However, Monte Carlo methods are known to have a highly-parallelisable structure, and very
significant speed-ups for statistical, chemical and economic applications have been demonstrated through the use of Graphical Processor Units (GPU's)~\citep{LYGDH10,FVR07,R00,SPFHTS07}.
%
Previous work by the authors also demonstrated a 98\% computational speed-up for implementing SMC on a GPU for \emph{en-route} conflict resolution
~\citep{EMCL13}. 

\section{Model Predictive Control}
\label{sec:MPC}
Model Predictive Control (MPC), also known as Receding Horizon Control, is an optimisation-based control strategy. A constrained, finite-horizon, optimal control problem is solved online at each iteration based on the available state/measurements. This solution yields the controls for the duration of a finite horizon.

\begin{figure*}[htbp]

\begin{center}
\begin{tikzpicture}

\tikzstyle{mpcblockmstr} = [draw, very thick, rounded corners, drop shadow, fill=white, font=\footnotesize, text width=1.75cm, text badly centered, minimum height=1cm]

\tikzstyle{mpcblock} = [mpcblockmstr]

\tikzstyle{mpcblockr} = [mpcblockmstr]
\tikzstyle{mpchl} = [fill=white];
\tikzstyle{mpcline} = [draw, very thick, ->];

\node[mpcblockr] (COST) at (0,0.5) {Objective function};
\node[mpcblock, left of=COST, node distance=2.5cm] (MODEL) {Prediction model};
\node[mpcblock, right of=COST, node distance=2.5cm] (CON) {Constraints};
\path (COST.south) -- +(0,-2) coordinate (MID1);

\node[mpcblock, node distance=2.5cm] at (MID1) (OPTPROB) {Optimal control problem};
\node[mpcblock, mpchl, right of=OPTPROB, node distance=2.5cm] (OPTALG) {Optimisation algorithm};

\node[mpcblock, right of=OPTALG, node distance=3cm] (PLANT) {Plant};

\node[mpcblock, above of=PLANT, node distance=1.5cm] (DIST) {Disturbances};

\path (OPTPROB) -- +(1.25,-1.5) coordinate (MID2);

\node[mpcblockr] at (MID2) (OBS) {Observer};

\draw[mpcline] (MODEL) -- (OPTPROB);
\draw[mpcline] (COST) -- (OPTPROB);
\draw[mpcline] (CON) -- (OPTPROB);
\draw[mpcline] (PLANT) |- (OBS);

\draw[mpcline] (OPTPROB) -- (OPTALG);
\draw[mpcline] (OPTALG) -- (PLANT);

\draw[mpcline] (OBS) -- +(-4,0) |- (OPTPROB);
\draw[mpcline] (DIST) -- (PLANT);

\begin{pgfonlayer}{background}
\draw[rounded corners, dashed, fill=yellow!20] (-4,-1) rectangle (4,-4.5);
\node[above right=0.25cm] at (-4,-4.5) {``Online tasks''};

\draw[rounded corners, dashed, fill=green!20] (4.25,-4.5) rectangle (6.75,0.5);
\node[above right=0.25cm] at (4.25,-4.5) {``Real world''};

\draw[rounded corners, dashed, fill=blue!20] (-4,-0.75) rectangle (4,1.25);
\node[above right=0.25cm] at (-4,-1) {``Setup''};

\end{pgfonlayer}

\end{tikzpicture}
\end{center}

\caption{The components of a Model Predictive Control (MPC) feedback system}
\label{fig:mpc_block_diagram}

\end{figure*}
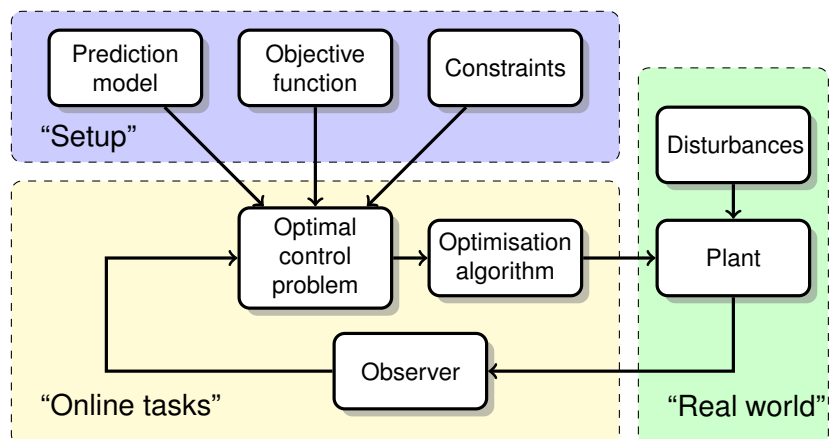


%


Figure \ref{fig:mpc_block_diagram} shows the components of a typical Model Predictive Control (MPC) feedback system~\cite{M02}. At the top, labelled ``Setup'', are the components that define the controller: a \emph{model} which will be used to generate the predictions needed by the optimiser, an \emph{objective function} which will be minimised by the optimiser, and \emph{constraints} which need to be satisfied by the solution generated by the optimiser. Below this the block labelled ``Online tasks'' shows what needs to be executed at each sampling/update interval:
\begin{enumerate}
\item Measurements are obtained from sensors located on the plant which is being controlled, and used by an observer to estimate the plant state vector (not always necessary).
\item A specific optimal control problem is defined (which depends on the current (estimated) plant state, using all the ingredients specified in the ``Setup''.
\item This optimal control problem is solved by an optimisation algorithm. Usually this computes a whole sequence of future input commands. Only the first of these is applied to the plant being controlled. Note that the optimisation problem is time-critical; it needs to be solved within the control-update interval.
\end{enumerate}
On the right of Fig.\ref{fig:mpc_block_diagram}, labelled ``Real world'', is the plant being controlled. Disturbances acting on this plant are shown explicitly, in order to emphasise that the plant never behaves exactly as predicted by the model.

In this work's application to air-traffic management, the `plant' will be the whole set of $N$ aircraft in the airspace being controlled. The state of each aircraft will be a vector of $n=6$ components: its position in 3-D space, its airspeed, its heading, and its mass (the latter varying because of fuel consumption). The control applied to each aircraft will be a vector of $m=3$ components: the engine thrust, the climb angle and the bank angle. Thus the complete `plant' will have a state vector with $nN$ components and an input vector with $mN$ components. The disturbances acting on the `plant' will be primarily the wind/turbulence acting on each aircraft, but these can also represent the effects of modelling errors when making predictions.

The prediction model used consists of standard, relatively simple, nonlinear difference equations for each aircraft due to the time discretisation within the problem. If measurements of all the aircraft states are taken at time step $k$, then we use our model to provide predictions at the following $H$ steps $k+1,k+2,\ldots,k+H$. $H$ is called the `prediction horizon'. These predictions depend on predicted values of the $m$ controls applied to each aircraft over $H$ steps, and it is the role of the optimiser to find the best values of these. Many alternative forms of optimisation have been applied to MPC and often, in simpler systems, linear or quadratic programming can be applied. In this report Sequential Monte Carlo optimisation is used as the method for solving for the applied control sequence.

\section{Sequential Monte Carlo Optimisation}
A Sequential Monte Carlo method approximates the probability distribution $p(\xi_k|\zeta_1,\zeta_2,\ldots,\zeta_k)$ of some (vector) variable $\xi_k$ at time step $k$, conditional on observations $\zeta_1,\zeta_2,\ldots,\zeta_k$ which have been made up to that time. It assumes that $\xi_k$ is the value of a Markovian random process: $p(\xi_k|\xi_1,\xi_2,\ldots,\xi_{k-1})=p(\xi_k|\xi_{k-1})$, and that the probability of an observation $p(\zeta_k|\xi_k)$ depends only on $\xi_k$. It does this by propagating, by simulation, a large number $L$ of samples of $\xi_k$, and approximating the distribution by a `histogram' of their values. Each of these samples is often referred to as a `particle', hence an alternative name for SMC is `particle filter', and in general it estimates the state of a `hidden Markov model'~\cite{DFG01},\cite[chapter~14]{RC99}. 

When SMC is used to solve the MPC problem described in section \ref{sec:MPC}, a `hidden state' $\xi_k$ is the optimal solution $\mathbf{u}_k^*$, and an `observation' $\zeta_k$ is a value $J_{k:k+H}(x_k,\mathbf{u}_k)$ of the cost that results from a trial solution $\mathbf{u}_k$, together with the corresponding set of aircraft trajectories $\hat{x}_{k+j|k}$, $j=1,\ldots,H$ --- these are needed to check whether any constraint violations are predicted to occur. $L$ needs to be chosen sufficiently large for the approximate distribution of $\xi_k$ to be sufficiently concentrated about the global optimum (or set of global optima) to be useful --- that is, that any sample drawn from this distribution is very likely to be very close to a global optimiser. See~\cite{VGS11} for a precise account of how an augmented statistical problem can be set up such that an optimisation problem is re-cast as an inference problem, in a context very similar to the one considered in this report.

\subsection{Algorithm}
\begin{algorithm*}[htb]
\caption{Sequential Monte Carlo}
\label{alg:SMC}
\begin{algorithmic}[1]
\STATE $J \leftarrow$ 0
\STATE Define a monotonically increasing SampleSchedule of length $J_{\max}$
\STATE Clone all aircrafts' current states to give each particle its own copy.
\STATE Set the weights of all aircraft to $1/L$ where $L$ is the number of particles.
\STATE For each particle and aircraft in the particle randomly generate controls over the entire horizon $H$ \label{alg:initial}
\WHILE {$J\leq J_{\max}$}
\STATE $j \leftarrow$ 0
\FOR {each particle $l$}
\WHILE {$j\leq$ SampleSchedule(J)} \label{alg:inner1}
\STATE Sample disturbance realisations for all agents and all time steps to horizon $H$.\label{alg:disturb}
\STATE Simulate trajectories for all agents till horizon $H$.\label{alg:plan}
\IF {an aircraft $i$ fails constraints} \label{alg:constraint}
\STATE aircraft's $i$'s weight set to 0.
\ENDIF
\STATE Calculate each aircraft's cost.
\STATE Scale each aircraft's weightings by their cost. \label{alg:weighting}
\STATE $j \leftarrow j+1$
\ENDWHILE \label{alg:inner2}
\ENDFOR
\STATE $J \leftarrow J+1$
\IF {$J \leq J_{\max}$}
\STATE Resample all particles \label{alg:resample}
\STATE Peturb all aircrafts'' controls with Gaussian white noise.\label{alg:perturb}
\STATE Set the weights of all agents to $1/L$ where $L$ is the number of particles.
\ENDIF
\ENDWHILE
\STATE Draw final sample from particle population. \label{alg:final}
\end{algorithmic}
\end{algorithm*}

Algorithm \ref{alg:SMC} summarises the implementation of SMC. Each step of the algorithm shall now be dealt with in further depth with specific reference to the application of multiple aircraft trajectory planning with conflict avoidance.

\subsection{Particles}
\label{sec:particle}
A particle represents a single instance of the problem with the data of all associated agents. The number of particles for optimisation is a design variable dependent on the complexity of the problem to be optimised. The more complex the problem, the larger the number of particles may need to be to adequately characterise the search space. Inside each particle there is:
\begin{itemize}
\item An initial state for each of the $N$ individual aircraft.  These initial states for aircraft are the same across all particles and are provided to the optimisation at the beginning of the algorithm.
\item The control inputs for every step of the MPC horizon $H$ for every aircraft. In this report's formulation there are 3 controls for each aircraft (thrust, bank angle and pitch angle) totalling $N*H*3$ control inputs. These control inputs are different for every particle.
\item A separate weighting for each of the $N$ individual aircraft in the scenario, used for resampling.
\end{itemize}  

The controls inside particles are initialised as random samples from a uniform distribution for each control between the maximum and minimum value. After this initialisation at line~\ref{alg:initial} the controls are then updated following resampling and perturbation in lines~\ref{alg:resample}-~\ref{alg:perturb}. The perturbation moves particles locally in the search space, whilst resampling causes particles to cluster around areas of the search space with positive attributes as determined by the objective function. The weights of all aircraft are initialised as $1/L$ where $L$ is the number of particles and these are updated in line~\ref{alg:weighting}.

\subsection{Trajectory Planning and Disturbance Realisations}
In lines~\ref{alg:disturb}-\ref{alg:plan} the future trajectory of each aircraft in each particle is simulated given the initial state and controls (from the particle's data) and the disturbance samples. The discretised aircraft dynamics model is a standard 6 DOF model with 6 states and 3 inputs
\begin{subequations}
\begin{align}
x_i(k+1)=&x_i(k)+\delta t (v_{s,i}(k) \cos(\chi_i(k))\cos(\gamma_i(k)))+w_x(k)\delta t \\
y_i(k+1)=&y_i(k)+\delta t(v_{s,i}(k) \sin(\chi_i(k))\cos(\gamma_i(k)))+w_y(k)\delta t \\
z_i(k+1)=&z_i(k)+\delta t(v_{s,i}(k) \sin(\gamma_i(k)))\\
v_{s,i}(k+1)=&v_{s,i}(k)+\delta t(\frac{T_i(k)-D_i(k)}{m_i(k)}-g\sin(\gamma_i(k))) \\
\chi_i(k+1)=&\chi_i(k)+\delta t(\frac{L_i(k)\sin(\phi_i(k))}{m_i(k)v_{s,i}(k)})\\
m_i(k+1)=&m_i(k)-\delta t(\eta T_i(k))
\end{align}%
\label{eqn:hor}%
\end{subequations}%
where: $\delta t$ is the step length; $x, y$ and $z$ are the Cartesian coordinates of the aircraft ($z$ acting as the altitude of the aircraft); $m$ is the total mass of the aircraft; $v_s$ is the true airspeed; $\chi$ is the heading angle; $\gamma$ the climb angle and $\phi$ the bank angle. The subscript $i$ denotes the $i^{\mathrm{th}}$ aircraft. The control variables are $\phi$ (bank), $T$ (thrust) and $\gamma$ (climb). Lift $L$ and drag $D$ are calculated using the standard aerodynamic relations. The fuel usage is controlled by the constant $\eta$.

Each aircraft is given an initial state $X_{0,i}=(x_{0,i},y_{0,i},z_{0,i},\chi_{0,i},v_{s,0,i},m_{0,i})$. Departing aircraft are modelled from an initial altitude of 400m, where the aircraft is assumed to have taken off from the runway and completed its initial climb. A departing aircraft is considered as having left the TMA once it has reached a set distance $D_{TMA}$ from the airport at which point it is assumed to enter a neighbouring control zone and is handled by alternative ATC.

\begin{figure}[ht]
\begin{center}
\def\picwidth{14cm}
\includegraphics[width=\picwidth,viewport=1.7in 4.0in 10.8in 6.0in,clip=true]{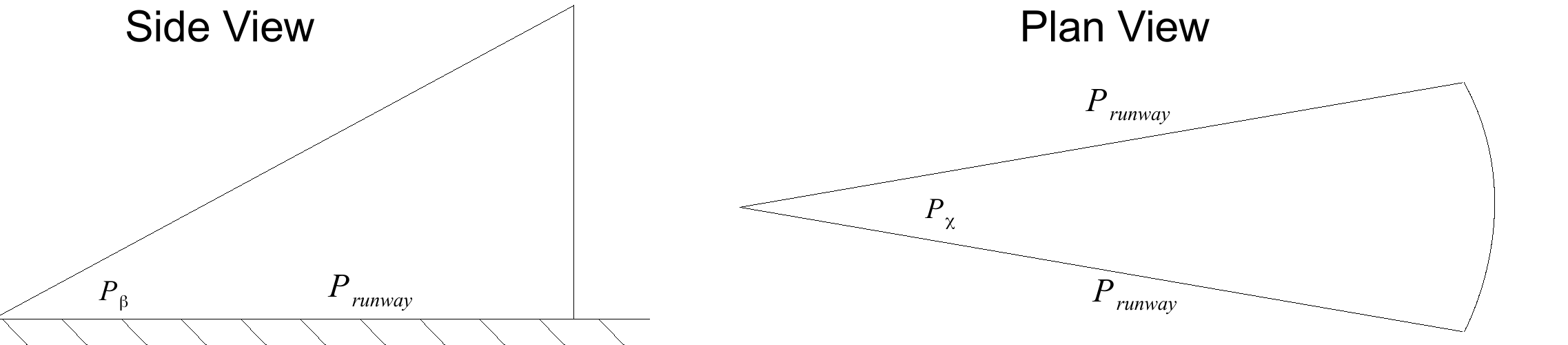}
\caption{Landing Envelope Shape and Size}
\label{fig:envelope}
\end{center}
\end{figure}

Arrival aircraft are allowed the flexibility to enter the problem at any point along the TMA's boundary. Arrival aircraft are deemed to have `landed' or passed from our control once they have satisfied the following constraints:
\begin{subequations}
\begin{align}
\sqrt{x_i(k)^2+y_i(k)^2}\leq& P_{\mathrm{runway}}\\
\beta_i(x_i(k),y_i(k),z_i(k))\leq& P_{\beta}\\
|\arctan(y_i(k)/x_i(k))|\leq&  P_{\chi}\label{eqn:cone}\\
|\chi_i(k)-180|\leq&  P_{\chi}\label{eqn:heading}\\
v_{s,i}(k)\leq& P_{v_s}\label{eqn:landspeed}
\end{align}
\label{eqn:TO_init}
\end{subequations}
The first 3 of these constraints define the landing sector, an arc with radius $P_{\mathrm{runway}}$ horizontal angle $ P_{\chi}$ and vertical angle  $P_{\beta}$ where $\beta_i$ is described later in Equation~\ref{eqn:flow}. Constraint \ref{eqn:heading} ensures that the aircraft is pointing towards the runway assuming an East to West landing without significant crosswind. Finally constraint \ref{eqn:landspeed} limits the airspeed ($P_{v_s}$)at which the aircraft can enter the landing sector for a successful landing. These are shown pictorially in Figure~\ref{fig:envelope}.

The disturbance realisations added to the system are the primary method of simulating uncertainty within the system. This uncertainty can come from wind, sensor noise, controller noise and human factors. Within the inner loop of optimisation (lines \ref{alg:inner1}-\ref{alg:inner2}) each aircraft within a particle will have drawn $\mathrm{SampleSchedule}(J)$ disturbance realisations and planned the same number of trajectories. This inner loop serves to simulate different disturbance scenarios on an aircraft to determine if the controls are valid with respect to constraints, and their fitness with respect to the objective function. This information is stored in the aircraft's weight within the particle as described in the next two subsections. In applications with no uncertainty there would be no need for repeated planning and simulation of the aircraft and the inner loop would only be executed once. In this report we have restricted the disturbances to being wind and this is discussed in greater depth in Section~\ref{sec:wind}.

\subsection{Constraint Handling}
\label{sec:constraints}
There are two types of constraint handled by line~\ref{alg:constraint}. The first are unary constraints which deal only with one aircraft at a time. Examples of these include flight envelope constraints and the minimum aircraft mass constraint.
The flight envelope constraints provide upper and lower bounds for both the state and control variables:
\begin{subequations}
\begin{align}
z_{\mathrm{min}}\leq z(t)&\leq z_{\mathrm{max}}\\
-\gamma_{\mathrm{max}}\leq \gamma(t)&\leq \gamma_{\mathrm{max}} \\
-\phi_{\mathrm{max}}< \phi(t) &< \phi_{\mathrm{max}}\\
T_{\mathrm{min}} \leq T(t) &\leq T_{\mathrm{max}}\\
v_{\mathrm{min}} \leq v_s(t) &\leq v_{\mathrm{max}}.
\end{align}
\end{subequations}
The minimum aircraft mass constraint enforces that the total aircraft mass must always remain above the mass of the aircraft alone without fuel:
\begin{equation}
m(t)\geq m_{\mathrm{Aircraft}}.
\end{equation}
This acts on a constraint on the total fuel use across the entire modelled flight. Further unary constraints can be added. For example to create exclusion zones where aircraft are not allowed to fly, or must fly above a certain altitude. These are not included in the existing model but can be implemented in future work.

The second type of constraint are the binary constraints. These hold between two aircraft in the same particle, such as in conflict avoidance. In this report the protection zone around each aircraft is modelled as a cylinder with horizontal radius $P_r$ and altitude separation of $P_h$. Two aircraft $i$ and $j$ avoid each other if they satisfy the following constraint for every time step in the MPC horizon:
\begin{align}
\nonumber (x_{i}(k)-x_{j}(k))^2+(y_{i}(k)-y_{j}(k))^2 &\geq (2P_r)^2 \\
\nonumber \vee |z_{i}(k)-z_{j}(k)| &\geq 2P_h \\
\forall k, \forall i,j &\in \{1,...,N\}: i\neq j.
\label{eqn:avoidance}
\end{align}

If an aircraft fails any of its unary constraints then its weight stored in the particle is set to 0. If a pair of aircraft fail a binary constraint then both aircraft have their weights set to 0. Any aircraft with a weighting of 0 will not be resampled in the resampling phase of the algorithm described in Subsection~\ref{sec:resample}.

The constraint handling method described in this section is a departure from the the generic formulation for a SMC optimisation. In the generic formulation a weight is linked to a single particle which would contain the entire simulation of all $N$ aircraft. In that case if a single aircraft fails any constraint the entire simulation would be weighted as 0 and all controls from that particle discarded when the particles are resampled in the next iteration. Conversely, by using our implementation, aircraft are weighted individually inside a particle. If a single aircraft fails a constraint the non-failing aircraft weights would still be non-zero. This would allow their controls to continue to the next iteration. The application under consideration is a highly constrained situation and this alteration gives greater flexibility in keeping valid control solutions for aircraft which would otherwise have been discarded. This in turn allows a smaller number of particles, as the distribution they are required to model is less complex than in the case of a multi-aircraft weighted particle.

\subsection{Particle Costing}
The SMC method requires a non-negative function as the objective function $J_{T}$ which is to be maximised. This objective function is made up of multiple elements and is different for arrival and departure aircraft due to their distinct priorities. However each element of the objective functions for both types of aircraft takes a basic form where we have negated the original minimisation function then normalised by adding the supremum of the function and dividing by the difference between the supremum and the infimum of the function. The supremum and infimum are readily available via geometric arguments as will be explained in greater depth case by case.

\subsubsection{Departure:}
Departing aircraft are given a desired target altitude $Z_{i}$,  a desired airspeed $v_{s,D}$  and bearing $\theta_{i,F}$. These are used to evaluate the aircraft's progress in the objective function. However it is not necessary for the aircraft to have reached these before leaving the TMA as it can take a while to reach cruise altitude and traffic in the TMA could make reaching the exact target bearing unrealistic. 

The departing aircrafts' objective function $J^D_{T,i}$ is made up of four elements: reaching the target bearing they need to achieve for onwards flight to their destination; a cost on fuel used; rewarding the aircraft for climbing towards their required cruise height and rewarding the aircraft for maintaining its desired speed. 
\begin{subequations}
\begin{eqnarray}
J^D_{1,i}(k:k+H)=&\frac{-A_i(k:k+H)+ \mathrm{sup} A_i(k:k+H) }{\mathrm{sup}A_i(k:k+H)- \mathrm{inf}A_i(k:k+H) }\\
J^D_{2,i}(k:k+H)=&\frac{(m_i(k)-m_i(k+H))+ \delta t H T_{\mathrm{max}}\eta}{\delta t H T_{\mathrm{max}}\eta}\\
J^D_{3,i}(k:k+H)=&\frac{-B_i(k:k+H)+ \mathrm{sup} B_i(k:k+H) }{\mathrm{sup}B_i(k:k+H)- \mathrm{inf}B_i(k:k+H)}\\
J^D_{4,i}(k:k+H)=&\frac{-C_i(k:k+H)+ \mathrm{sup} C_i(k:k+H) }{\mathrm{sup}C_i(k:k+H)- \mathrm{inf}C_i(k:k+H)}
\end{eqnarray}
\end{subequations}
where:
\begin{subequations}
\begin{align}
A_i(k:k+H)=&\sum^{k+H}_{j=k} \frac{|\theta_i(j)-\theta_{i,F}|}{H}\\
B_i(k:k+H)=&\sum^{k+H}_{j=k} \frac{\sqrt{(z_{t_f,i}-z_i(j))^2}}{H}\\
C_i(k:k+H)=&\sum^{k+H}_{j=k} \frac{|v_{s,i}(j)-v_{s,D}|}{H}
\end{align}
\end{subequations}
A is the absolute angle between the desired bearing of the aircraft leaving the airport and the actual bearing of the aircraft from the airport. B is the distance between the desired altitude for cruise and the altitude of the aircraft. C is the difference between the desired airspeed and actual airspeed of the aircraft. The sup and inf of A are approximated as the maximum deflection from the desired bearing, 180 degrees and the minimum deflection 0 degrees respectively. Similarly for B the sup is the maximum distance the aircraft can be from its target altitude given its current altitude, obtained by either climbing or descending with $\gamma_{\mathrm{max}}$ and $V_{\mathrm{max}}$. The inf is then the minimum distance the aircraft can be from its target altitude similarly using the $\gamma_{\mathrm{max}}$ and $V_{\mathrm{max}}$ until it is at the desired altitude at which point it would be 0. By using this formulation, the cost for aircraft $i$ which is on its desired bearing for all steps of the horizon would have a cost of 1 for $J^D_{1,i}$, while an aircraft $j$ which is 90 degrees away from its desired bearing for all steps of the horizon would have a cost of 0.5 for $J^D_{1,j}$. The cost $J^D_2$ is the fuel minimisation recast as a maximisation; $\delta t H T_{\mathrm{max}}\eta$ is the maximum amount of fuel that could be used in a step, the minimum is approximated to be 0.

The four parts of the cost $J^D$ are all normalised and then scaled in importance by the weighting factors $\alpha_j$. These weighting coefficients have been tuned empirically and have a large effect on the behaviour of aircraft. Choice of $\alpha_j$ affects policy e.g. fuel use, noise abatement etc.

\begin{subequations}
\begin{eqnarray}
J^D_{T,i}(k:k+H)=&\sum^4_{j=1} \frac{\alpha_j J^D_{j,i}(k:k+H)}{H} \label{eqn:maxcost}\\
0 \leq& J^D_{T,i}(k:k+H) \leq 1 \\
\sum^{4}_{j=1}\alpha_j =& 1
\end{eqnarray}
\end{subequations}

\subsubsection{Arrival}
In contrast to departures, arriving aircraft have far stricter terminal constraints with a specific requirement that to land the aircraft must be: travelling towards the runway; be in line with the runway given a heading tolerance; within a certain distance of the runway; travelling below a given speed; and be within a set altitude trajectory to accomplish landing. The model does not extend to fully landing the plane and continues to assume pilot control for this aspect of the flight, similar to that of the initial take off.

\begin{figure}[ht]
\begin{center}
\def\picwidth{10cm}
\includegraphics[width=\picwidth,viewport=0.7in 0.1in 6.3in 4.8in,clip=true]{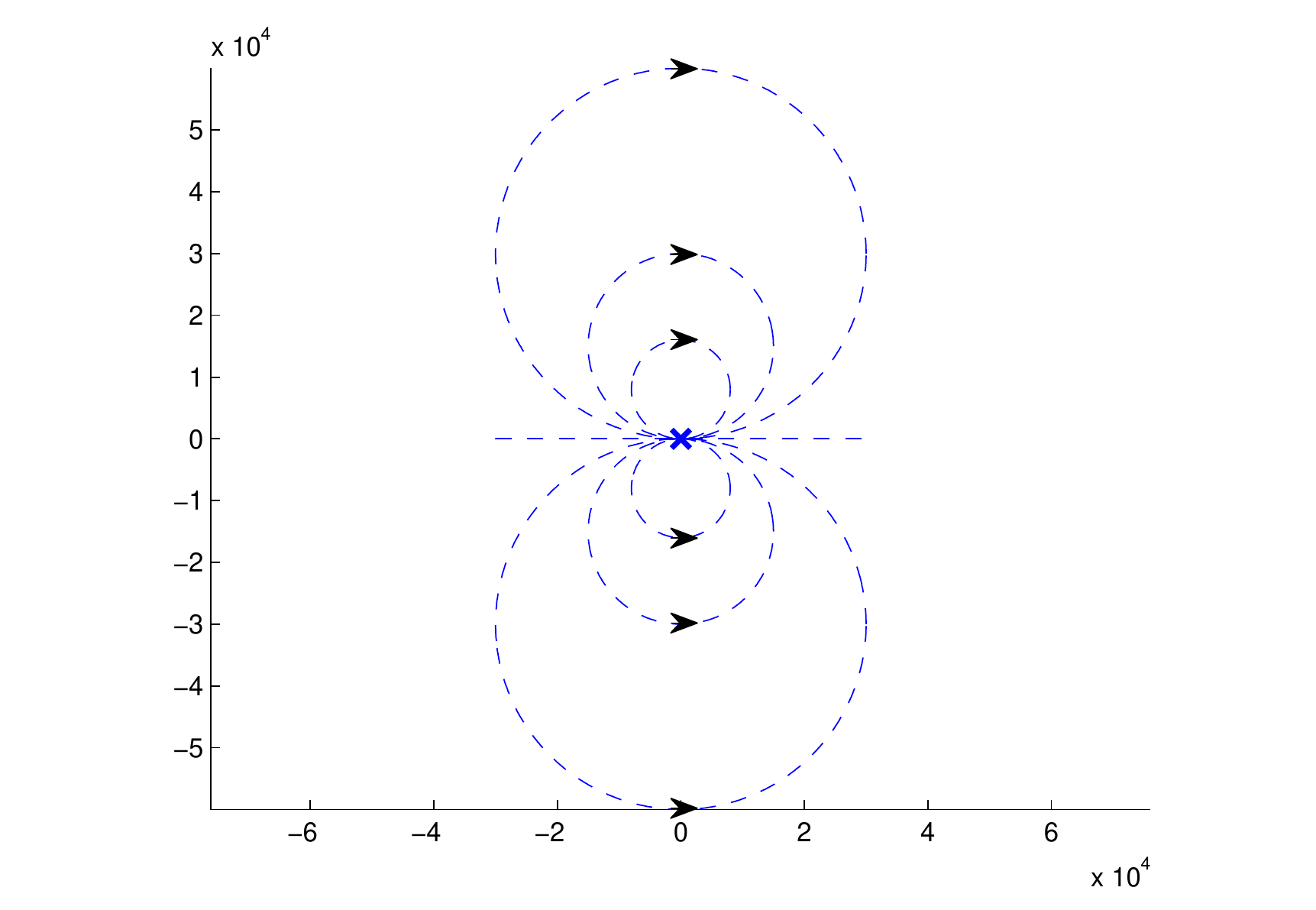}
\caption{Example flowfield for a single runway airport with East to West landing pattern}
\label{fig:flowfield}
\end{center}
\end{figure}

Arriving aircraft will head towards the airport following an objective function made up of 3 elements. The first is a cost on the fuel used to encourage minimal fuel use; the second rewards the aircraft for following a nominal altitude descent trajectory; and the third rewards the aircraft for keeping its heading close to a nominal flow field designed to encourage the aircraft to turn smoothly before approaching the runway. Such a flow field is entirely designable based on the airport or goals of the ATC for the airport. For this work the airport flow field has been formulated as shown in Figure~\ref{fig:flowfield} and the cost function penalises aircraft for their heading not aligning with the flow field. The use of a nominal altitude descent trajectory rather than relying purely on the minimisation of fuel or a target altitude is justified by the finite horizon of planning which might not be long enough for an aircraft to reach the landing sector. By providing a nominal trajectory the aircraft's descent is smoother and less likely to drop in altitude sharply before having to fly close to the ground till it reaches the landing sector.
The arrival problem is then optimised for the maximisation of the cost $J^A_{T}$ which is the sum of the weighted costs for heading offset from flowfield, fuel usage and nominal altitude descent trajectory following:
\begin{subequations}
\begin{eqnarray}
J^A_{1,i}(k:k+H)=&\frac{-D_i(k:k+H)+ \mathrm{sup} D_i(k:k+H) }{\mathrm{sup}D_i(k:k+H)- \mathrm{inf}D_i(k:k+H) }\\
J^A_{2,i}(k:k+H)=&\frac{(m_i(k)-m_i(k+H))+ \delta t H T_{\mathrm{max}}\eta}{\delta t H T_{\mathrm{max}}\eta}\\
J^A_{3,i}(k:k+H)=&\frac{-E_i(k:k+H)+ \mathrm{sup} E_i(k:k+H) }{\mathrm{sup}E_i(k:k+H)- \mathrm{inf}E_i(k:k+H) }
\end{eqnarray}
\end{subequations}
where:
\begin{subequations}
\begin{align}
D_i(k:k+H)=&\sum^{k+H}_{j=k} \frac{|\chi_i(j)-\hat{\chi}_{i,f}(x_i(j),y_i(j))|}{H}\\
E_i(k:k+H)=&\sum^{k+H}_{j=k} \frac{|\beta_i(x_i(j),y_i(j),z_i(j))-\beta_{i,f}|}{H}
\end{align}
\end{subequations}
and
\begin{equation}
\hat{\chi}_{i,f}(x_i(j),y_i(j))=2\arctan\left (\frac{x_i(j)}{y_i(j)}\right )+90\\
\end{equation}
\begin{equation}
\beta_i(x_i(j),y_i(j),z_i(j))=\arctan\left(\frac{z_i(j)}{180-2\arctan\left(\frac{x_i(j)}{y_i(j)}\right)\frac{||x_i(j),y_i(j)||}{\cos(\arctan(x_i(j)/y_i(j)))}      }\right)
\label{eqn:flow}
\end{equation}
$D$ is the normalised difference between the flow field's heading and the aircraft's, $E$ is the normalised difference between the angle to the ground origin for the aircraft and the desired angle of approach. The geometric terms of $\beta$ and $\hat\chi$ are based on the flow field shown in Figure~\ref{fig:flowfield} with an east west runway centred on the origin with arrival flights landing on the eastern end of the runway. This example flow field is used throughout the report for its elegant geometric properties. It is simple to change the flow field and would result in a change of equation~\ref{eqn:flow}. The distance from the aircraft to the origin for $E$ is measured not as the Euclidean distance but rather as the distance remaining on the arc of the flow field the aircraft is currently on. The point of this design choice is to encourage a constant descent angle for the aircraft throughout its travel to the runway. The Euclidean distance would have lead to variations in the angle of descent over the full path for aircraft starting at positions greater than 90 degrees to the desired landing bearing.

Like the departures, the arrivals objective function is designed to be non-negative and the total cost $J^A_T$ and all its constituent costs are normalised between 0 and 1.  The weightings on the terms of the cost function $\tilde\alpha$ are subject to priorities of the desired trajectory.
\begin{subequations}
\begin{eqnarray}
J^A_{T,i}(k:k+H)=&\sum^3_{j=1} \frac{\tilde\alpha_j J^A_{j,i}(k:k+H)}{H} \label{eqn:maxcost}\\
0 \leq& J^A_{T,i}(k:k+H) \leq 1 \\
\sum^{3}_{j=1}\tilde\alpha_j =& 1
\end{eqnarray}
\end{subequations}

\subsection{Weight Scaling}
\label{sec:weighting}
As previously mentioned in subsection~\ref{sec:particle} each aircraft in a particle has a weight associated with it. This weight records the degree of success an aircraft has in the simulations with various noise realisations. To do this in line~\ref{alg:weighting} the cost of aircraft $i$ in particle $l$ is multiplied by the weight of aircraft $i$ in particle $l$ each time the inner loop \ref{alg:inner1}-\ref{alg:inner2} is executed. An aircraft only needs to violate one constraint in any execution of the inner loop to have its weighting set to 0. This weighting would remain as 0 until resampling (described in the next section) removes that aircraft's proposed controls from the population.

Ideally at each iteration of the inner loop the sum across all particles $L$ of weights associated with aircraft $i$ would be normalised to 1. This would reduce numerical issues as the inner loop progresses. Without such a normalisation step the weight would always be decreasing
\begin{subequations}
\begin{align}
W_{i,l}^{j+1}&=W_{i,l}^j J_{T,i,l}^j\\
W_{i,l}^0 &= 1/L \\
0 &\leq J_{T,i,l}^j \leq 1
\end{align}
\end{subequations}
where $W_i^j$ represents the weight and $J_{T,i,l}^j$ the objective function value of aircraft $i$ in particle $l$ at inner loop iteration $j$ and $L$ is the number of particles. The cost is always between 0 and 1 and the weight starts at less than 1.

\subsection{Particle Resampling and Perturbation}
\label{sec:resample}
In lines \ref{alg:resample} - \ref{alg:perturb} the aircraft are individually resampled to generate a new population of particles. These new particles  have their aircraft controls perturbed by Gaussian white noise to separate particles with similar controls across the local search space. This perturbation also allows particles to explore away from points which were in the original randomly sampled population of controls.

Resampling is done on the basis of sequential importance sampling using a method such as Kitagawa resampling~\citep{K96}. The higher a weight aircraft $i$ has in particle $l$ compared to the weight of aircraft $i$ in all other particles, the more likely it is to have its controls resampled into the new population.

As the aircraft are resampled separately to form the new particles no conflict avoidance guarantees can hold once resampling has taken place. This arises from the phenomena that the controls from aircraft $i$ from particle $l$ could now be in a new particle with the controls from aircraft $j$ from particle $m$. This again arises from the departure from the generic method of one weight per particle mentioned in Section~\ref{sec:constraints}, justified by the significant simplification of dimensionality of the search space and the effect on both computation time and particle populations needed. As long as the final sample from the distribution of particles is drawn correctly with regards to the binary conflict constraints the lack of conflict guarantees at the beginning of each inner loop is not problematic. 

\subsection{Final Sample Selection}
This is the final step of the SMC optimisation (line~\ref{alg:final}) before it hands back to the MPC update process. This step determines the controls to be used by all aircraft for the next update step. In many generic SMC applications this final sample is taken as the mean of the values of the particles. In multi modal distributions this would give an inaccurate estimation of the global optimisers (consider the mean of a control where half of the population steers left around an obstacle and the other half right). An alternative statistical measure to offset this problem would be to select the mode of the distribution. However this is also unsuitable in our implementation. The presence of binary constraints across aircraft control distributions could lead to constraint violation. Therefore throughout this report we select the best performing particle in the final iteration as our estimate of the maximiser.

To find the best performing particle the weights of all aircraft in the particle are multiplied together and the particle with the greatest overall weight is selected as the estimator of the maximiser
\begin{equation}
\max_{l} \prod_{i=0}^N W_{i,l}.
\end{equation}
Any particle where an aircraft has failed any constraint in simulation would have a weight 0 associated with the failing aircraft and thus could not be selected as the best performing particle.

\subsection{Rolling Window}
To allow for longer scenarios where aircraft both enter and leave the problem, as would happen in a real airport, the whole application is considered in the form of a rolling window. The length of the window is the same length of the horizon $H$ for the MPC and the window moves with every completion of an MPC update. Any aircraft entering the problem (either through take off or landing) does so at a pre-determined time step, though this can be relaxed to consider uncertainty on arrival into the problem. Only aircraft which are considered active in the problem have their trajectories optimised by the algorithm previously described. 

As we are using MPC, aircraft may become active in the problem midway through a horizon. Their trajectories are optimised from the point they enter the problem. Similarly for arrival aircraft there is the chance for an aircraft to complete mid horizon, thus leaving the problem. During the algorithm when trajectories are planned (lines~\ref{alg:disturb}-\ref{alg:plan}) arrival aircraft active in the problem are tested to see if they've entered the landing sector. If they have they will be noted as finished and future time steps within that horizon will not be planned for that aircraft. Departure aircraft are less constrained and not tested for early completion during the trajectory planning stage. Arrival aircraft which have met their landing constraints mid horizon are given the best possible cost, 1, for all remaining steps of the horizon as a bonus to that aircraft for landing.

\begin{figure}[htb]
\begin{center}
\def\picwidth{14cm}
\includegraphics[width=\picwidth,,viewport=0in 0.5in 10.5in 8.4in,clip=true]{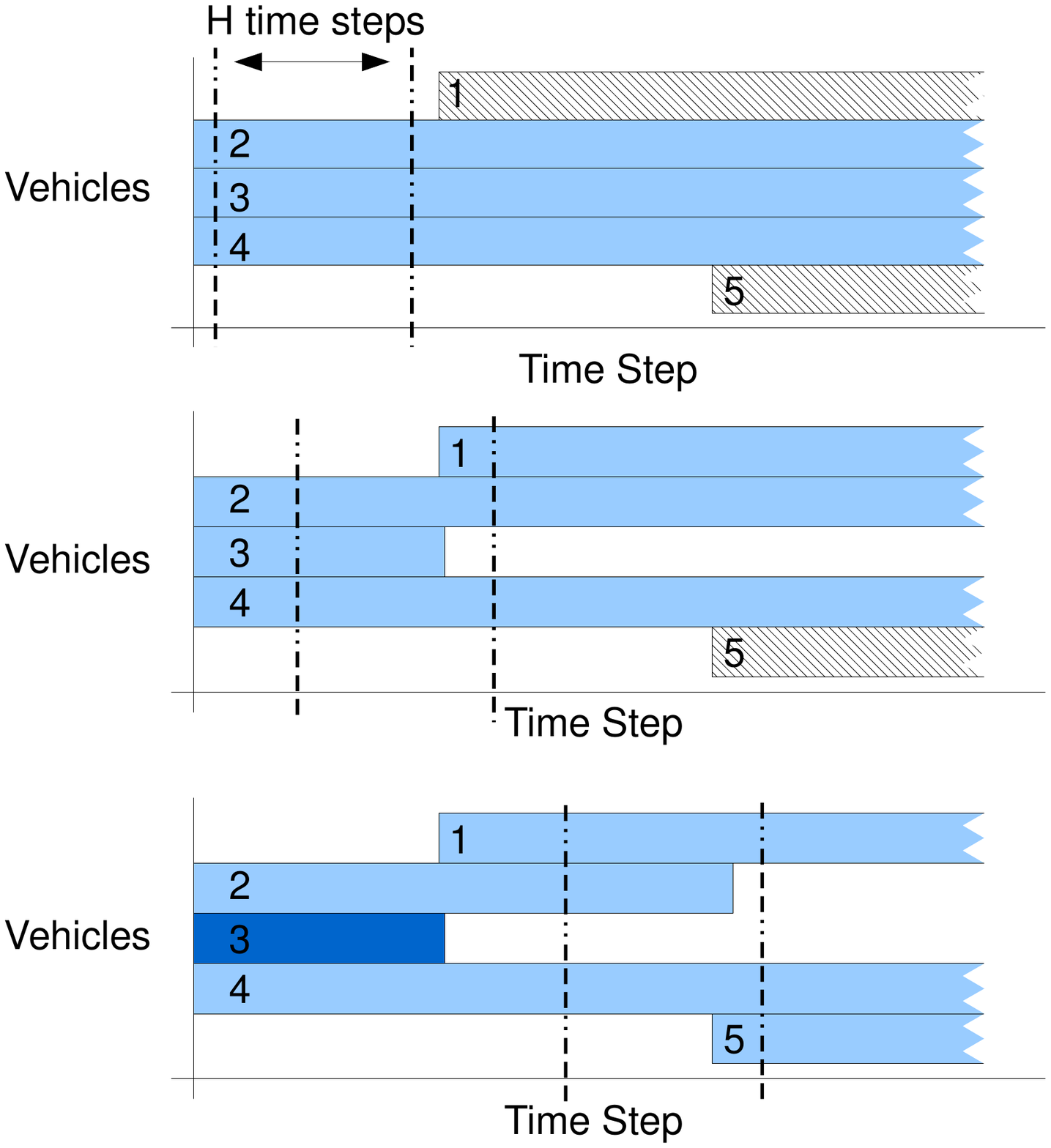}
\caption{Graphical representation of active aircraft (light blue), completed aircraft (dark blue) and aircraft which haven't joined the problem (hatched) for rolling window simulation at three different points in simulation}
\label{fig:rollingwindow}
\end{center}
\end{figure}
Figure~\ref{fig:rollingwindow} shows an example problem where aircraft enter and leave the problem during simulation. Active aircraft are shown in light blue, finished aircraft which have left the problem are shown in dark blue and aircraft yet to enter the problem are shown in hatched lines. Aircraft will not be deactivated until they no longer appear within the bounds of the rolling window (MPC horizon), likewise they are not activated until they have appeared in the MPC horizon. The aircraft which have not finished are shown planning for an unknown time into the future as finishing can only be detected within the MPC horizon.

\subsubsection{Wind Model}
\label{sec:wind}
Large portions of the uncertainty about aircraft flight plans results from the inherent uncertainty in meteorological forecasts. Wind speed and direction is probably the most important meteorological factor affecting aircraft trajectories in the TMA. The wind disturbance is assumed to be made up of two components, a nominal component representing the wind forecast and a stochastic component representing the forecast errors. This report uses a spatio-temporal wind model for the forecast errors based on the model used by Lymperopoulos and Lygeros \citep{LL10}. Their model approximations were designed for wind prediction over a far larger geographical scale than that used by the TMA example. However, they appear to be suitable for our application, following experimentation.

The model of Lymperopoulos and Lygeros assumes an isotropic wind-field (invarient under rotations) with uncorrelated wind speeds in the South-North and East-West directions. Vertical wind is neglected in the model however correlation between wind at different altitudes is considered. The wind field is generated from a random field where each point is Gaussian with zero mean and a covariance matrix $R(t,P,t',P')$ where $t$ and $t'$ are the points in time and $P$ and $P'$ are the points in three dimensional space. Through these assumptions and by using an exponential function the correlation can be approximated as:
\begin{equation}
\rho_{xy}(t,P,t',P')=\sigma(z)\sigma(z')\mathrm{exp}(-\lambda|t-t'|)\mathrm{exp} \left(-\beta \left |\begin{array}{ccc}
{x-x'}\\
{y-y'} \end{array}\right | \right )\mathrm{exp}(-\gamma| z-z'|).
\label{eqn:cov}
\end{equation}
Where $\sigma(z)$ is the standard deviation of the wind error in m/s at altitude $z$. The figures used for $\lambda, \beta$ and $\gamma$ used by Lygeros and Glover\citep{LG04} were $\lambda=6*10^{-6}m^{-1}$, $\beta=1.6*10^{-6}m^{-1}$ and $\gamma=1.5*10^{-5}m^{-1}$ along with $\sigma(z)$ are based on the data reported in \citep{CRKB98}.

To use the correlation function above Lymperopoulos and Lygeros proposed the following method for generating wind realisations.
The wind field is divided into a grid with $N_x$ points in the South-North direction, $N_y$ points in the East-West direction and $N_z$ vertically. For each point in the three dimensional grid 2 random numbers are generated from a zero mean Gaussian distribution. The first of these random numbers in each pair relates to the South-North wind error and the second to the East-West wind error, these numbers are stored in two vectors $W_X(k)$ and $W_Y(k)$ at time step $k$.

The covariance matrix $\hat{R}\in \mathbb{R}^{N_X N_Y N_Z \times N_X N_Y N_Z}$ where each entry is a comparison between two grid points using equation~\ref{eqn:cov} is the covariance matrix for both $W_X(k)$ and $W_Y(k)$ given the isotropic assumption. Lygeros and Lymperopoulos assume this covariance matrix to be constant in time and thus generate wind samples using the following linear Gaussian model:
\begin{eqnarray}
\begin{array}{cc}
W_X(0)=\hat{Q}v_X(0), & W_X(k+1)=aW_X(k)+Qv_X(k+1),\\
W_Y(0)=\hat{Q}v_Y(0), & W_Y(k+1)=aW_Y(k)+Qv_Y(k+1),\\
\end{array}
\end{eqnarray}
where $v_X(k),v_Y(k)\in\mathbb{R}^{N_X N_Y N_Z}$ are standard independent Gaussian random variables, $Q$ and $\hat{Q}$ are derived by Cholesky Decomposition from the covariance matrix $\hat{R}$ using
\begin{equation}
\begin{array}{ccc} QQ^T=(1-a^2)\hat{R} &\mathrm{and} & \hat{Q}\hat{Q}^T=\hat{R}\end{array}
\end{equation}
where $a=e^{-\delta t/G_t}$ with $G_t$ a parameter of time correlation obtained from \citep{CRKB98}. The wind error for individual aircraft is then obtained by tri-linear interpolation between the grid points.

\section{Parallelisation}
\label{sec:parallel}

It had been previously been proposed there is a large degree of speed up that can be obtained by parallelising the SMC algorithm\citep{LYGDH10,FVR07,R00,SPFHTS07}. The authors previous work~\citep{EMCL13} focused on implementation on a graphics processing unit (GPU) using the CUDA programming language provided by NVIDIA and demonstrated a 98\% computational improvement. CUDA is a scalable parallel programming language similar to C and C++ with the libraries and utility to design kernels to implement on a GPU. These kernels are repeatable functions representing the code that the user wishes to be executed in parallel with different data. GPUs can use thousands of threads at the same time, whereas a conventional computer  will only use a few. The execution and near instantaneous switching between threads is how a GPU achieves high efficiency on parallel applications.

\subsection{GPU Structure}
A GPU differs from a standard CPU in the number and type of cores it has. A CPU will have a few large cores optimised for sequential operations whilst a GPU will have many hundreds of smaller cores efficient at running different tasks in parallel. In GPU programming the structure of the problem is very a data driven model of computation. Typically each thread executes the same operation on different elements of the data in parallel. Threads have an id number which is used to compute memory addresses for where that thread's data is stored and the id number can also be used for controlling decisions inside the kernel. To simplify handling many thousands of threads and to take advantage of potential structures in the data threads are grouped together into blocks containing up to 512 threads (larger numbers of threads can be handled in blocks with advances in the programming architecture) in 1D, 2D or 3D arrays. These blocks can then also be structured in grids. Threads can also bundled together into blocks for cooperation purposes, if threads need to share information they can use a shared memory resource accessible by all threads in the same block. This shared memory is both quite small by comparison to the global GPU memory but is also significantly faster as it is based on-chip with the computation cores themselves. Therefore the primary decisions when implementing an existing algorithm in parallel are: firstly identifying the code to parallelise inside a kernel; and deciding on intelligent use of shared memory for the threads executing that kernel. 

\subsection{Kernel Design}

\begin{figure}[htb]
\begin{center}
\def\picwidth{14cm}
\includegraphics[width=\picwidth,viewport=0.1in 1in 10in 6.8in,clip=true]{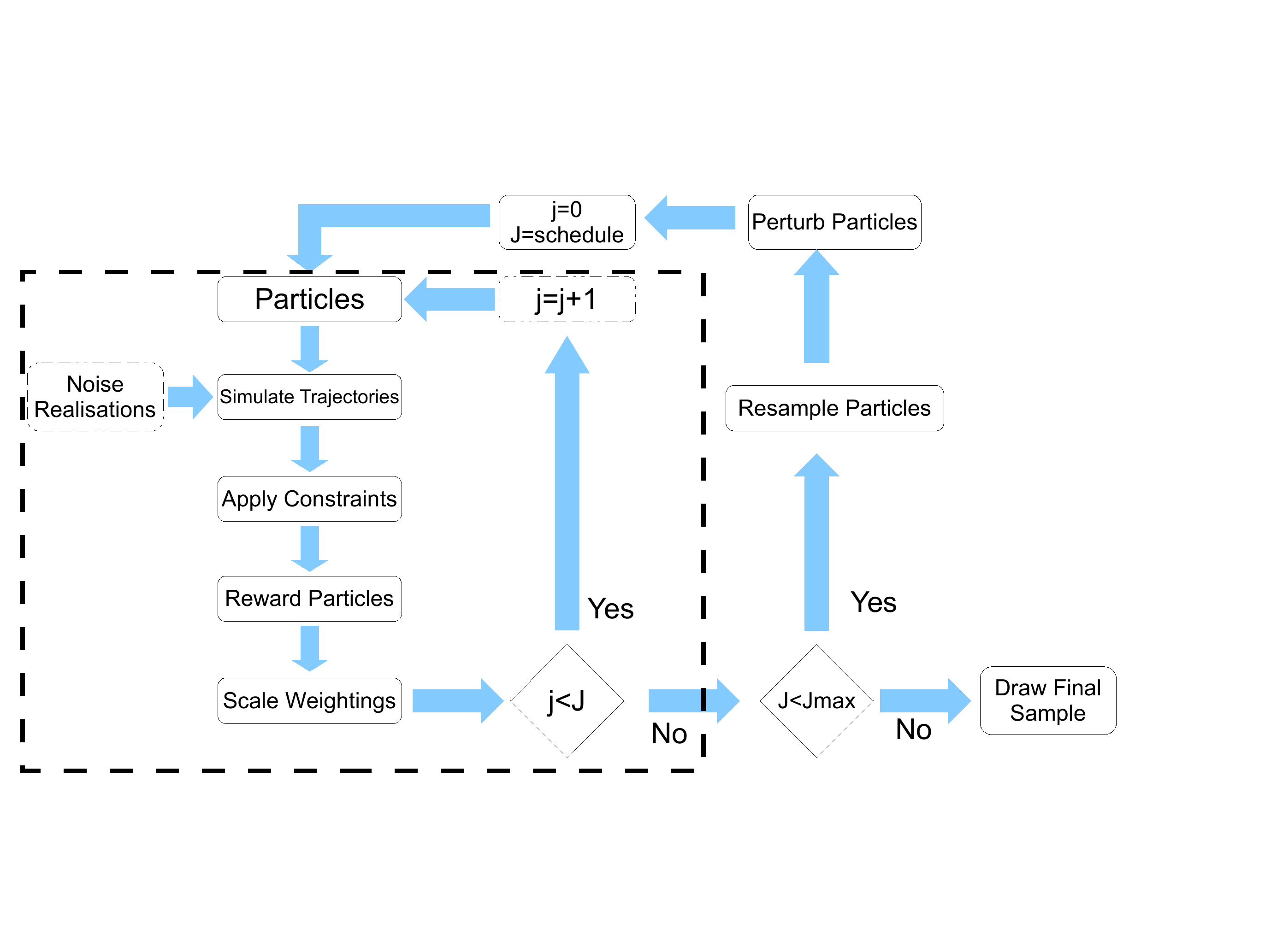}
\caption{Graphical Representation of SMC Algorithm with Kernel Highlighted}
\label{fig:paralleldiag}
\end{center}
\end{figure}

Figure \ref{fig:paralleldiag} shows a simplified graphical representation of the SMC algorithm focusing on the operations taking place in a single computation step of the MPC. There is also a boxed region of the flowchart demonstrating the GPU kernel. This kernel incorporates the inner loop of simulation specifically the lines \ref{alg:inner1}-\ref{alg:inner2} of the algorithm. This code was a bottleneck in a sequential implementation as it must be executed multiple times for each particle where the only difference between each execution is the particle's individual data. This made it a perfect candidate for parallelisation with one particle per thread, as there is minimal thread interaction and use of shared memory. 

In the algorithm used there is no need for the threads to synchronise or communicate and thus the problem can be considered as `embarrassingly parallelisable'. The remaining implementation decision was to select an appropriate random number generator for use within the kernel. In this report's implementation we have focused on using the NVIDIA provided library CURAND, specifically the XORWOW (Xor Shift added with Weyl sequence generator~\citep{MX03}) which seeds each thread and maintains the state of each random generator between kernel calls. This generator was used for its superior speed compared to that of the others included in CURAND.

\subsection{Computational Savings}
\label{sec:compsave}
Previously the MPC-SMC method was demonstrated for fast computation of many aircraft problems in cruise-like conditions\citep{EMCL13} with computation time per MPC update step varying up to 33 seconds for a 20 vehicle problem. This implementation had the simplification of fixed scenario lengths with all aircraft entering the problem at the start and no aircraft leaving the problem before the end. As such all aircraft were active during the entire scenario and computation time per step of the algorithm was deterministic based on the number of aircraft in the problem. This is not true for this report's implementation where the number of active aircraft in the problem varies step by step as aircraft both enter and leave the scenario.

\begin{figure}[htb]
\begin{center}
\def\picwidth{12cm}
\includegraphics[width=\picwidth]{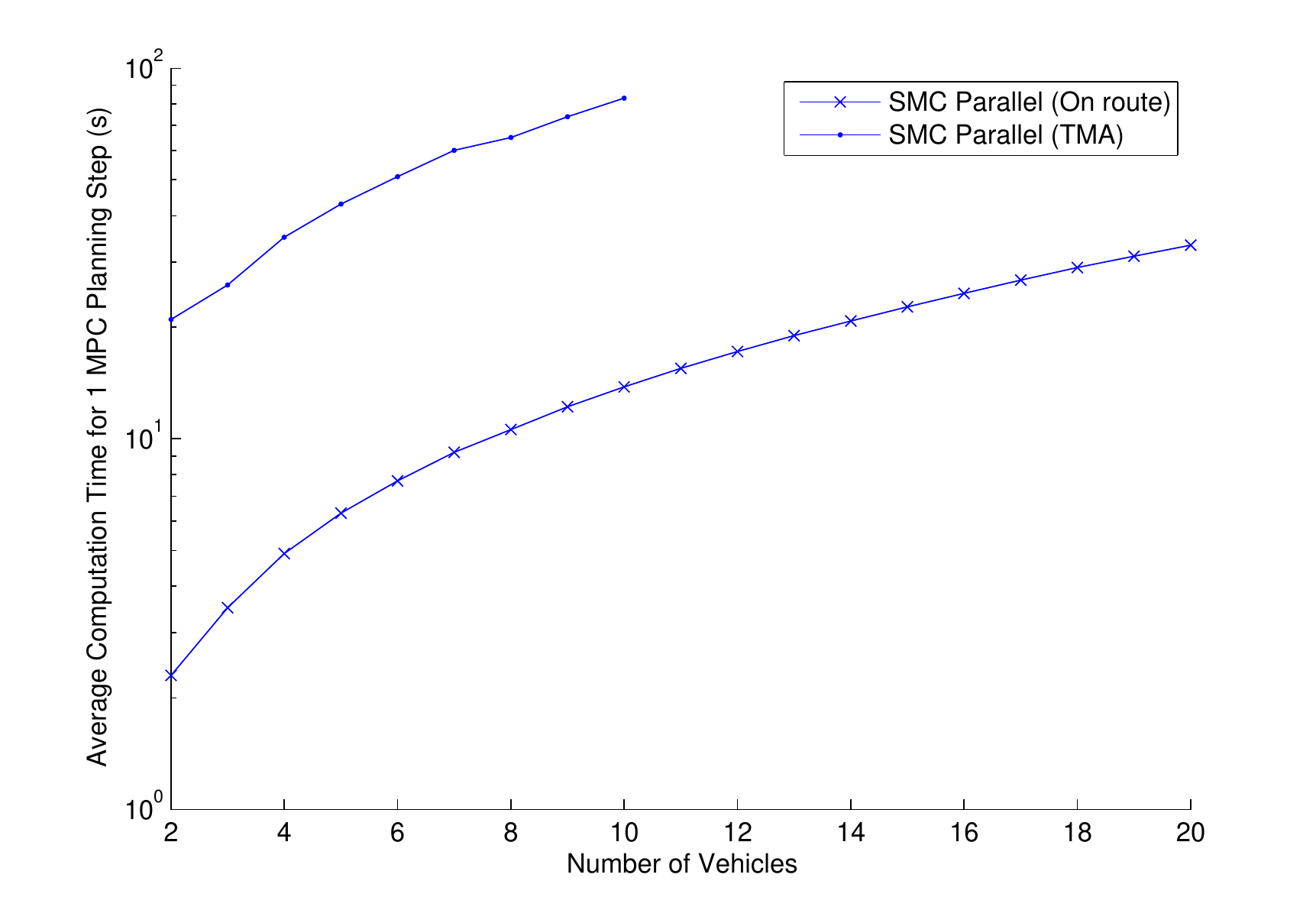}
\caption{Comparison of previous and current computation time spent per time step for problems with different number of vehicles}
\label{fig:times}
\end{center}
\end{figure}

Figure~\ref{fig:times} shows the comparison of our previous implementation's time per update step compared to the average of this report; fewer points are available for this report's implementation due to the nature of the scenarios being solved however it offers an insight into the magnitudes involved.

This implementation is significantly slower than the previous incarnation, most notably due to the 10-fold increase in number of particles needed to adequately represent the search space. The problem of landing vehicles into a small terminal area whilst obeying all flight constraints is significantly harder than the previous cruise like situation and 1024 particles was not enough to reliably find feasible solutions. A larger or faster GPU would offer direct benefit to the computation speed through allowing more cores to process the particle threads and each core to work faster.
\subsection{Structural Considerations}
Given the layout of parallelisation previously defined is such that no individual thread of simulations needs to communicate with any other thread during their time in the kernel there is a high level of flexibility on how we lay out our threads and blocks for computation. In GPU computation individual runs of the kernel are threads which are then bundled together into blocks for structural reasons. All threads in a block are run on the same microprocessor on the GPU and the current maximum number of threads a block can hold is 512. Threads inside the same block can make use of low latency memory onboard their specific microprocessor, which is important if they need to communicate or synchronise (not needed in our application). Groups of blocks are held together in a grid and there is no real limit on the number of blocks that can be held in a grid. There can be more blocks than there are microprocessors on the GPU and in this case blocks which cannot be immediately assigned a processor are held in queue until they can be.

With a fixed number of particles used throughout our TMA simulations ($10240 = 5*2^{11}$) there are 17 combinations of threads and blocks which directly factorise the total number of particles and do not exceed the maximum number of threads one block can hold. These 17 combinations are shown in Table \ref{tab:comp_sense}. The computation speed of a time step is mostly dependant on how many aircraft are in the problem and how many of them are active at at a given time. Therefore to test the variation of computational speed based on the layout of threads in blocks it was key to set up a standard problem where a fixed number of aircraft were always active for enough time steps to observe the fixed computational speed. The set problem used to test computational speed was chosen as a 10 aircraft problem where all aircraft entered the problem at the 5th MPC time step (with a horizon of 6 steps so aircraft were optimised from time step one over increasing parts of the horizon until the 5th step was reached). The times for the first 4 time steps were disregarded as aircraft were not active across the entire horizon and thus the computation time varied step to step. After the 5th time step regardless of structural layout fixed time step lengths were observed. The average computation time for a given structural layout was then obtained by averaging 3-4 time steps after the settling period.

\begin{table}[h]
  \centering
\begin{tabular}{ c | c | c  }
  \hline                        
No. Threads in Block & No. Blocks & Average Computation Time (s)\\
\hline
  1 & 10240 & 1497\\
  2 & 5120 & 761\\
  4 & 2560 & 394\\
  5 & 2048 & 330\\
  8 & 1280 & 201\\
  10 & 1024 & 180\\
  16 & 640 & 109\\
  20 & 512 & 98\\
  32 & 320 & 72\\
  40 & 256 & 70\\
  64 & 160 & 60\\
  80 & 128 & 61\\
  128 & 80 & 58\\
  160 & 64 & 57\\
  256 & 40 & 58\\
  320 & 32 & 56\\
  512 & 20 & 64\\
\hline
\end{tabular}
\caption{Average Computation Time for a single update in different Thread*Block layouts}
  \label{tab:comp_sense}
\end{table}

\begin{figure}[htb]
\begin{center}
\def\picwidth{12cm}
\includegraphics[width=\picwidth]{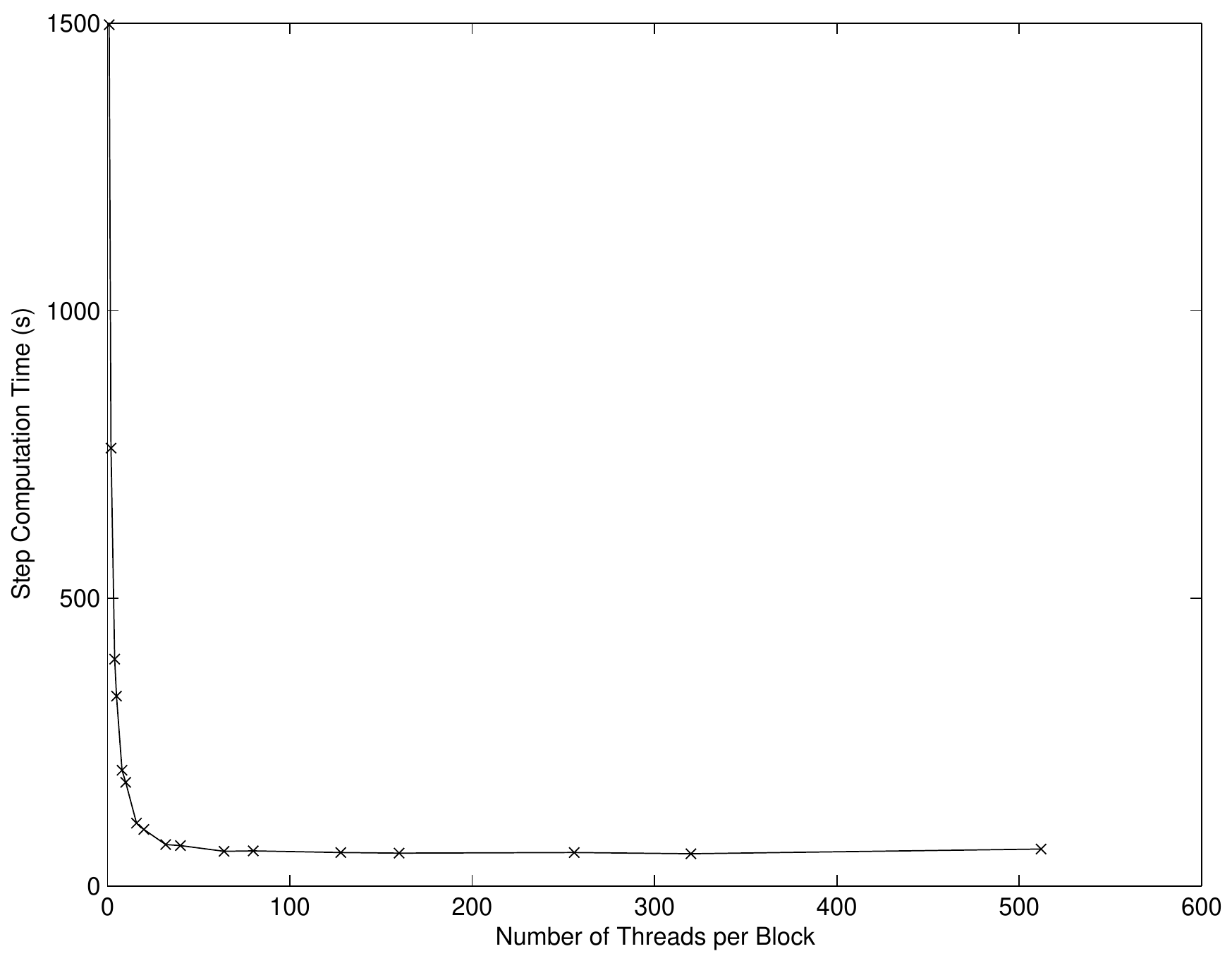}
\caption{Average Computation Time for a single update in different Thread*Block layouts}
\label{fig:comp_sense}
\end{center}
\end{figure}

All tests were done on the same GPU unit and it would be expected that individual GPU setups (specifically the number of cores) would affect the outcome of the average computational time and best thread*block layout. The GPU used for these tests was an NVIDIA GeForce GTX 580 which has 512 cores. Table \ref{tab:comp_sense} and Figure \ref{fig:comp_sense} show the average step computation time in different layouts. Clearly sending a single thread per block to one of the GPU's processors suffers too greatly from the overheads of memory writing and instructing the processor. More surprisingly sending 20 threads to every processor on the GPU is also still quite slow. The best observed step times seem to strike a balance to sending a large group of threads to each processor (potentially reducing the effect of computational overhead) vs using more of the GPU simultaneously. 

\section{Simulation Results}
This section presents the results of simulations performed on the previously outlined algorithm and application setup. These simulations vary in complexity from single aircraft paths, used to tune the performance of the objective function coefficients to tests designed to see how far we can push the implementation with our current setup. The final section of the results discusses the computational speed of the method compared to that of a previous implementation\citep{EMCL13}.

\subsection{Simulation Setup}
For all simulations within this report the following parameters were kept constant. MPC time steps were 10 seconds in length ($\delta t = 10$) and had a horizon length of $H=6$ time steps. $\delta t$ is arguably far shorter than needed in the application where standard ATC would be looking at 30 seconds to 60 seconds updates. For simplicity aircraft have been standardised with the aerodynamic properties of an Airbus A320 obtained from BADA\citep{BADA}. The TMA is considered as a circle of radius 30km where there is a single runway at the origin. Aircraft both land and depart East to West from this one runway. There is no additional runway scheduler running alongside the simulations so separation near the landing envelope is maintained by the distance separation used for conflict avoidance from Equation~\ref{eqn:avoidance}.

The algorithm used a $J_{\max}=100$, with a schedule function of $\mathrm{SampleSchedule}(J)=\lfloor(3+5e^{0.05J})\rfloor$ and $L=10240$ particles (which have been naively arranged as 80 blocks of 128 threads for the GPU implementation). The simulations were performed on an NVIDIA Tesla C2070, which has 448 cores.

The disturbances from wind for all scenarios were simulated using the wind model outlined in Section~\ref{sec:wind}. A three dimensional grid of eight points was used to generate the wind at the outlying reaches of the TMA before trilinear interpolation found the disturbances at the aircrafts' individual locations. No prevailing wind was used in the scenarios so the aircraft are only subject to the random disturbances of the wind field itself. A denser grid of points can be used for sampling the wind speed however this slows computation down significantly with the need to generate many more random numbers within each inner loop of the GPU kernel.

\subsection{Single Aircraft Trajectories}
\label{sec:SAT}
\begin{figure}[htb]
\begin{center}
\def\picwidth{12cm}
\includegraphics[width=\picwidth]{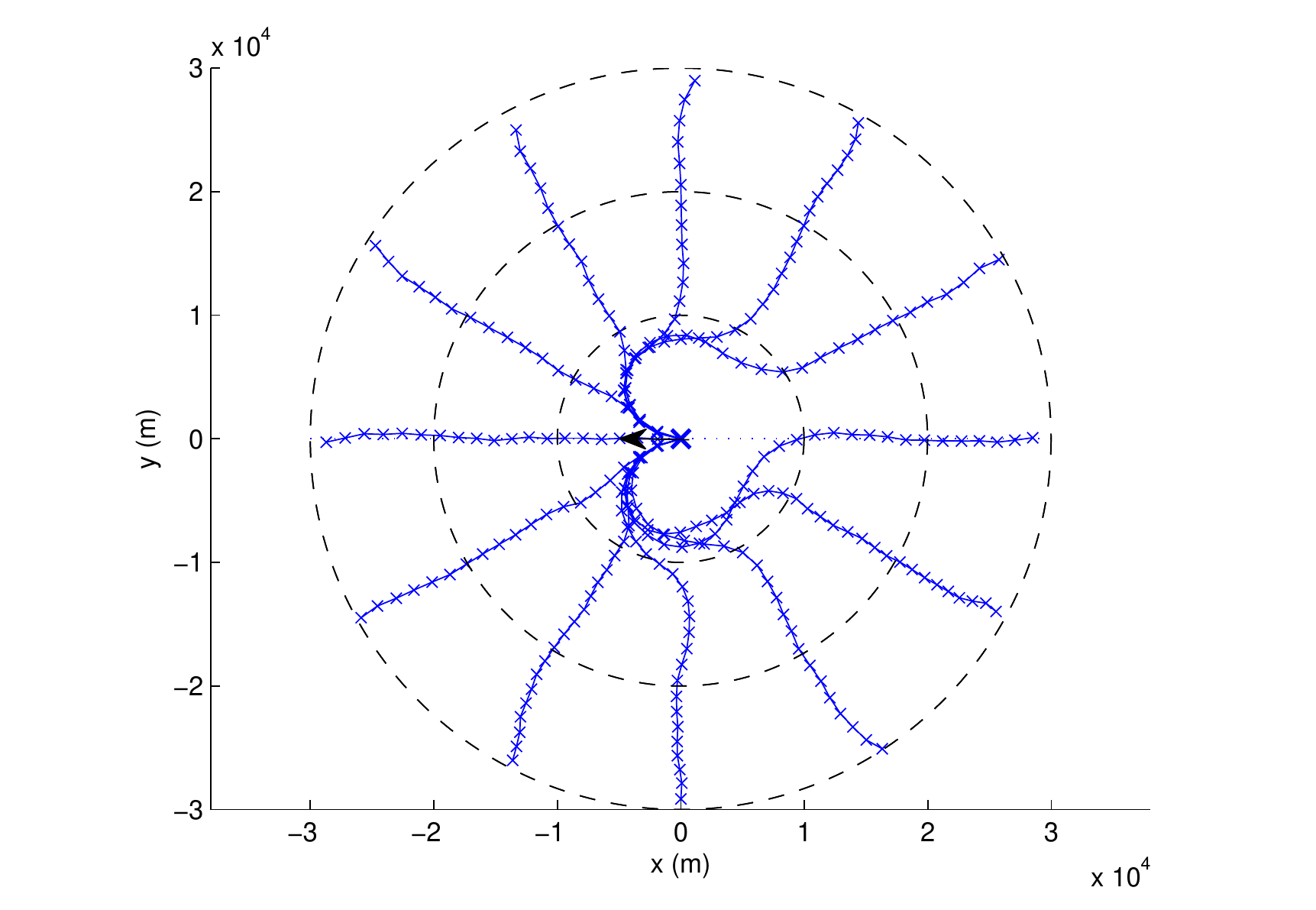}
\caption{Departure trajectories for 12 single aircraft taking off from a runway at the origin in a Westerly direction before heading to their desired bearing}
\label{fig:depart}
\end{center}
\end{figure}

\begin{figure*}[h]
\begin{center}
\def\picwidth{13cm}
\includegraphics[width=\picwidth]{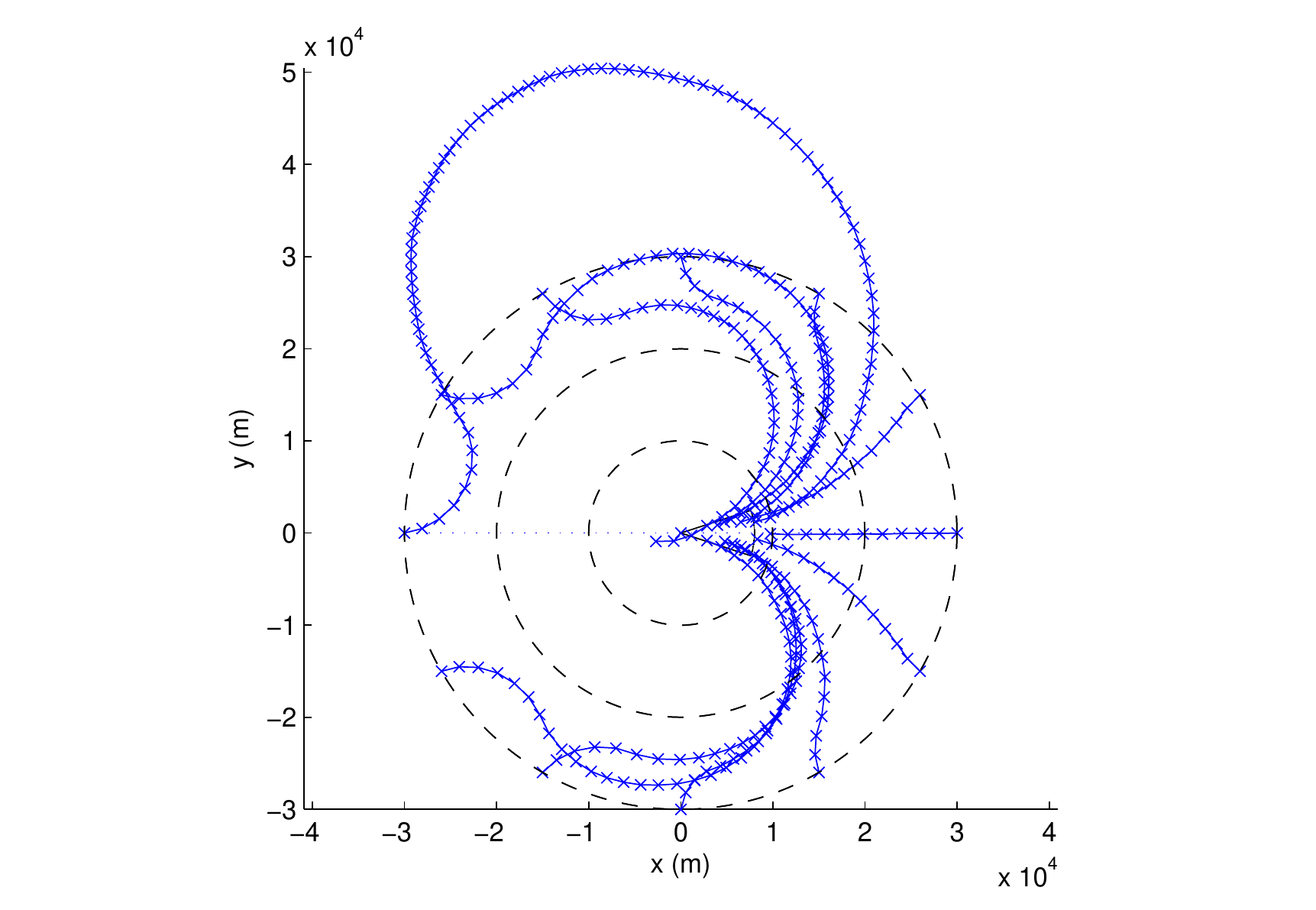}
\caption{Arrival trajectories for 12 single aircraft entering the TMA at different angles before proceeding to the landing envelope to land on an runway from the East}
\label{fig:arrive}
\end{center}
\end{figure*}

The presented results were used to tune the coefficients in the objective functions for both departure and arrival aircraft and evaluate the behaviour of arrival aircraft when exposed to the flow field. In each case for departures and arrivals 12 aircraft have been independently simulated and then plotted together on Figures~\ref{fig:depart} and~\ref{fig:arrive} to give an overview on behaviour. The arrival aircraft were started at different places around the airport and simulated until they reached the same landing zone. Departing aircraft were given different desired bearings and simulated until they reached the edge of the TMA zone.


The effect of the flow field on arrival aircraft is clear in Figure~\ref{fig:arrive}, on aircraft approaching from the Eastern side of the runway minimal diversion is taken prior to arrival. However the vehicle approaching from due west of the runway is given a long diversion outside the TMA before re-approaching from a more suitable bearing. In reality aircraft are allowed to approach an airport from a restricted number of directions and flying directly into departing aircraft traffic would be unlikely.

From this tuning the weightings shown in Table~\ref{tab:coeff} were used for landing and take off aircraft in the remainder of the scenarios in the report. Tuning these parameters differently will have significant effects on the trajectories generated.
\begin{table}[h]
  \centering
\begin{tabular}{ c | c | c | l }
  \hline                        
Aircraft Type & Coefficient & Value & Associated Cost \\
\hline
  \multirow{4}{*}{Departure}&$\alpha_1$ & 0.4 & Difference in bearing from desired to current\\
  &$\alpha_2$ & 0.1 & Fuel minimisation\\
  &$\alpha_3$ & 0.25 & Difference in altitude from desired to current\\
  &$\alpha_4$ & 0.25 & Difference in airspeed from desired to current\\
\hline
\multirow{4}{*}{Arrival} & $\tilde\alpha_1$ & 0.25 & Difference from aircraft heading to flowfield\\
 & $\tilde\alpha_2$ & 0.65 & Difference from aircraft altitude to nominal altitude\\
 & $\tilde\alpha_3$ & 0.1 & Fuel minimisation\\
  \hline  
\end{tabular}
\caption{Example objective function weightings used for optimisation}
  \label{tab:coeff}
\end{table}

\subsection{Standard Arrivals and Departures}
\begin{figure}[htb]
\begin{center}
\def\picwidth{12cm}
\includegraphics[width=\picwidth]{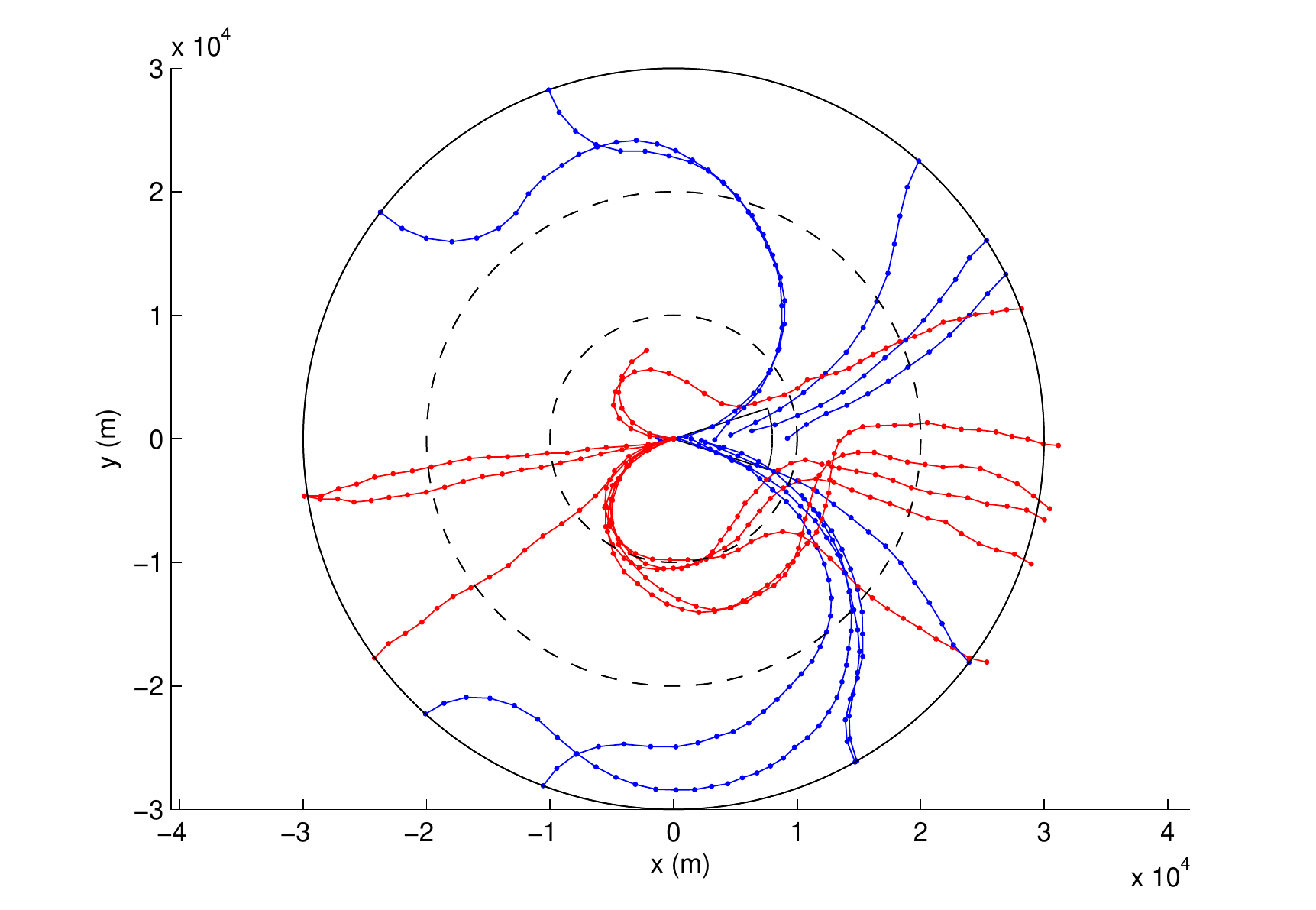}
\caption{Trajectories of 10 arrival (blue) and 10 departure (red) vehicles over a 100 time step simulation}
\label{fig:LandDepart}
\end{center}
\end{figure}

The simulations presented here demonstrate the unified running of the departure and arrival aircraft trajectories for the airport. Scenarios were randomly generated such that both arrival and departure aircraft joined the simulation at fixed times evenly spread across the entire simulation interval with a variance of a few time steps. The simulations were run for 100 time step updates (total scenario length of around 15 minutes). Between 5-10 vehicles of both arrival and departure types were included in each scenario. Figure~\ref{fig:LandDepart} shows an example of final trajectories for a 10 landing 10 departure scenario. Not all scenarios completed successfully. In cases where departing aircraft joined the problem when a landing aircraft was due to enter the landing envelope separation constraints were violated and no valid solution could be found. This is one of the restrictions of having fixed departure times. These could be amended by allowing flexibility in scheduling like a real airport. Whether this flexibility can be directly implemented in the SMC optimisation without an outside scheduling agent is a topic of future work.

\subsection{Low Fuel Scenario}
\begin{figure}[htb]
\begin{center}
\def\picwidth{12cm}
\includegraphics[width=\picwidth]{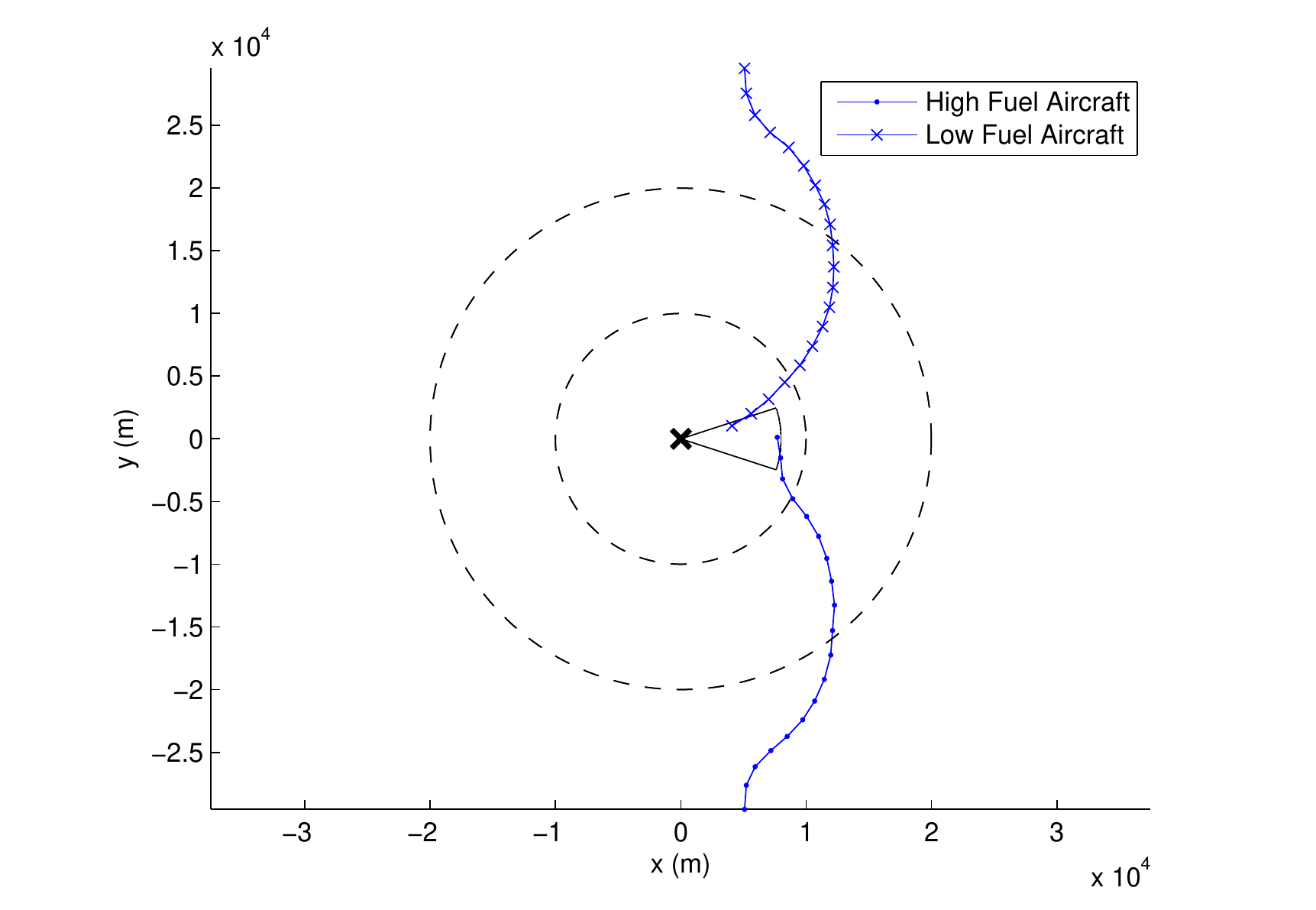}
\caption{Example trajectory for an aircraft with low fuel and an aircraft with high fuel approaching the landing envelope at the same time (low fuel landing first)}
\label{fig:lowfuel}
\end{center}
\end{figure}

In the low fuel scenario two arriving aircraft with the same distance till landing are initialised, one with plenty of fuel and the other with very little fuel remaining. This scenario is then run 20 times till both aircraft have landed to see how the aircraft behave. In each run of the scenario the aircraft are subject to a different set of disturbances from wind as well as the stochastic nature of the algorithm. Figure~\ref{fig:lowfuel} shows one of the generated trajectories. In the case of a human operator we would expect the aircraft with low fuel to be prioritised for landing with the aircraft with more fuel being diverted or held. However as Table~\ref{tab:fuel} shows the order in which the aircraft landed was actually more split. The cases where the low fuel aircraft landed second was due to the aircraft using its fuel sparingly for acceleration; most likely due to the fuel minimisation term in the objective function. Thus the high fuel aircraft proceeded faster towards the landing envelope and landed first. In both situations the low fuel aircraft lands before it runs out of fuel. To alter this behaviour an outside observer with some form of priorities would be needed to be added to the system to schedule the order which aircraft may enter the landing envelope.
\begin{table}[htb]
  \centering
\begin{tabular}{ c | c | c }
  \hline                        
Aircraft which landed first & Number of times & Average fuel (\kilogram) remaining on low fuel aircraft \\
\hline
High Fuel & 11 & 22.6 \\
\hline
 Low Fuel& 9 & 6.2 \\
  \hline  
\end{tabular}
\caption{Breakdown of which aircraft landed first in low fuel scenario and average fuel remaining }
  \label{tab:fuel}
\end{table}

\subsection{Congested Arrival Schedule}
\begin{figure}[htb]
\begin{center}
\def\picwidth{12cm}
\includegraphics[width=\picwidth]{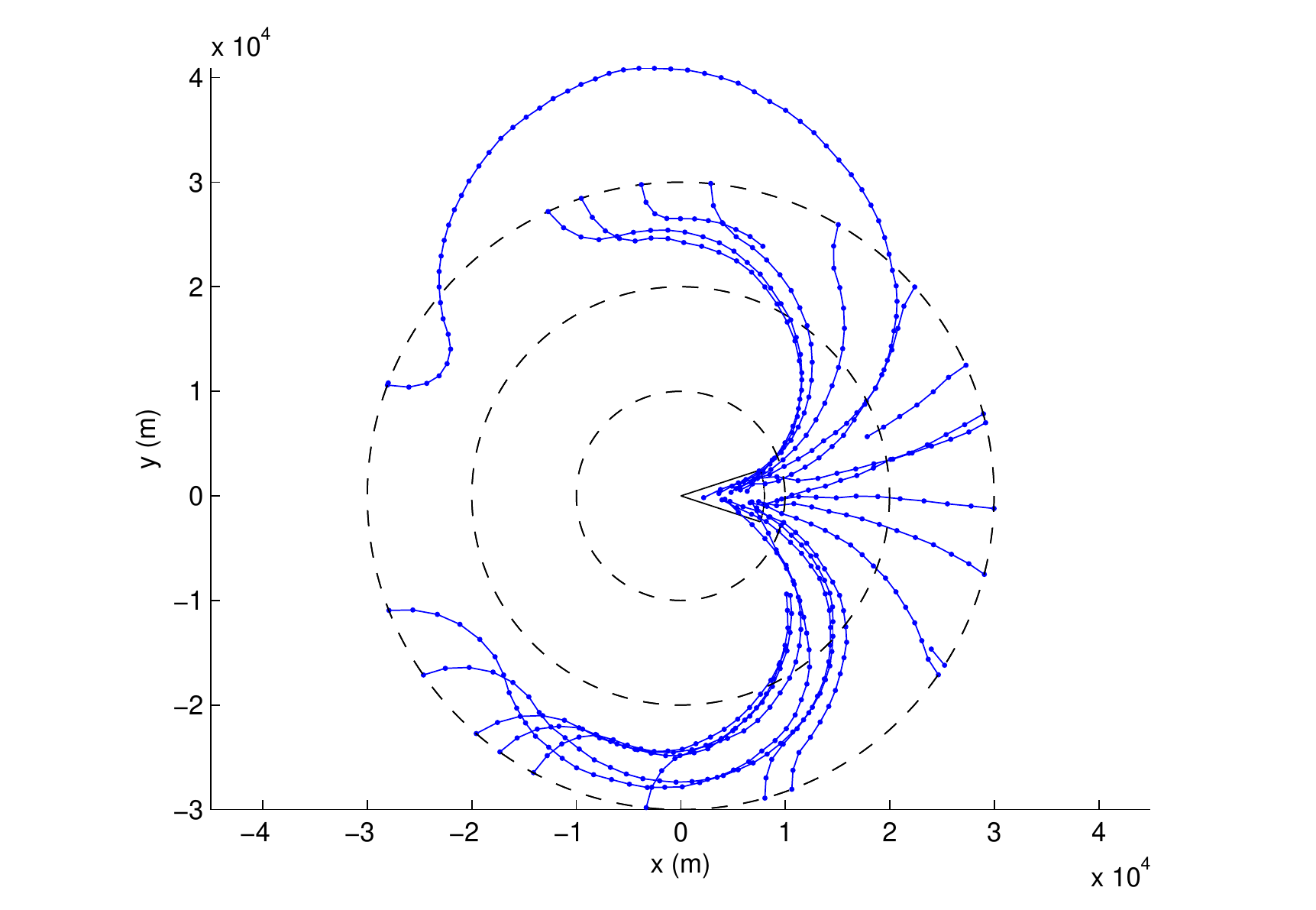}
\caption{Aircraft trajectories from a 100 step simulation with 24 arriving aircraft from random locations}
\label{fig:congest}
\end{center}
\end{figure}

The congested arrival schedule is a demonstration of the number of aircraft the system can handle concurrently. Previous work by the authors did demonstrate the method for up to 20 aircraft simultaneously in cruise-like situations. However with the limited resource of the single runway it is not necessarily sensible to initialise 20 aircraft in the first time step and run until all have landed. Instead aircraft are introduced to the scenario gradually as would happen in a real airport and removed from the scenario as they land. No departure aircraft are considered in this scenario as this would necessitate flexible scheduling of take offs which is currently neglected for simplicity (but as mentioned is the subject of future work). Figure~\ref{fig:congest} shows 24 arrival aircraft landing trajectories over time. The maximum concurrency observed in this scenario was 8 vehicles.

\section{Comparison to Real Airport Data}
\subsection{Collection of Real Data}

Data was collected from FlightRadar24. Using the open broadcast system automatic dependent surveillance-broadcast (ADS-B). ADS-B signals are picked up by a network of receivers which forward the information to the FlightRadar24 servers. The information that each aircraft broadcasts is a mixture of its identity, GPS position, heading, airspeed and timestamps. In Europe roughly 75\% of passenger aircraft contain ADS-B equipment. Relatively accurate information about aircraft can be obtained in populated areas (99\% of Europe is covered with ADS-B receivers) though coverage can be limited as aircraft fly further over oceans. As this work focuses on the terminal manoeuvring area of an airport such as Gatwick, the area around this airport does have good receiver coverage and thus aircraft with the equipment can be reliably tracked.

Data from the FlightRadar24 servers is stored in one minute update intervals which can be queried back as far as a month. The data used for comparison in this work was minute updates from the 14th of September 2013 over the full 24 hours, resulting in 1440 snapshots of aircraft across the globe. This data was then filtered to only include aircraft which were within 50km of London Gatwick airport and had Gatwick as either their destination or arrival airport.. Other traffic within the TMA not related the airport was neglected to avoid the need for the modelling of a third class of aircraft (aircraft which are neither arrival nor departure). Aircraft can continue to broadcast their ADS-B signals whilst on the ground at the airport so a minimum altitude was also imposed to allow detection of aircraft which were actually in flight.

The remaining data was then parsed into Matlab data files into three key data structures:
\begin{itemize}
\item The numerical state data across the entire 1440 time intervals containing, x, y co-ordinates translated from the latitudes and longitudes, the heading of the aircraft, the altitude of the aircraft and the true airspeed.
\item The non numerical data for each individual aircraft detected throughout the day which contains their flight code, aircraft type, registration, departure airport and arrival airport
\item The final contained over-arching numerical data about each aircraft including the first time step it was detected in the area, the last time step it was detected in the area, a binary value to indicate if it was an arrival or departure aircraft and 7 key aerodynamic coefficients and aircraft measures relating to the model of aircraft (obtained by cross referencing with BADA using the aircraft model).
\end{itemize}

The data obtained from Flightradar24 was relatively raw data in that there were data corruption issues or errors which needed to be filtered down before simulations of fuel usage could take place. Around 600 aircraft were registered as appearing in the 50km TMA of LGW and related to the airport, this was reduced to 556 once any aircraft which only appeared for fewer than 3 intervals were removed from the data set. These 556 aircraft were then iteratively reduced down to 528 as aircraft which vanished before leaving the TMA or landing or had clear errors in their GPS data were removed from the data set. One aircraft on the day did a go-around to have a second attempt at landing and as such was removed from the data set in order to provide a fair comparison of fuel use between our simulations and the real data was removed from the data set for that. Of the remaining aircraft 13  incorrectly identified whether they were arriving or departing. This was corrected in the data directly and the aircraft remained as part of the data set.

Weather reports for the 14th of September were obtained in order to provide estimates of the nominal wind directions throughout the day. These initial estimates of wind directions and strengths were updated at 20-30 minute intervals which were then interpolated for a minute by minute estimate. In cases where no clear wind direction was recorded on the weather data a random wind bearing was chosen. There were no unusual or remarkable weather conditions on the day in question, the prevailing wind was West to East so no changes in runway direction were recorded over the day.

Since the data obtained from FlightRadar24 contains no estimate of the aircraft's weight the several assumptions are made. All arrival aircraft are operating at their max take-off load and have 20\% of their fuel reserve left when they enter the TMA. All departure aircraft are operating at their max take-off load and have their full fuel reserve at take-off. Maximum loads were obtained from each model of aircraft's specs.

\subsection{Analysis of Real Data}

\begin{figure}[htb]
\begin{center}
\def\picwidth{12cm}
\includegraphics[width=\picwidth]{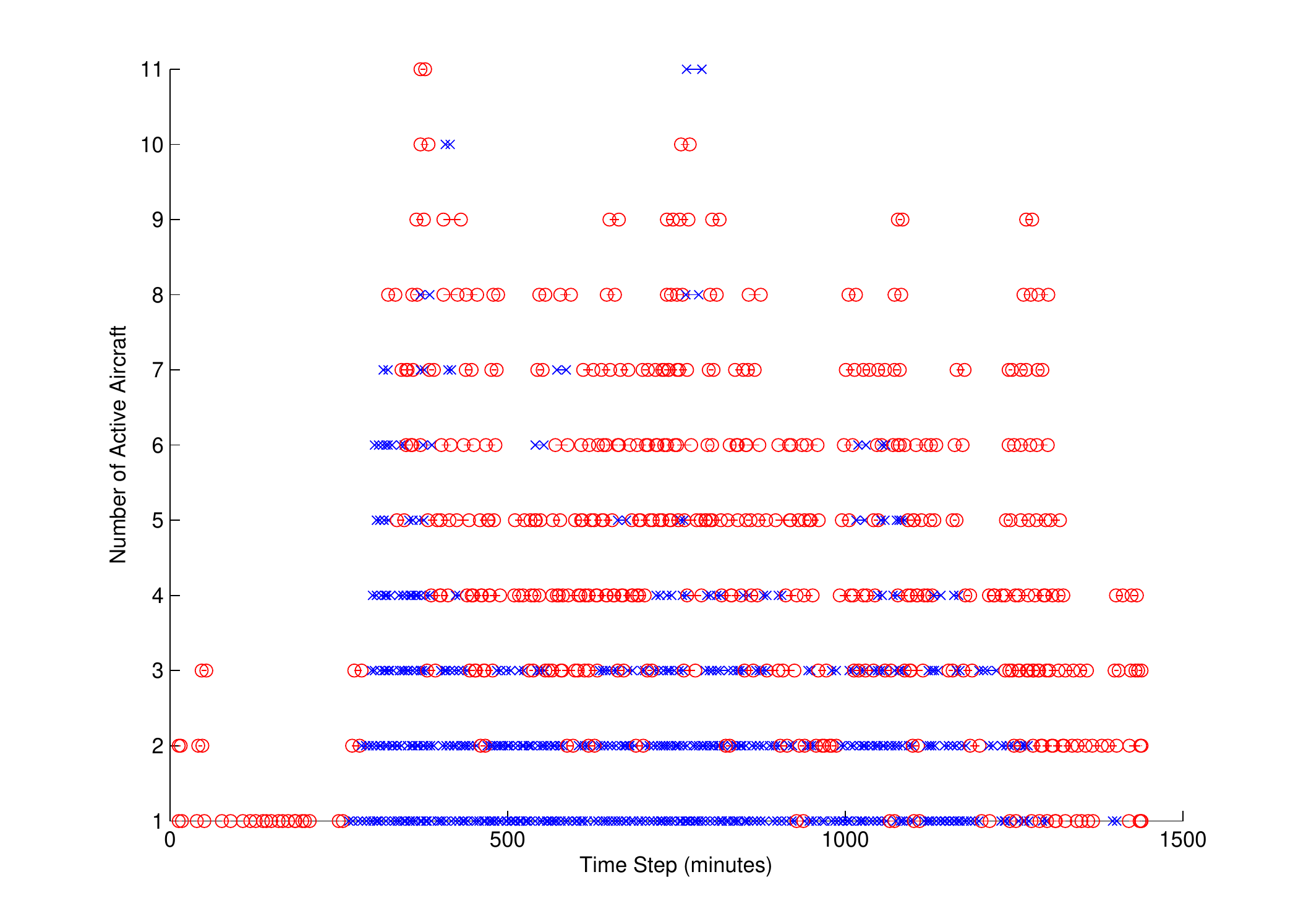}
\caption{Pattern of 528 real aircraft (arrival red circles, departure blue cross) active across 24 hours at Gatwick airport}
\label{fig:flight_concurrency}
\end{center}
\end{figure}

With the remaining 528 vehicles' data available we processed the information to get an overview of the airport operations prior to fuel use calculations. Figure\ref{fig:flight_concurrency} shows a concurrency plot of aircraft in the TMA 40km radius of the airport (although aircraft were recorded up to 50km away they were only considered to enter the TMA at 40km). The beginning and end points for each aircraft are marked by either a circle (for arrival aircraft) or cross (for departure aircraft) connected by a line between the two markers. The y axis counts how many aircraft are concurrently active in the TMA and the x axis displays the time in minutes across the day (1440 minutes = 24 hours). Departures from the airport don't begin until around 4am at which point activity at the airport increases rapidly and stays high across the daytime hours until it begins to curtail later into the evening. Arrivals can occur at all times of the day and night. This plot allows us to identify key times to simulate the airport for various traffic demands over the 24 hours.
\begin{figure}[htbp]
\begin{center}
\def\picwidth{11cm}
\includegraphics[width=\picwidth]{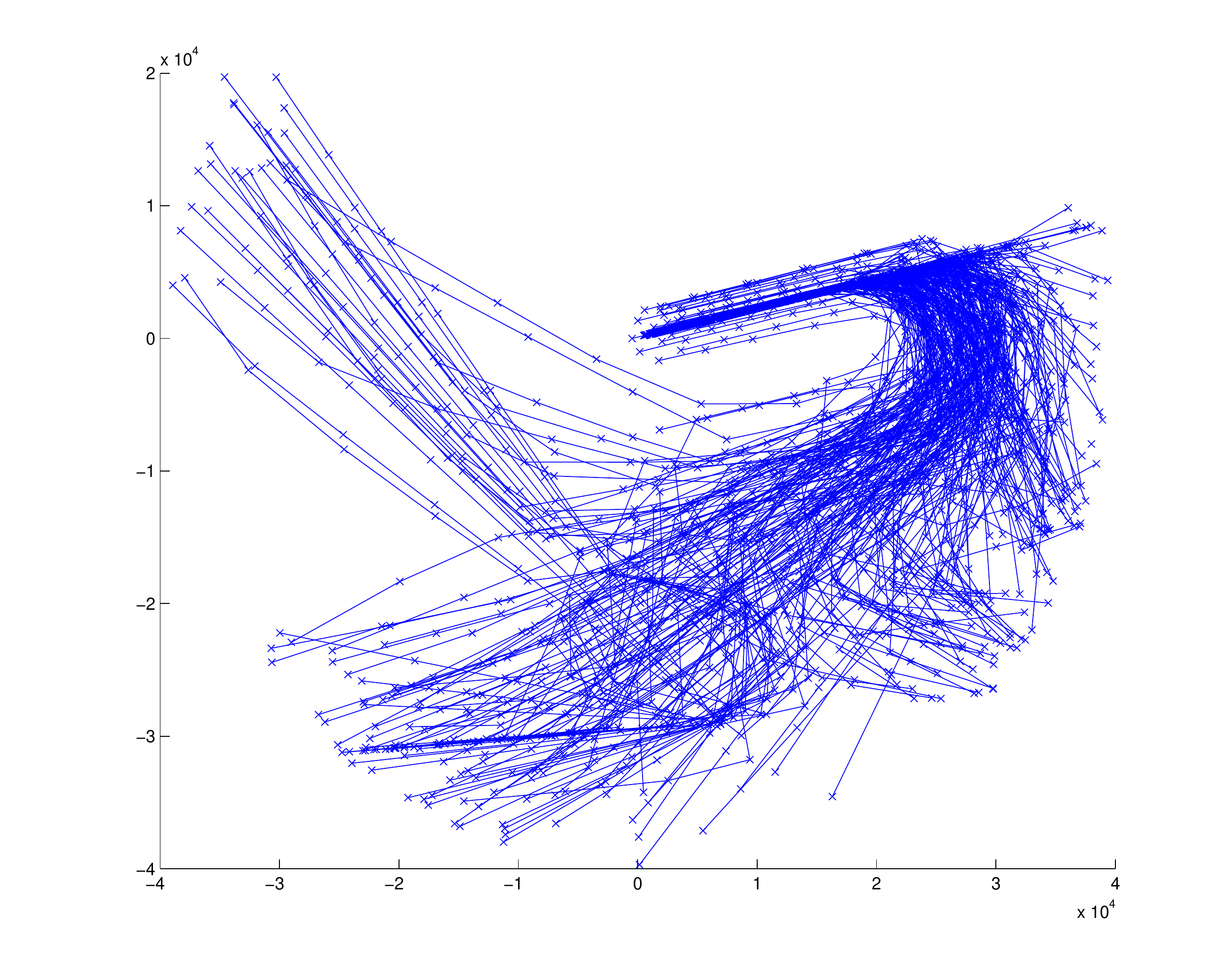}
\caption{Paths of all arrival aircraft within the 24 hour period at Gatwick airport}
\label{fig:gat_arrival}
\end{center}
\end{figure}

\begin{figure}[htbp]
\begin{center}
\def\picwidth{11cm}
\includegraphics[width=\picwidth]{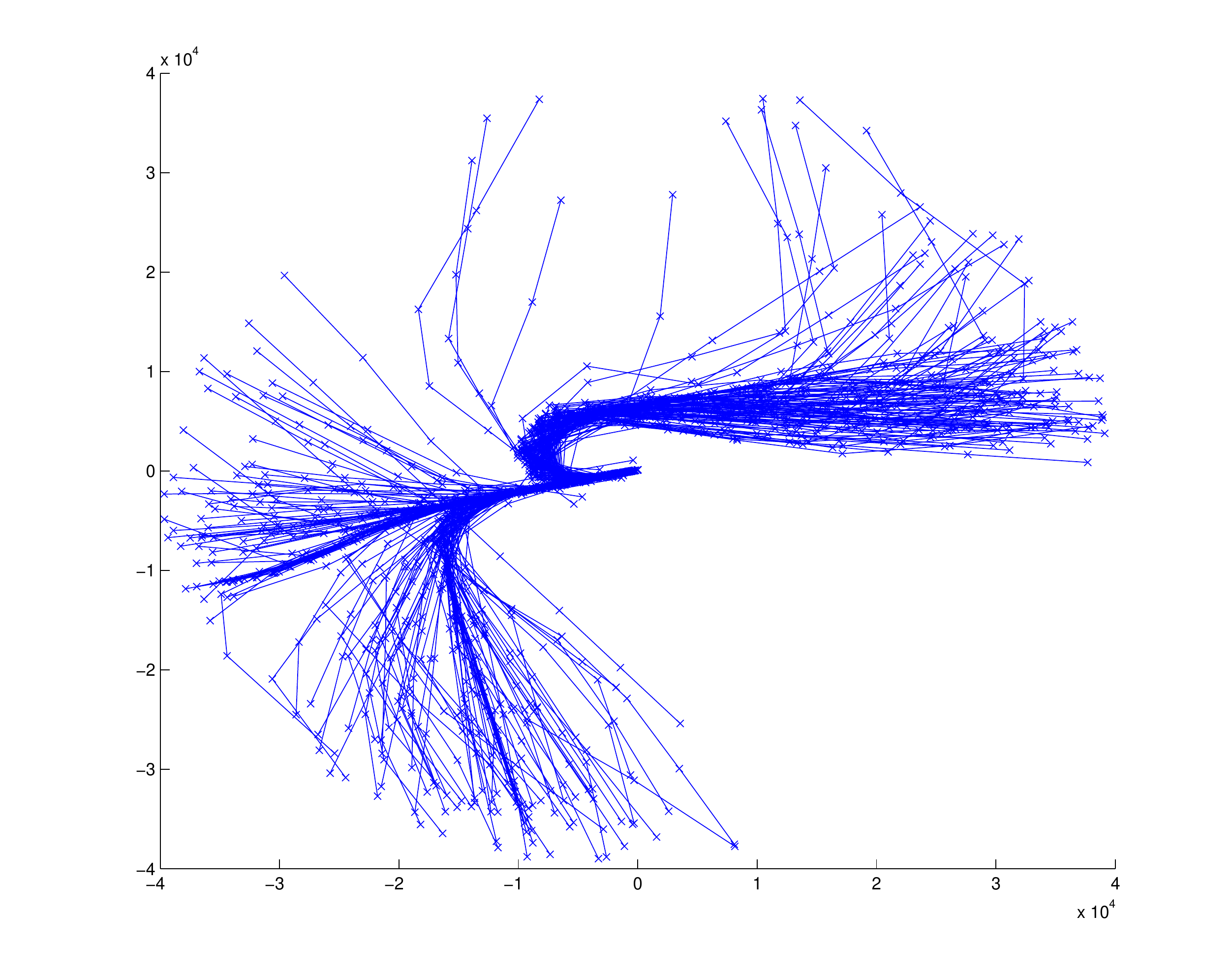}
\caption{Paths of all departure aircraft within the 24 hour period at Gatwick airport}
\label{fig:gat_departure}
\end{center}
\end{figure}
Figures\ref{fig:gat_arrival} and \ref{fig:gat_departure} display the paths taken by all arrival and all departure aircraft respectively across the day. These paths give an indication of the way aircraft are managed at the airport in defined streams to avoid arrival and departure aircraft interacting. 

\subsubsection{Fuel Estimate 1}
\begin{figure}[hptb]
\begin{center}
\def\picwidth{12cm}
\includegraphics[width=\picwidth]{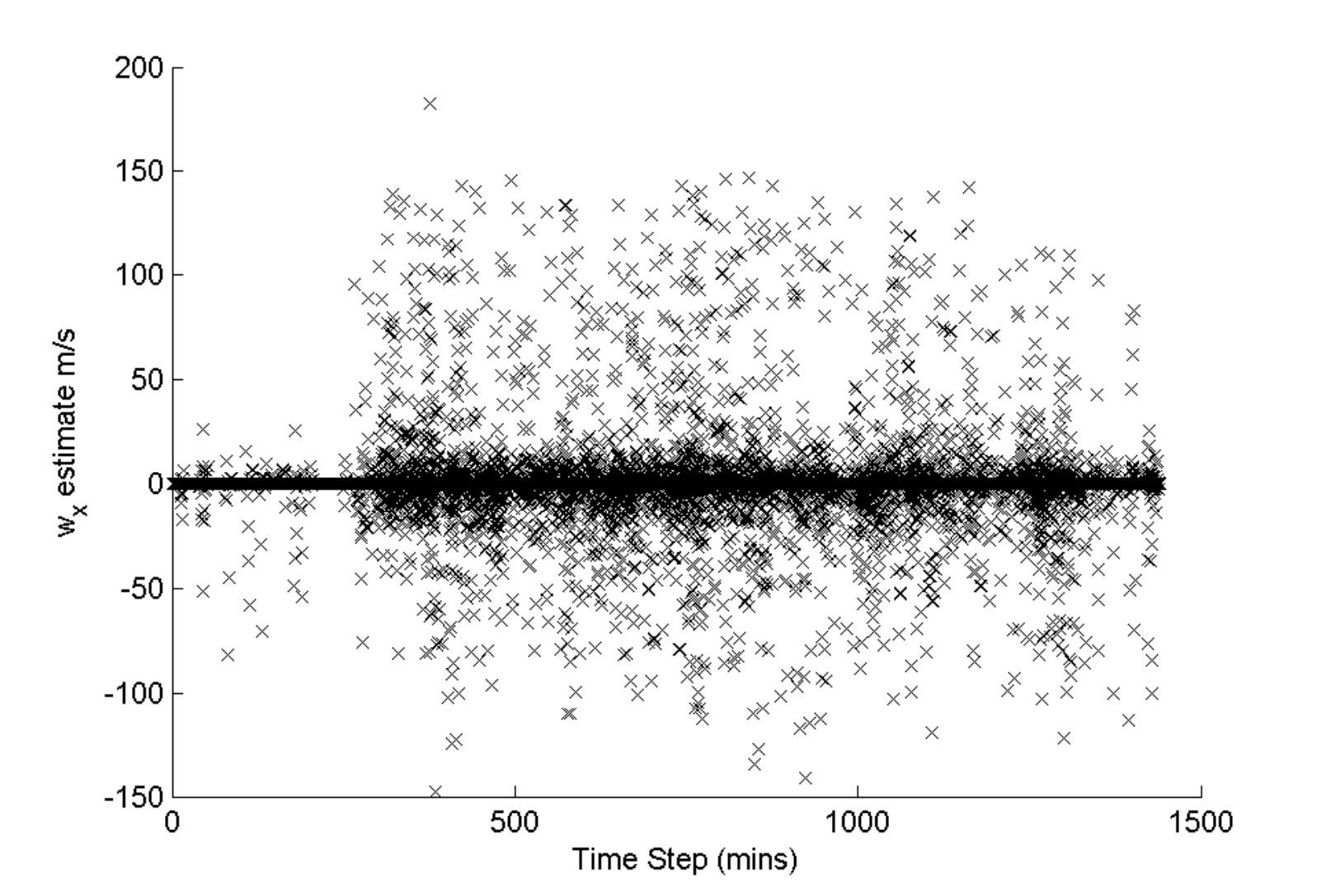}
\caption{Estimates of $w_x$ over the 24 hour data period}
\label{fig:w_x}
\end{center}
\end{figure}

\begin{figure}[hptb]
\begin{center}
\def\picwidth{12cm}
\includegraphics[width=\picwidth]{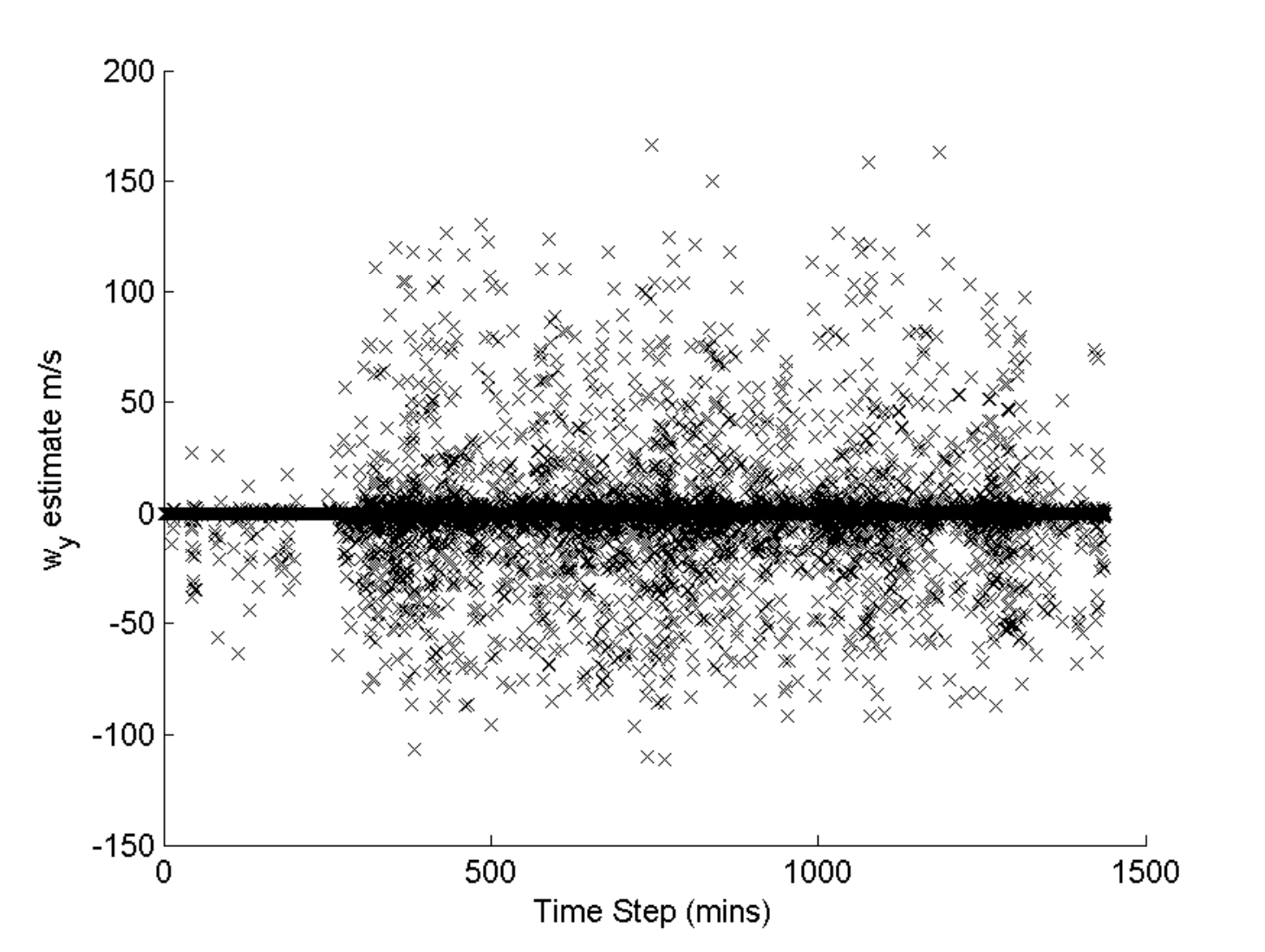}
\caption{Estimates of $w_y$ over the 24 hour data period}
\label{fig:w_y}
\end{center}
\end{figure}

Fuel estimate 1 assumes that aircraft fly on the heading that they are recorded at for each minute interval, accelerating or decelerating such that their airspeed will match the next interval's value. Any difference between the new x,y position and the second interval's x,y position is considered to be caused by disturbance or 'wind'. The aircraft dynamics model is effectively run in reverse to calculate the change in mass of the aircraft and thus the fuel burned. 
\begin{subequations}
\begin{align}
\eta_i(k)=&C_{f1,i}\left(1+{{v_{s,i}(k)}\over{C_{f2,i}}}\right) \\
\gamma_i(k)=& \sin^{-1}\left({{z_i(k+1)-z_i(k)}\over{\delta t v_{s,i}(k)}}\right) \\
w_{x,i}(k) =&{{x_i(k+1)-x_i(k)}\over{\delta t}}- (v_{s,i}(k) \cos(\chi_i(k))\cos(\gamma_i(k))) \\
w_{y,i}(k) =&{{y_i(k+1)-y_i(k)}\over{\delta t}}- (v_{s,i}(k) \sin(\chi_i(k))\cos(\gamma_i(k))) \\
T_i(k+1)=& {{m_i(k)(v_{s,i}(k+1)-v_{s,i}(k))}\over{\delta t}}+D_i(k) +m_i(k)g\sin(\gamma_i(k)) \\
m_i(k+1)=&m_i(k)-\delta t(\eta_i(k) T_i(k))
\end{align}
\label{eqn:blah}
\end{subequations}
Where $\eta_i(k)$ is the fuel burn rate of aircraft $i$ at step $k$, based on the true airspeed and two aircraft model specific aerodynamic coefficients $C_{f1,i},C_{f2,i}$ obtained from BADA.
The following values are known from the data: $x_i(k), y_i(k), z_i(k), v_{s,i}(k), \chi_i(k), m_i(1), x_i(k+1), y_i(k+1), z_i(k+1), v_{s,i}(k+1), \chi_i(k+1)$ with $m_i(k+1)$ needing calculation. For this method of fuel estimation additional parameters $w_x$ and $w_y$ are obtained. These paramters act as an estimate of how accurate the fuel usage might be. When $w_x$ and $w_y$ are in the same order of magnitude as the nominal wind conditions at that time then the fuel estimate is likely reasonable. Figures \ref{fig:w_x} and \ref{fig:w_y} show all the estimates of $w_x$ and $w_y$ obtained at each time step for every aircraft. The vast majority of $w_x$ and $w_y$ estimates are around the origin which indicates good general performance, however there are also many cases where wind estimates are clearly too high for the fuel estimates associated with them to be entirely accurate. This motivates the need for a secondary method of fuel estimation which can be used for comparison both with this first estimate and with data obtained from our own fuel simulations.

\subsubsection{Fuel Estimate 2}
To provide a second estimate of fuel use from the real data an alternative method was used. This other method uses dead reckoning of total distance travelled between the two points $d_i$ and assumes no wind disturbance was present. This distance is used to get an estimate for true airspeed of the aircraft $\hat{v}_{i,s}$. From the airspeed estimate, the fuel burn coefficient $\eta_i$ is obtained along with the climb angle $\gamma_i$. An estimate of track angle $\hat{\chi}$ is obtained from the two positions of the aircraft, which when compared to the aircraft's reported heading, leads to the bank angle needed to change the heading by that margin. Finally these values can be used to work out the thrust needed to get the velocity change between the airspeed reported by the aircraft at the next time step and the estimate of true airspeed.
\begin{subequations}
\begin{align}
d_i(k)=& \sqrt{(x_i(k+1)-x_i(k))^2+(y_i(k+1)-y_i(k))^2+(z_i(k+1)-z_i(k))^2}\\
\hat{v}_{s,i}(k)=&{{d_i(k)}\over{\delta t}}\\
\eta_i(k)=&C_{f1,i}\left(1+{{\hat{v}_{s,i}(k)}\over{C_{f2,i}}}\right) \\
\gamma_i(k)=& \sin^{-1}\left({{z_i(k+1)-z_i(k)}\over{\delta t \hat{v}_{s,i}(k)}}\right) \\
\hat{\chi}_i(k)=&\tan^{-1}\left({{y_i(k+1)-y_i(k)}\over{x_i(k+1)-x_i(k)}}\right)\\
\delta \chi_i(k)=&\hat{\chi}_i(k)-\chi_i(k)\\
\phi_i(k)=&\tan^{-1}\left({{\delta \chi_i(k) \hat{v}_{i,s}(k)}\over{g \delta t}}\right)\\
T_i(k+1)=& {{m_i(k)(v_{s,i}(k+1)-\hat{v}_{s,i}(k))}\over{\delta t}}+D_i(k) +m_i(k)g\sin(\gamma_i(k)) \\
m_i(k+1)=&m_i(k)-\delta t(\eta_i(k) T_i(k))
\label{eqn:blah}
\end{align}
\end{subequations}

This method can suffer from issues where the x-y distance change is actually very small since the aircraft was circling and the path was poorly captured by the 1 minute sample time. This leads to poor estimates of airspeed which can result in large thrusts needed to reach the airspeed at the next interval and thus artificially high fuel usage estimates in some cases.

\subsection{Comparison Between Fuel Estimates}
\begin{figure}[hptb]
\begin{center}
\def\picwidth{12cm}
\includegraphics[width=\picwidth]{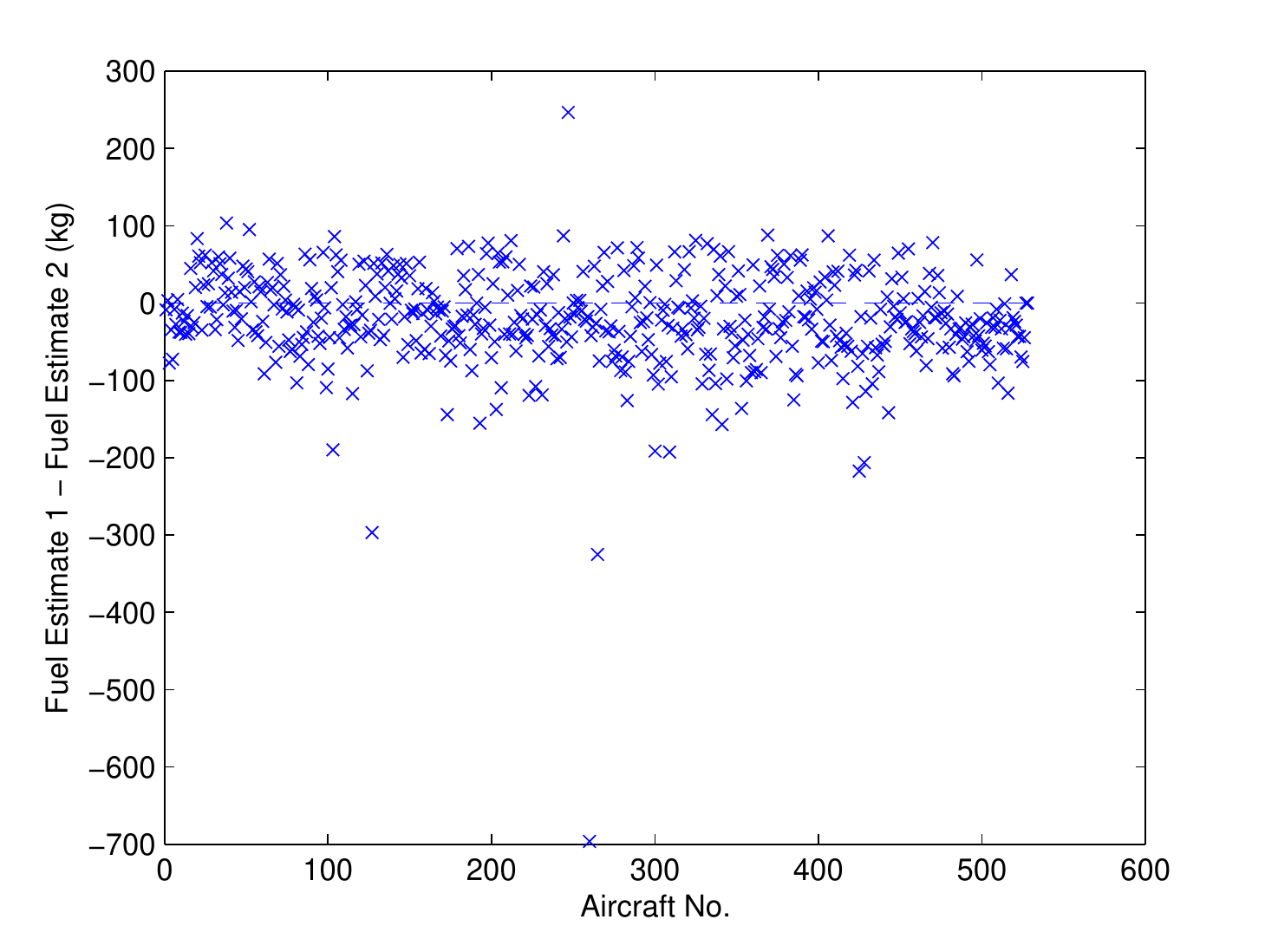}
\caption{Fuel Estimate 1 - Fuel Estimate 2 for each aircraft}
\label{fig:fuel_subtract}
\end{center}
\end{figure}

\begin{figure}[hptb]
\begin{center}
\def\picwidth{12cm}
\includegraphics[width=\picwidth]{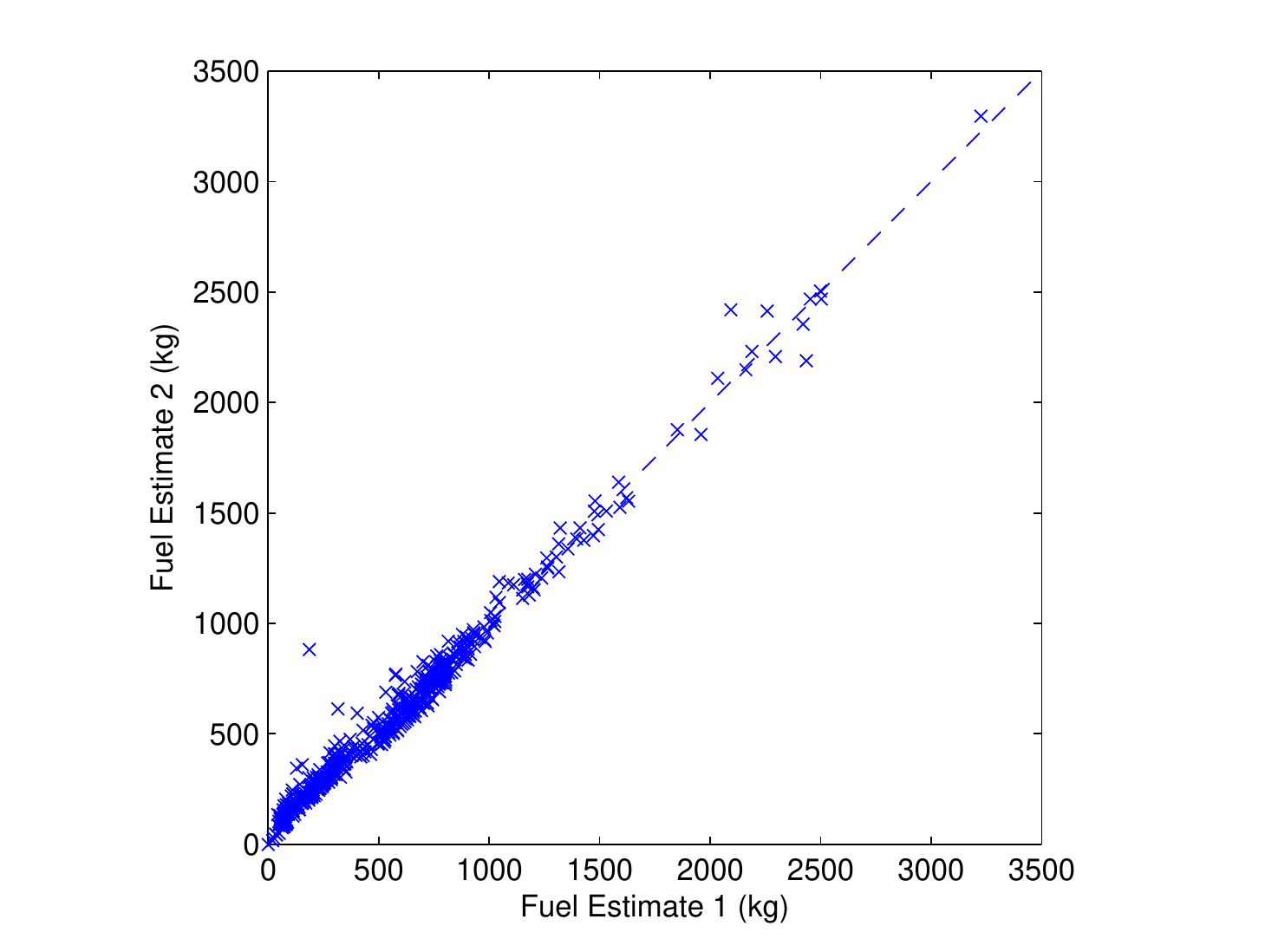}
\caption{Fuel Estimates 1 and 2 plotted against each other}
\label{fig:fuel_comp}
\end{center}
\end{figure}

For both fuel estimates aircraft were restricted from burning negative fuel and gaining weight. The two estimates were used as guidelines for how much fuel the aircraft might have used, though neither method can reliably claim to always upper or lower bound the fuel use. Figures \ref{fig:fuel_subtract} and \ref{fig:fuel_comp} show the difference between fuel estimate 1 and fuel estimate 2 for each of the 528 aircraft. The vast majority of cases are within 100kg of fuel each. The significant outliers are typically caused by a very poor fuel estimate from fuel estimate 2. This poor fuel estimate is due to a very short Euclidean distance apparently being travelled by the aircraft in the sampling time of a minute. In reality the outlier aircraft were still circling in stack with close to one revolution of the stack per minute. A short Euclidean distance leads to a severe under-estimate of $\hat{v}_{i,s}(k)$ which then leads to a significant over-estimate of thrust needed to increase airspeed to match the observed airspeed $v_{i,s}(k+1)$.

\subsection{Comparison Between Algorithm and Real Data}
Running a single simulation from our algorithm to simulate across a full 24 hour period of airport activity would be ideal but would cause issues due to the raw size of the problem, memory issues and stability issues. Several adjustments will therefore need to be made when comparing to the real data. The most important of these is to shorten the sampling time from every 60 seconds from the recorded data to time steps of 20 seconds. As highlighted in the various fuel estimates with a linear model between time steps longer sampling times do have negative effects on the accuracy of aircraft trajectories. Additionally although we only have recorded data from the aircraft every minute they are running a sophisticated inner loop controller in the form of a pilot. A pilot's decisions are not limited to minute by minute sampling times. The shorter sampling time of 20 seconds was chosen as a compromise between the additional computation load and the ability of the model to cope. Extending our sampling time would ideally require some form of inner loop controller and a distinct change of structure of our simulations.

Given the real world data available at every minute interval it is possible to initialise our own simulations starting at any particular minute of the day to run for a set number of steps. Computational speed is highly impacted by the number of aircraft within the simulation window. In high density traffic simulations taking place at the peaks of the day the length of simulation window was limited to one hour (180 simulation time steps). In low density operations (e.g. very early in the day) simulation window length was increased to take advantage of the reduced number of aircraft.

Fuel comparison simulations were classified into two types based on traffic density ranging from low (under 5 aircraft active per time step), to high density. These traffic densities were identified from the concurrency analysis shown on Figure \ref{fig:flight_concurrency}.

\subsubsection{Low Density Traffic}
Low density traffic mainly occurs in the very early hours of the morning between midnight and 4.30 am (data time steps 0 to 230). Traffic in low density was periods was dominated by arrivals as take offs do not begin in the airport until close to dawn in order to reduce the impact of noise pollution on residents. A similar period of low density arrival traffic takes place in the 100 minutes leading up to midnight. In this period the majority of aircraft are arrivals however 1 delayed departure does take place as well. Examining the low traffic density periods is a useful technique for fine tuning arrival aircraft objective function weightings and bug testing the simulations.

\begin{figure}[p]
\begin{center}
\def\picwidth{12cm}
\includegraphics[width=\picwidth]{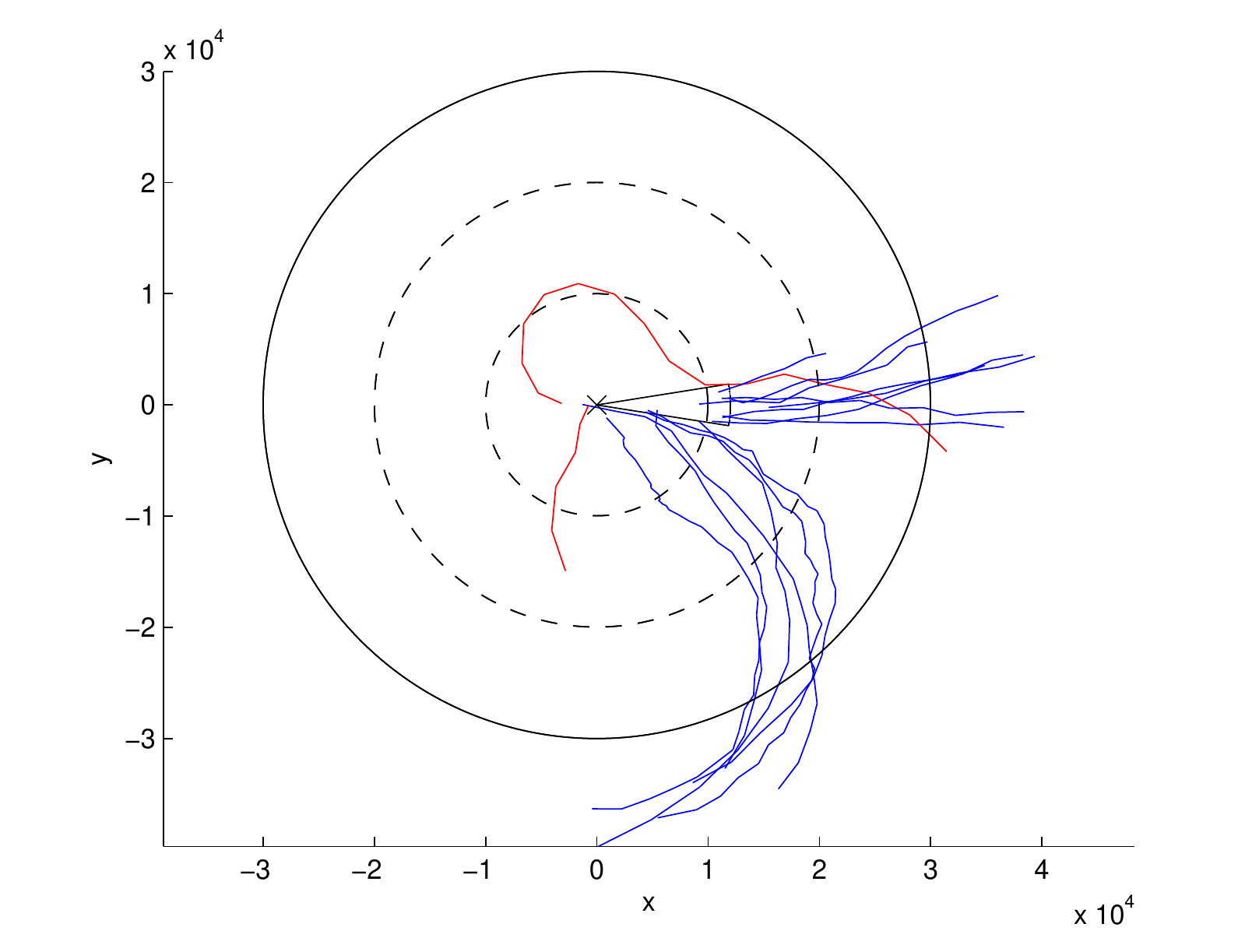}
\caption{Simulated trajectories of aircraft in first 230 minutes of Sept 14th}
\label{fig:low_1_traj}
\end{center}
\end{figure}
\begin{figure}[p]
\begin{center}
\def\picwidth{12cm}
\includegraphics[width=\picwidth]{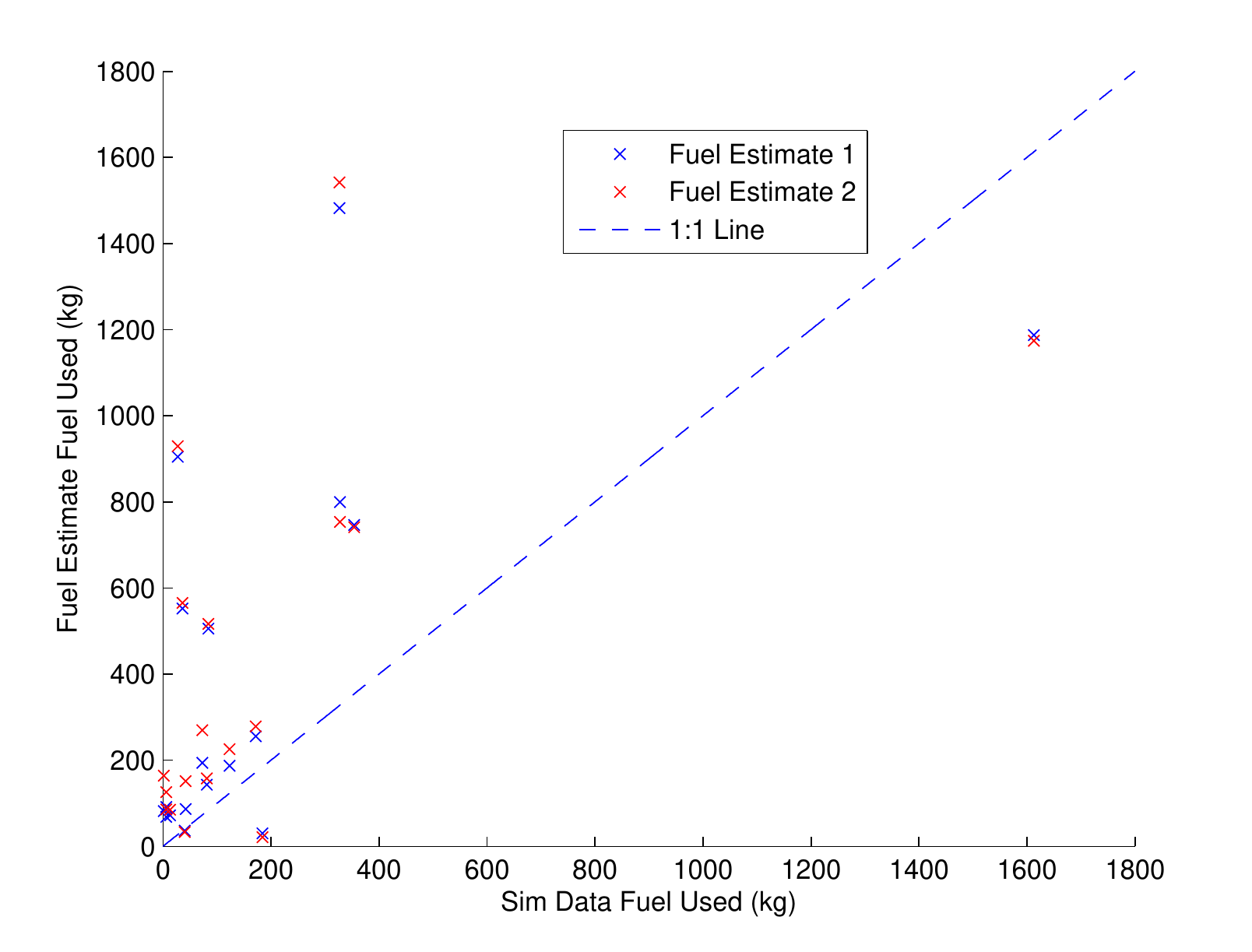}
\caption{Fuel comparison between fuel estimates and simulations for first 230 minutes of Sept 14th}
\label{fig:low_1_fuel}
\end{center}
\end{figure}

\begin{figure}[p]
\begin{center}
\def\picwidth{12cm}
\includegraphics[width=\picwidth]{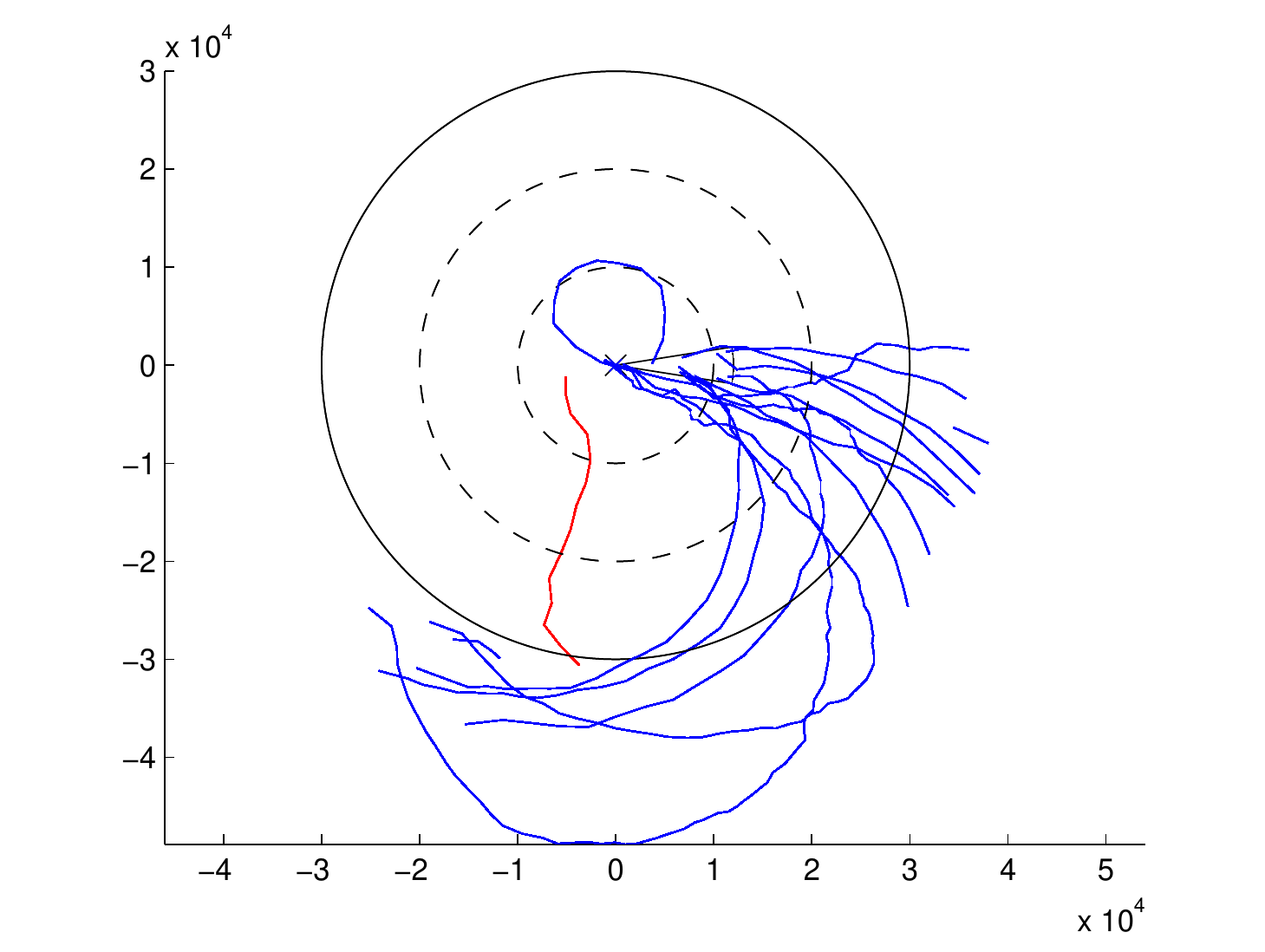}
\caption{Simulated trajectories of aircraft in last 100 minutes of Sept 14th (Blue = arrival aircraft, Red = departure aircraft)}
\label{fig:low_2_traj}
\end{center}
\end{figure}
\begin{figure}[p]
\begin{center}
\def\picwidth{12cm}
\includegraphics[width=\picwidth]{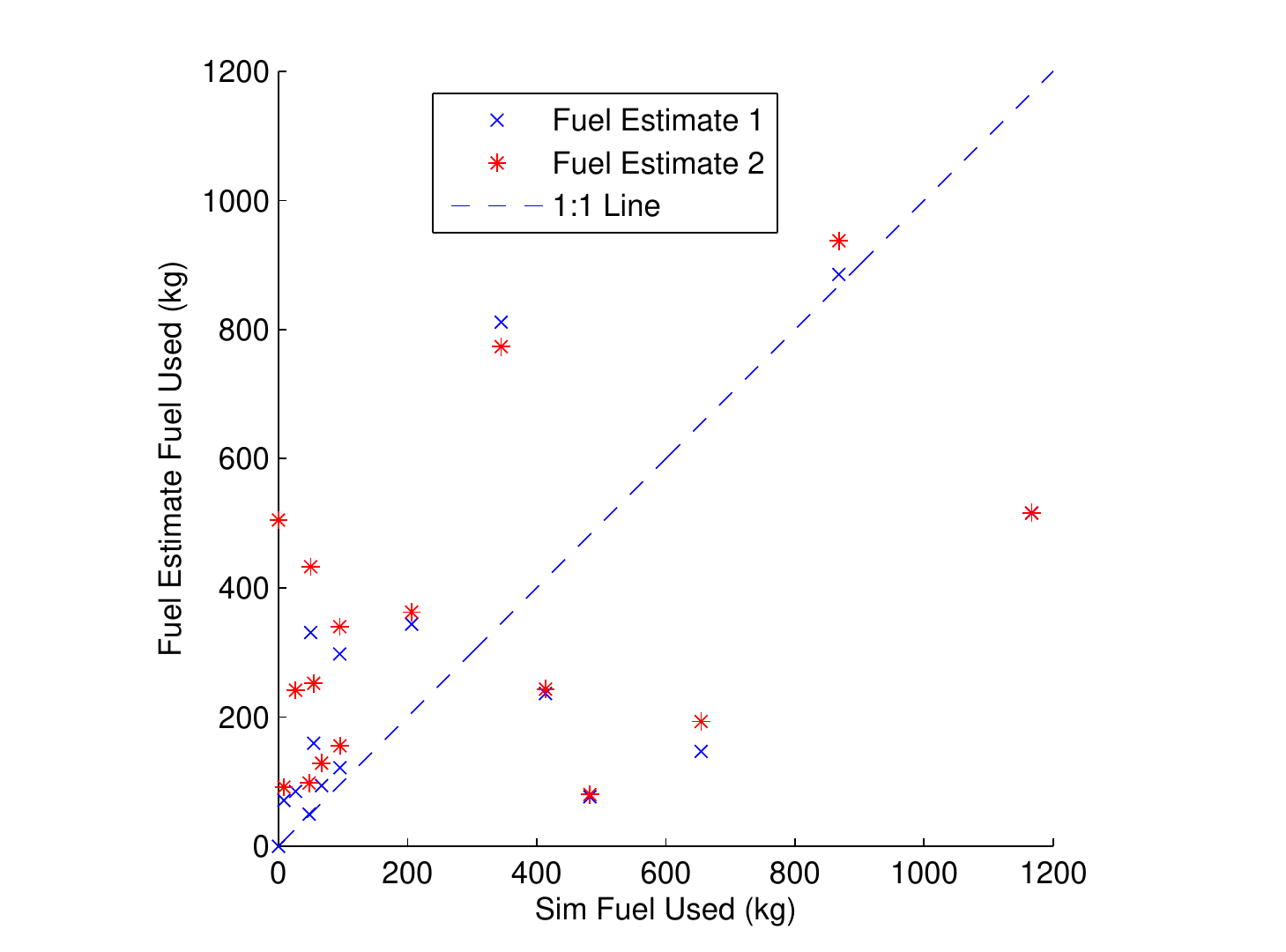}
\caption{Fuel comparison between fuel estimates and simulations for last 100 minutes of Sept 14th}
\label{fig:low_2_fuel}
\end{center}
\end{figure}

The final trajectories of the low density aircraft traffic are shown in Figures \ref{fig:low_1_traj} and \ref{fig:low_2_traj}. In one case of the last 100 minutes of Sept 14th a single arrival aircraft was unable to land successfully in the landing envelope on its first pass needing to do a go around before landing successfully the second time. The fuel cost of a go-around will be significant and should show up as an outlier in the fuel comparison between the real data and the simulations.

Figures \ref{fig:low_1_fuel} and \ref{fig:low_2_fuel} and Tables \ref{tab:low_1_fuel} and \ref{tab:low_2_fuel} show the fuel comparisons between the fuel estimates obtained from the data and the simulated fuel use. In both cases the majority of aircraft do make fuel savings from the simulation compared to the fuel estimates. In all tables $F_s$ is the fuel used in the simulation, $F_1$ is the first fuel use estimate, $F_2$ is the second fuel estimate and the mean $F$ is the mean of the two fuel estimates. 

\begin{table}[htb]
  \centering
\begin{tabular}{ c | c | c | c | c| c}
  \hline                        
Aircraft No. & $F_S$ (kg) & $F_1$ (kg) & $F_2$ (kg) & Mean $F$ (kg) & Fuel Saving \% \\
\hline
1 &40.23	&36.47	&33.25	&34.86 & -15.40\\
2 &184.20	&30.08	&21.02	&25.55 & -620.97\\
3 &73.02	&193.63	&269.61	&231.62 & 68.48\\
4 &13.31	&72.11	&84.99	&78.55 &83.05\\
5 &5.91	&91.17	&125.85	&108.52 &94.55\\
6 &83.87	&505.48	&516.39	&510.93 &83.58\\
7 &123.20	&187.12	&225.91	&206.51 &40.34\\
8 &353.95	&746.28	&741.06	&743.67 &52.41\\
9 &41.77	&86.49	&151.63	&119.06 &64.92\\
10 &5.90	&68.14	&85.73	&76.94 &92.34\\
11 &36.21	&552.18	&565.02	&558.60 &93.52\\
12 &171.88	&255.35	&278.81	&267.08 &35.64\\
13 &1.53	&81.74	&163.87	&122.81 &98.76\\
14 &80.95	&142.83	&157.48	&150.15 &46.09\\
15 &1612.85	&1187.06	&1173.86	&1180.46 &-36.63\\
16 &326.73	&1482.06	&1541.67	&1511.87 &78.39\\
17 &27.35	&904.74	&929.20	&916.97 &97.02\\
18 &327.89	&799.50	&753.69	&776.59 &57.78\\
\hline
Total &3510.7 &7422.4 &7819.0 & 7620.7 &53.93\\
\hline  
\end{tabular}
\caption{Fuel use comparison between simulation and estimates for first 230 minutes of Sept 14th}
  \label{tab:low_1_fuel}
\end{table}

\begin{table}[htb]
  \centering
\begin{tabular}{ c | c | c | c | c| c}
  \hline                        
Aircraft No. & $F_S$ (kg) & $F_1$ (kg) & $F_2$ (kg) & Mean $F$ (kg) & Fuel Saving \% \\
\hline
1 &67.15	&93.90	&128.59	&111.24	&39.64 \\
2 &1166.72	&515.56	&516.18	&515.87	&-126.16\\
3 &54.88	&159.55	&252.38	&205.96	&73.36\\
4 &654.54	&146.76	&192.92	&169.84	&-285.38\\
5 &868.04	&885.26	&937.35	&911.31	&4.75\\
6 &206.54	&343.97	&362.41	&353.19	&41.52\\
7 &345.05	&811.48	&773.74	&792.61	&56.47\\
8 &482.51	&76.57	&80.11	&78.34	&-515.93\\
9 &413.82	&236.52	&243.02	&239.77	&-72.59\\
10 &50.04	&330.67	&432.73	&381.70	&86.89\\
11 &95.04	&297.55	&340.01	&318.78	&70.19\\
12 &26.60	&84.90	&241.64	&163.27	&83.71\\
13 &95.42	&121.57	&155.55	&138.56	&31.13\\
14 &8.34	&70.90	&91.49	&81.20	&89.73\\
15 &47.48	&49.43	&97.68	&73.55	&35.45\\
16 &0.079	&0		&504.90	&252.45	&99.97\\
\hline
Total &4582.2    &4224.6   & 5350.7    &4787.7 & 4.29\\
  \hline  
\end{tabular}
\caption{Fuel use comparison between simulation and estimates for last 100 minutes of Sept 14th}
  \label{tab:low_2_fuel}
\end{table}

\begin{figure}[hptb]
\begin{center}
\def\picwidth{12cm}
\includegraphics[width=\picwidth]{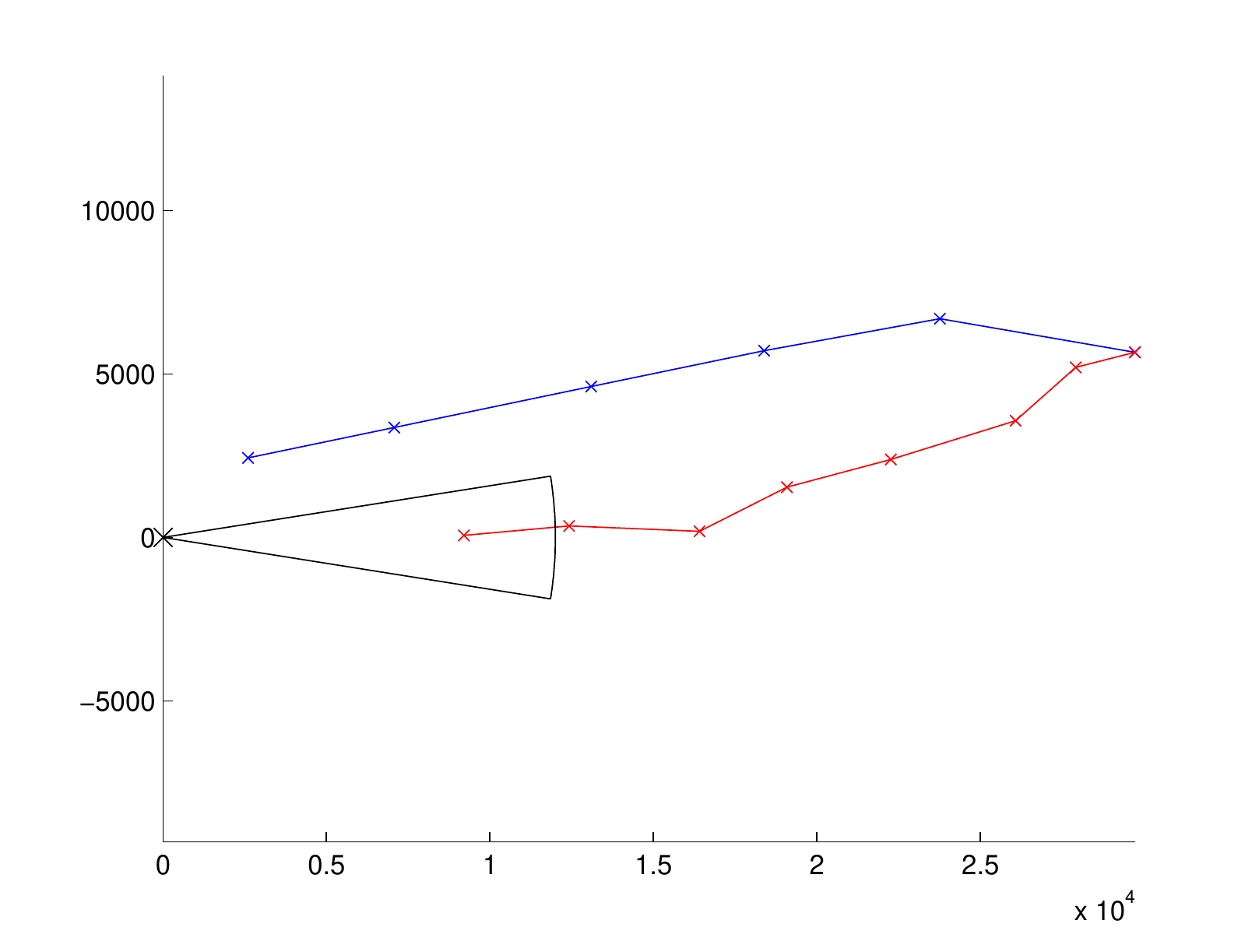}
\caption{Trajectory comparison of 2nd vehicle in first 230 minutes of Sept 14th (blue= real aircraft, red= simulated)}
\label{fig:outlier_1}
\end{center}
\end{figure}

In the first low density simulations (first 230 minutes of Sept 14th) there is one main outlier case where a significant percentage more fuel was used by the simulations compared to the fuel estimates and a second case where both fuel estimates and simulation fuel used were very high. Figure \ref{fig:outlier_1} shows the largest outlier where 620\% more fuel was used by the simulation aircraft. Both paths appear superficially similar though, given the sampling rate, three points for the simulation aircraft (red) should cover a similar distance to one point for the real data (blue). The simulated aircraft is going much faster than the real data. Given the distance needed to travel before the landing envelope is small the objective function in the simulation is likely at fault encouraging an earlier completion time so that the aircraft leaves the problem over the importance of the fuel minimisation itself. Notably in the fuel estimations the fuel estimates for 3 of the 5 steps taken were calculated as negative and recorded as 0, it is likely that this has also impacted on the significant difference in fuel estimates.

\begin{figure}[hptb]
\begin{center}
\def\picwidth{12cm}
\includegraphics[width=\picwidth]{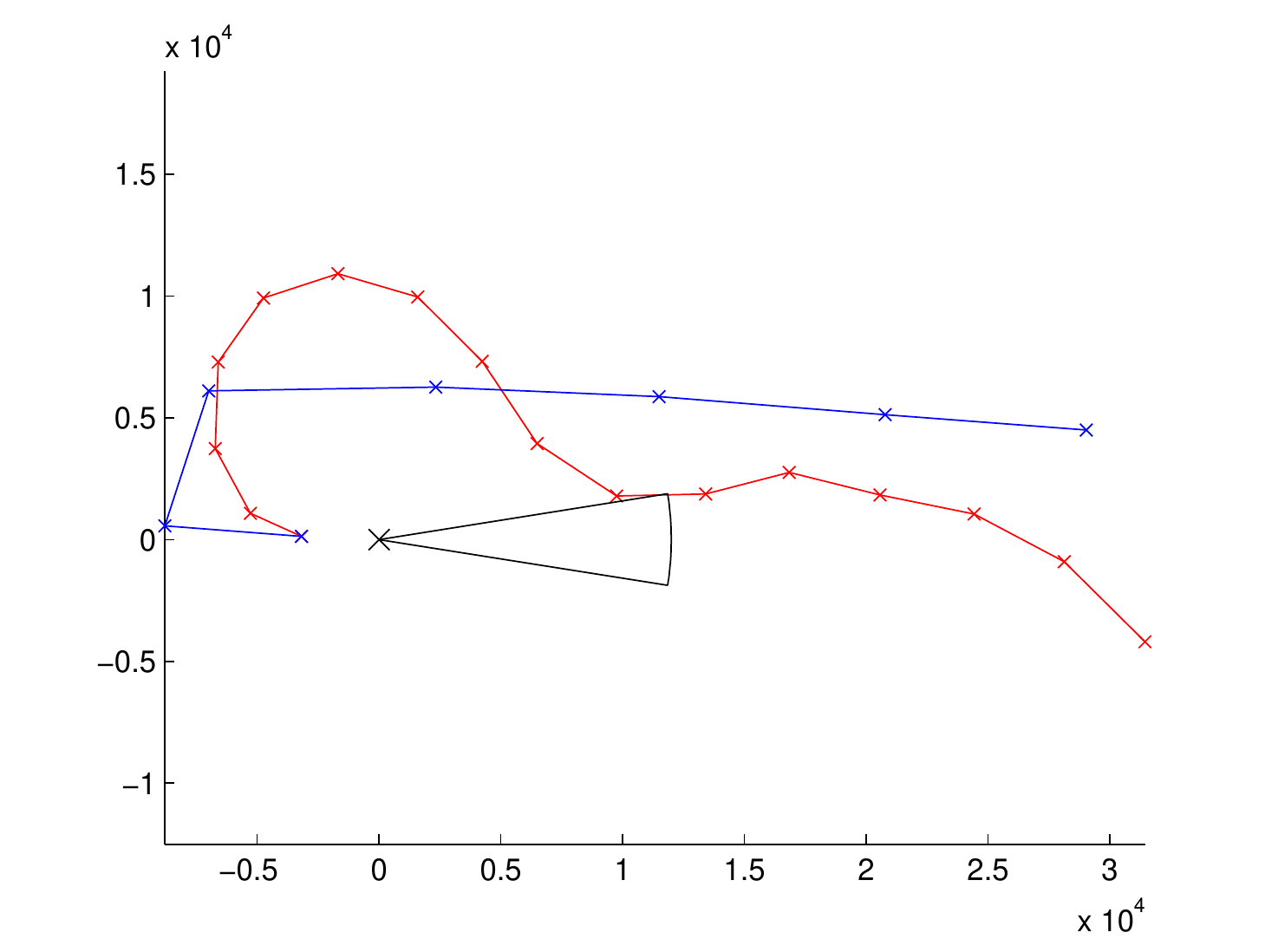}
\caption{Trajectory comparison of 15th vehicle in first 230 minutes of Sept 14th (blue= real aircraft, red= simulated)}
\label{fig:outlier_15}
\end{center}
\end{figure}

The second case where all the fuel estimates were high and the simulation fuel was also high is shown in Figure \ref{fig:outlier_15}. This is the first take-off of the day and does take a longer trajectory than the real data contributing to the increased fuel use.


\begin{figure}[hptb]
\begin{center}
\def\picwidth{12cm}
\includegraphics[width=\picwidth]{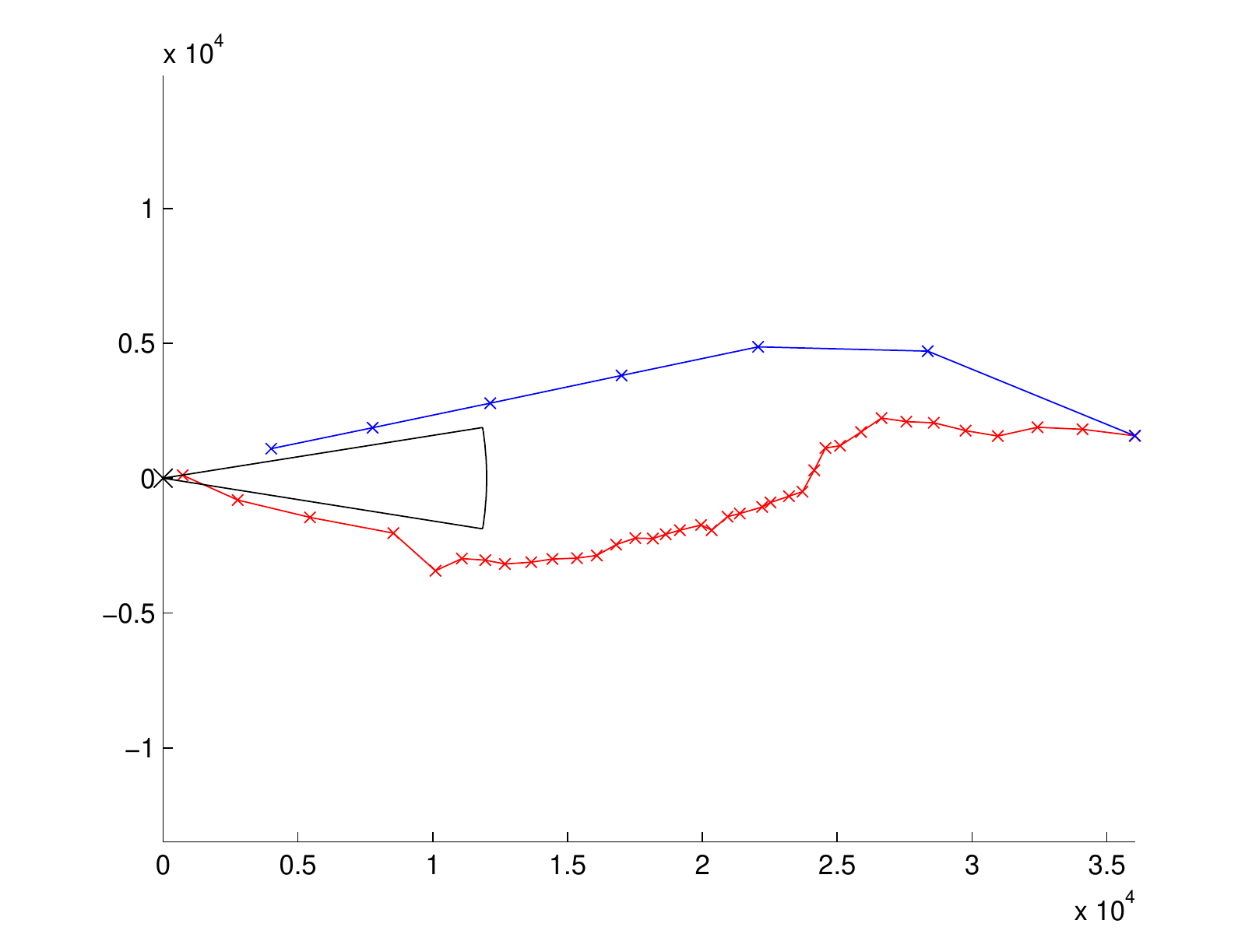}
\caption{Trajectory comparison of 8th vehicle in last 100 minutes of Sept 14th (blue= real aircraft, red= simulated)}
\label{fig:outlier_8}
\end{center}
\end{figure}

Considering the second low density simulations (last 100 minutes of Sept 14th) there are 2 outlier cases where significantly more fuel was used by the simulations than the fuel estimates. Figure \ref{fig:outlier_8} shows the largest outlier where 600\% more fuel was used by the simulation than the mean of the fuel estimations. This appears to be a straightforward path requirement where the aircraft is entering the TMA directly in line with the landing envelope. The simulated aircraft significantly slows down whilst following the nominal descent altitude and appears to keep this slow speed until close to the end of its trajectory where it rapidly descends and accelerates into the flight envelope very close to the airport. The significantly longer time that the simulated aircraft took to land compared to that of the real aircraft (12 minutes vs 6 minutes) might explain the increased fuel use. The behaviour of the simulated aircraft could be due to congestion around the flight envelope from other arrival traffic causing this aircraft to slow down in order to allow others to land before it.

\begin{figure}[hptb]
\begin{center}
\def\picwidth{12cm}
\includegraphics[width=\picwidth]{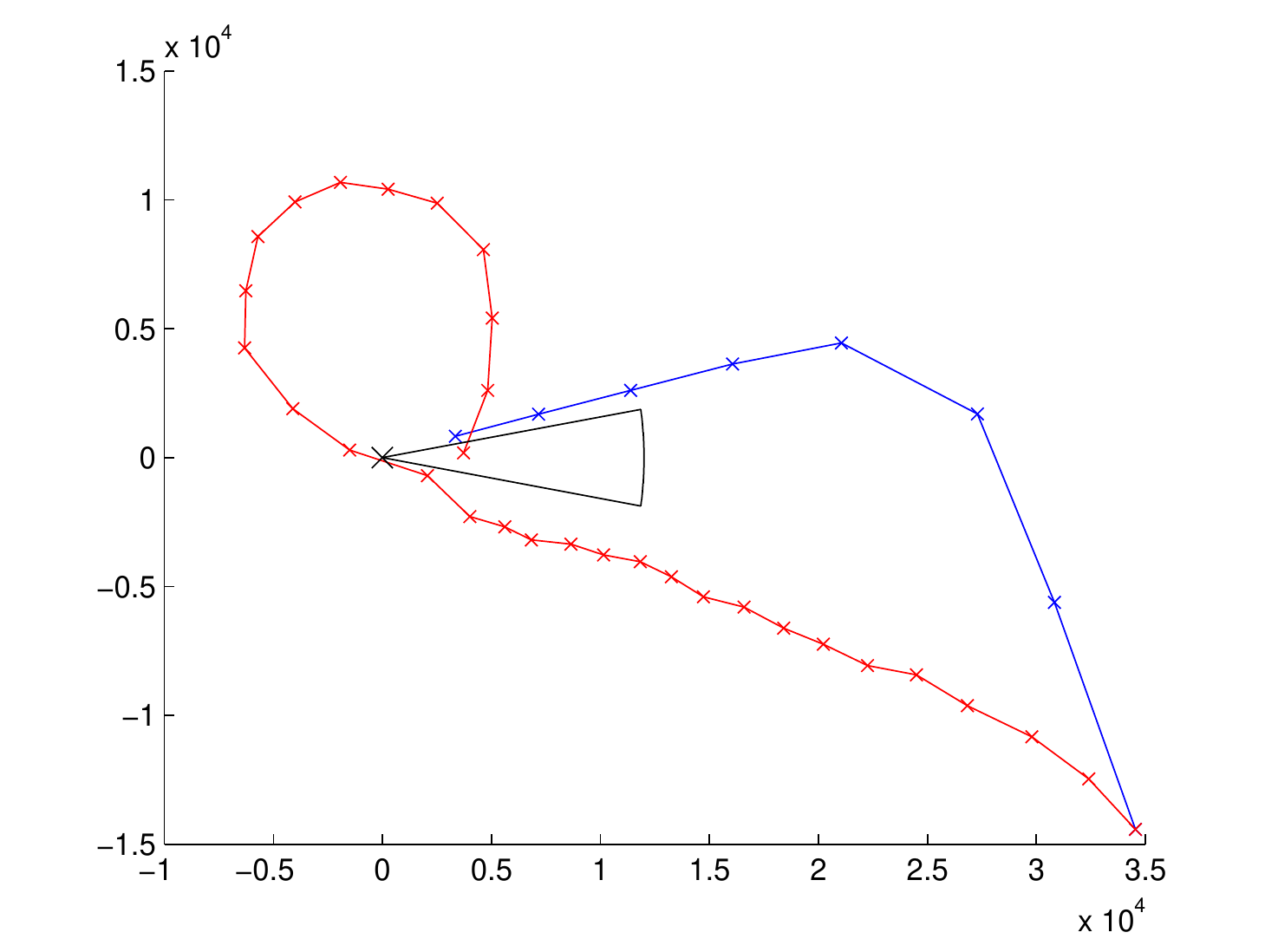}
\caption{Trajectory comparison of 4th vehicle in last 100 minutes of Sept 14th (blue= real aircraft, red= simulated)}
\label{fig:outlier_4}
\end{center}
\end{figure}

The second significant outlier in the second low density simulations is shown in Figure \ref{fig:outlier_4}. Here the failure of the simulated arrival aircraft to enter the landing envelope prior to reaching the airport has resulted in the need for a go around. This would obviously require significantly more fuel than the real aircraft's successful approach.

\subsubsection{High Density Traffic}
The high density traffic scenarios are selected in order to try and push the simulation to the limit. In these cases the airport is operating at maximum capacity and handling an arrival or departure runway movement every minute.

Given that the simulation has no built in runway control, occasional problems with collision avoidance feasibility occur when multiple departure aircraft take off every minute (every 3 simulation time steps). In the highest density scenario considered 2 aircraft were delayed by up to 2 time steps so as to increase the gap between departure aircraft such that they don't break the collision avoidance constraints. Similarly 2 aircraft in the second high density scenario had be delayed outside of the simulation envelope to allow enough space for the remaining aircraft to take off without collision avoidance constraint violations. Aircraft in the simulation tend to turn quickly to their target bearing once they've taken off and the cylindrical avoidance constraints are conservative. In reality the aircraft cannot reach one another due to their turn rates and speeds. A more sophisticated collision avoidance system would likely fix this problem.

\begin{figure}[p]
\begin{center}
\def\picwidth{12cm}
\includegraphics[width=\picwidth]{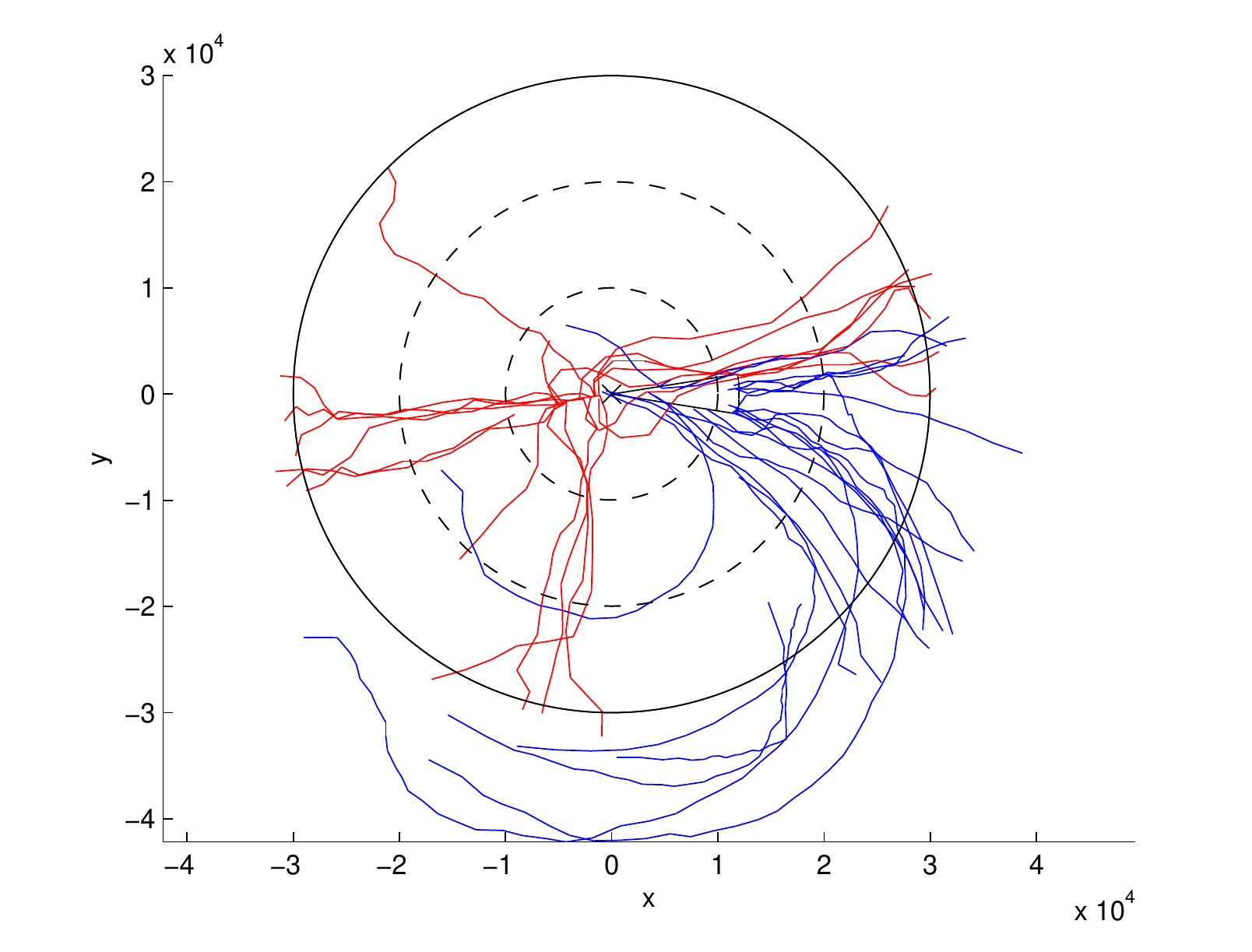}
\caption{Simulated trajectories of aircraft in highest density traffic flow 11.40am-12.40pm of Sept 14th (Blue = arrival aircraft, Red = departure aircraft)}
\label{fig:high_1_traj}
\end{center}
\end{figure}
\begin{figure}[p]
\begin{center}
\def\picwidth{12cm}
\includegraphics[width=\picwidth]{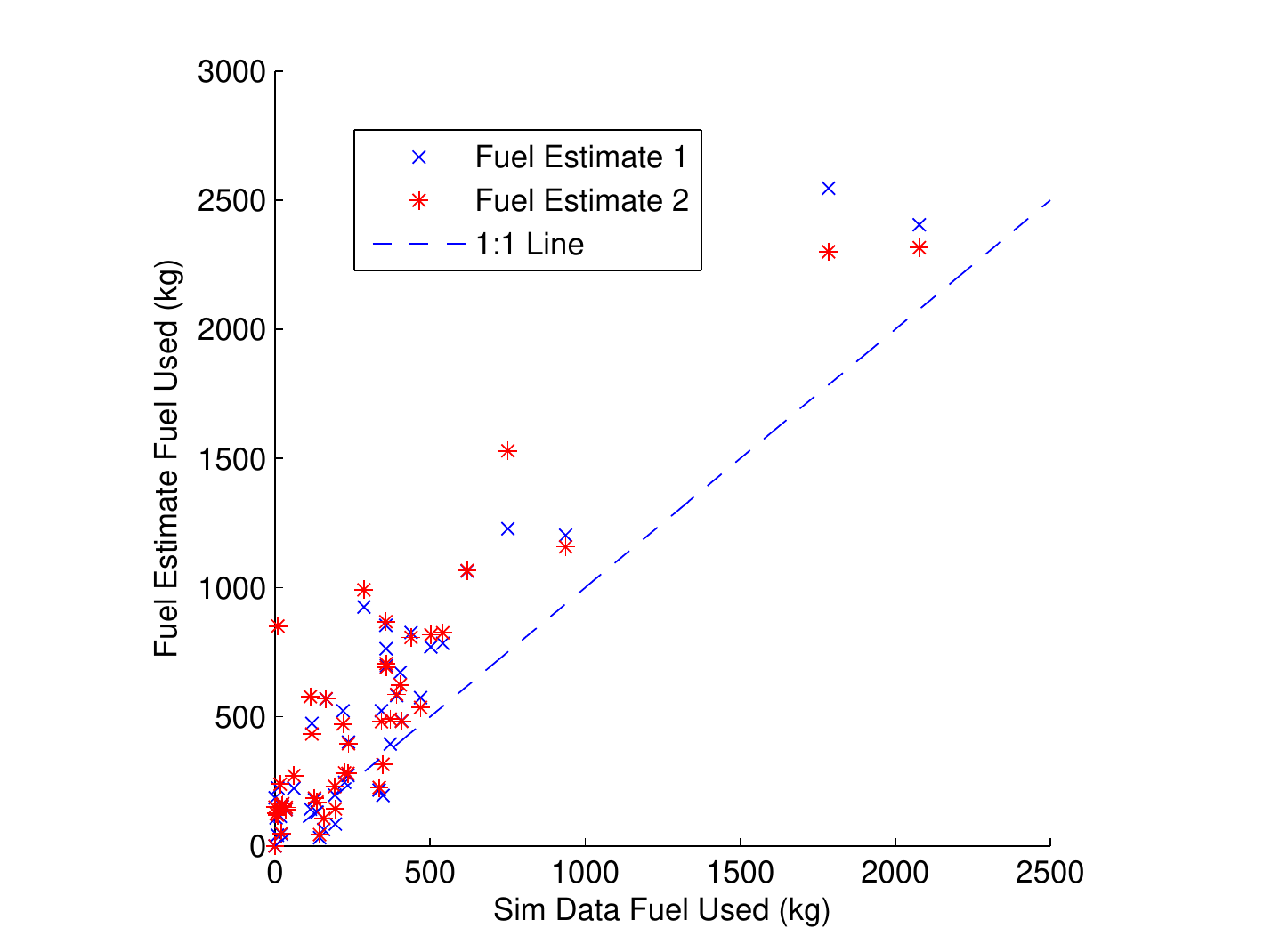}
\caption{Fuel comparison between fuel estimates and simulations for highest density traffic flow 11.40am-12.40pm of Sept 14th}
\label{fig:high_1_fuel}
\end{center}
\end{figure}

\begin{figure}[p]
\begin{center}
\def\picwidth{12cm}
\includegraphics[width=\picwidth]{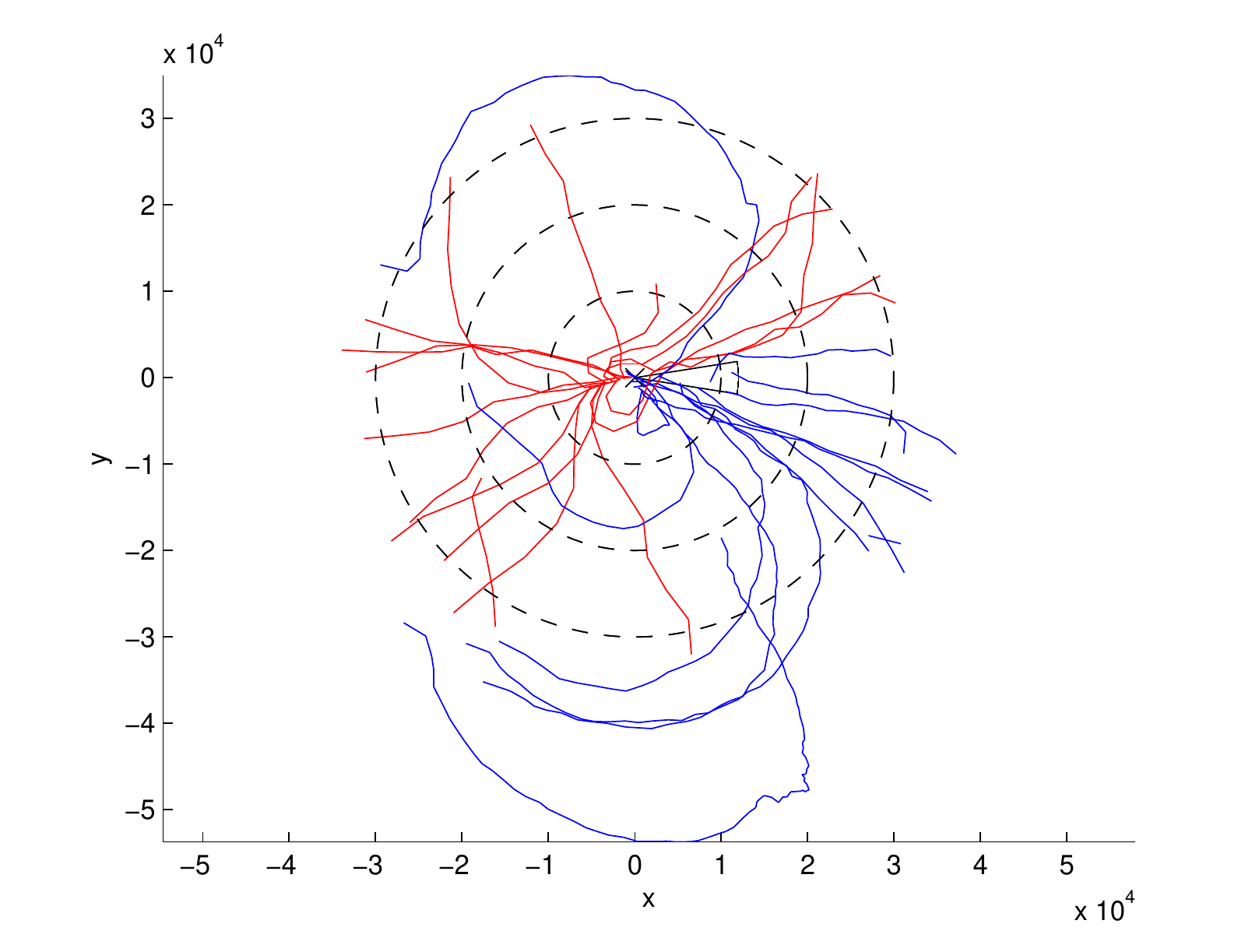}
\caption{Simulated trajectories of aircraft in 4.40pm-5.40pm of Sept 14th (Blue = arrival aircraft, Red = departure aircraft)}
\label{fig:high_2_traj}
\end{center}
\end{figure}
\begin{figure}[p]
\begin{center}
\def\picwidth{12cm}
\includegraphics[width=\picwidth]{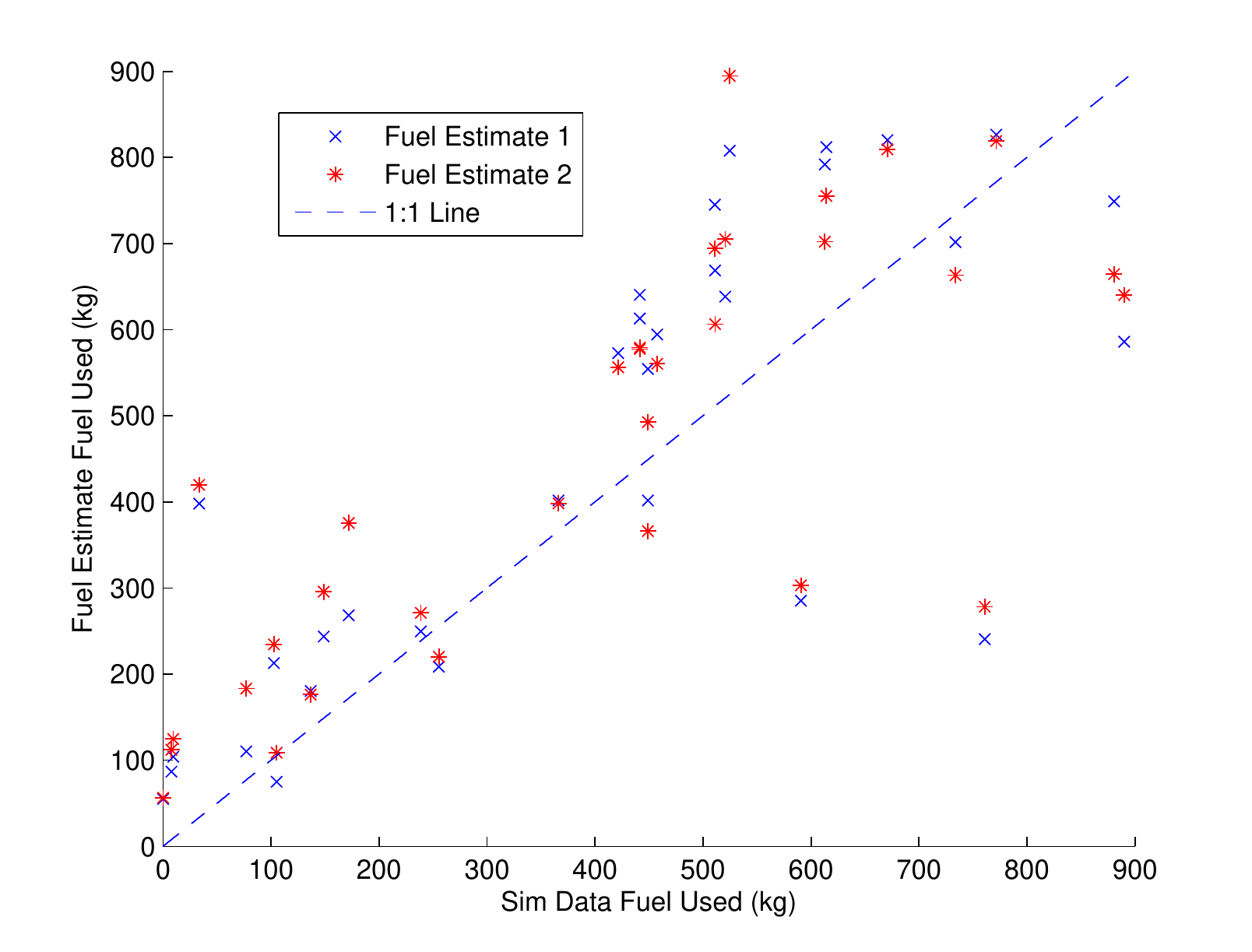}
\caption{Fuel comparison between fuel estimates and simulations in 4.40pm-5.40pm of Sept 14th}
\label{fig:high_2_fuel}
\end{center}
\end{figure}

The final trajectories of the high density aircraft traffic are shown in Figures \ref{fig:high_1_traj} and \ref{fig:high_2_traj}. Figures \ref{fig:high_1_fuel} and \ref{fig:high_2_fuel} show the fuel comparisons between the fuel estimates obtained from the data and the simulated fuel use along with Tables \ref{tab:high_1_fuel} and \ref{tab:high_2_fuel}. In both cases the majority of aircraft do make fuel savings from the simulation compared to the fuel estimates. The first high density scenario (11.40am - 12.40pm) is especially good for fuel savings with only some very short paths failing to make a saving over the estimate. The second high density scenario shows a greater variation of fuel savings with a number of paths failing to make fuel savings over estimates but only two paths show significant outlying behaviour from the norm.

\begin{table}[phtb]
  \centering
\begin{tabular}{ c | c | c | c | c| c}
  \hline                        
Aircraft No. & $F_S$ (kg) & $F_1$ (kg) & $F_2$ (kg) & Mean $F$ (kg) & Fuel Saving \% \\
\hline
1& 6.75	&41.20	&118.71	&79.96	&91.55\\
2& 20.09	&45.35	&47.99	&46.67	&56.94\\
3& 143.47	&31.78	&43.06	&37.42	&-283.43\\
4& 60.46	&222.96	&271.58	&247.27	&75.55\\
5& 370.88	&394.79	&491.19	&442.99	&16.28\\
6& 236.31	&401.19	&395.79	&398.49	&40.70\\
7& 35.04	&141.03	&148.90	&144.96	&75.83\\
8& 937.13	&1203.14	&1158.40	&1180.76	&20.63\\
9& 15.94	&114.13	&238.54	&176.34	&90.96\\
10& 347.35	&194.67	&315.22	&254.94	&-36.25\\
11& 391.88	&581.55	&587.00	&584.28	&32.93\\
12& 224.07	&245.99	&282.36	&264.18	&15.18\\
13& 469.25	&573.37	&536.30	&554.83	&15.42\\
14& 286.07	&924.53	&991.72	&958.13	&70.14\\
15& 135.32	&130.35	&168.95	&149.65	&9.57\\
16& 357.87	&763.70	&705.73	&734.71	&51.29\\
17& 540.73	&784.48	&824.65	&804.56	&32.79\\
18& 113.22	&142.08	&577.63	&359.86	&68.54\\
19& 2077.19	&2405.27	&2317.04	&2361.16	&12.03\\
20& 502.14	&770.08	&817.34	&793.71	&36.74\\
21& 335.33	&215.76	&225.92	&220.84	&-51.85\\
22& 437.88	&826.45	&807.23	&816.84	&46.39\\
23& 193.77	&84.12	&143.55	&113.83	&-70.22\\
24& 407.62	&484.46	&482.49	&483.48	&15.69\\
25& 156.55	&63.27	&106.31	&84.79	&-84.63\\
26& 358.29	&697.57	&691.05	&694.31	&48.40\\
27& 127.01	&179.45	&184.48	&181.97	&30.20\\
28& 8.36	&223.60	&851.40	&537.50	&98.44\\
29& 1785.12	&2546.99	&2300.89	&2423.95	&26.35\\
30& 192.42	&196.04	&229.73	&212.88	&9.61\\
31& 356.64	&853.68	&868.37	&861.02	&58.58\\
32& 751.03	&1227.82	&1529.48	&1378.65	&45.52\\
33& 234.34	&272.95	&280.80	&276.87	&15.36\\
34& 619.16	&1064.15	&1066.75	&1065.45	&41.89\\
35& 342.84	&522.50	&481.20	&501.85	&31.68\\
36& 34.08	&140.89	&139.81	&140.35	&75.72\\
37& 164.20	&568.39	&570.20	&569.30	&71.16\\
38& 403.58	&671.65	&623.49	&647.57	&37.68\\
39& 23.49	&163.07	&163.19	&163.13	&85.60\\
40& 23.23	&151.44	&155.20	&153.32	&84.85\\
41& 218.78	&522.63	&471.82	&497.23	&56.00\\
42& 118.15	&474.08	&433.16	&453.62	&73.95\\
43& 2.97	&107.74	&125.47	&116.60	&97.46\\
\hline
Total &14566 &22370 &23970 & 23170 &37.13\\
\hline  
\end{tabular}
\caption{Fuel use comparison between simulation and estimates for 11.40am to 12.40pm of Sept 14th}
  \label{tab:high_1_fuel}
\end{table}

\begin{table}[phtb]
  \centering
\begin{tabular}{ c | c | c | c | c| c}
  \hline                        
Aircraft No. & $F_S$ (kg) & $F_1$ (kg) & $F_2$ (kg) & Mean $F$ (kg) & Fuel Saving \% \\
\hline
1	&9.42		&104.07	&124.51	&114.29	&91.76\\
2	&105.33	&74.91	&108.67	&91.81	&-14.73\\
3	&366.33	&401.63	&398.24	&399.94	&8.40\\
4	&102.73	&212.76	&234.45	&223.61	&54.06\\
5	&148.83	&243.54	&295.60	&269.57	&44.79\\
6	&612.76	&791.94	&702.23	&747.09	&17.98\\
7	&255.55	&208.69	&219.98	&214.34	&-19.23\\
8	&511.21	&668.67	&606.37	&637.51	&19.81\\
9	&520.56	&638.37	&705.13	&671.75	&22.51\\
10	&171.92	&268.15	&375.57	&321.86	&46.59\\
11	&77.09	&110.13	&183.20	&146.67	&47.44\\
12	&524.75	&807.87	&894.76	&851.31	&38.36\\
13	&511.04	&744.91	&694.32	&719.62	&28.98\\
14	&457.58	&594.28	&560.27	&577.28	&20.74\\
15	&448.99	&554.31	&492.59	&523.45	&14.23\\
16	&136.85	&180.33	&176.21	&178.27	&23.24\\
17	&670.78	&820.13	&809.46	&814.80	&17.68\\
18	&590.81	&285.33	&303.03	&294.18	&-100.84\\
19	&238.68	&249.56	&270.99	&260.27	&8.30\\
20	&614.25	&811.95	&755.17	&783.56	&21.61\\
21	&441.54	&640.34	&577.50	&608.92	&27.49\\
22	&760.94	&240.40	&278.13	&259.27	&-193.49\\
23	&421.44	&572.56	&556.21	&564.39	&25.33\\
24	&733.80	&701.65	&663.05	&682.35	&-7.54\\
25	&7.79		&86.75	&112.10	&99.42	&92.17\\
26	&771.69	&826.40	&819.08	&822.74	&6.21\\
27	&33.47	&397.96	&419.55	&408.75	&91.81\\
28	&880.66	&748.89	&664.47	&706.68	&-24.62\\
29	&890.05	&586.08	&640.08	&613.08	&-45.18\\
30	&441.65	&612.92	&579.21	&596.07	&25.91\\
31	&449.05	&401.71	&366.04	&383.87	&-16.98\\
32	&0.13		&55.08	&56.32	&55.70	&99.78\\
\hline
Total 	&12908 	&14642 	&14643 	&14642 	&11.85\\
\hline  
\end{tabular}
\caption{Fuel use comparison between simulation and estimates for 4.40pm to 5.40pm of Sept 14th}
  \label{tab:high_2_fuel}
\end{table}

\begin{figure}[hptb]
\begin{center}
\def\picwidth{12cm}
\includegraphics[width=\picwidth]{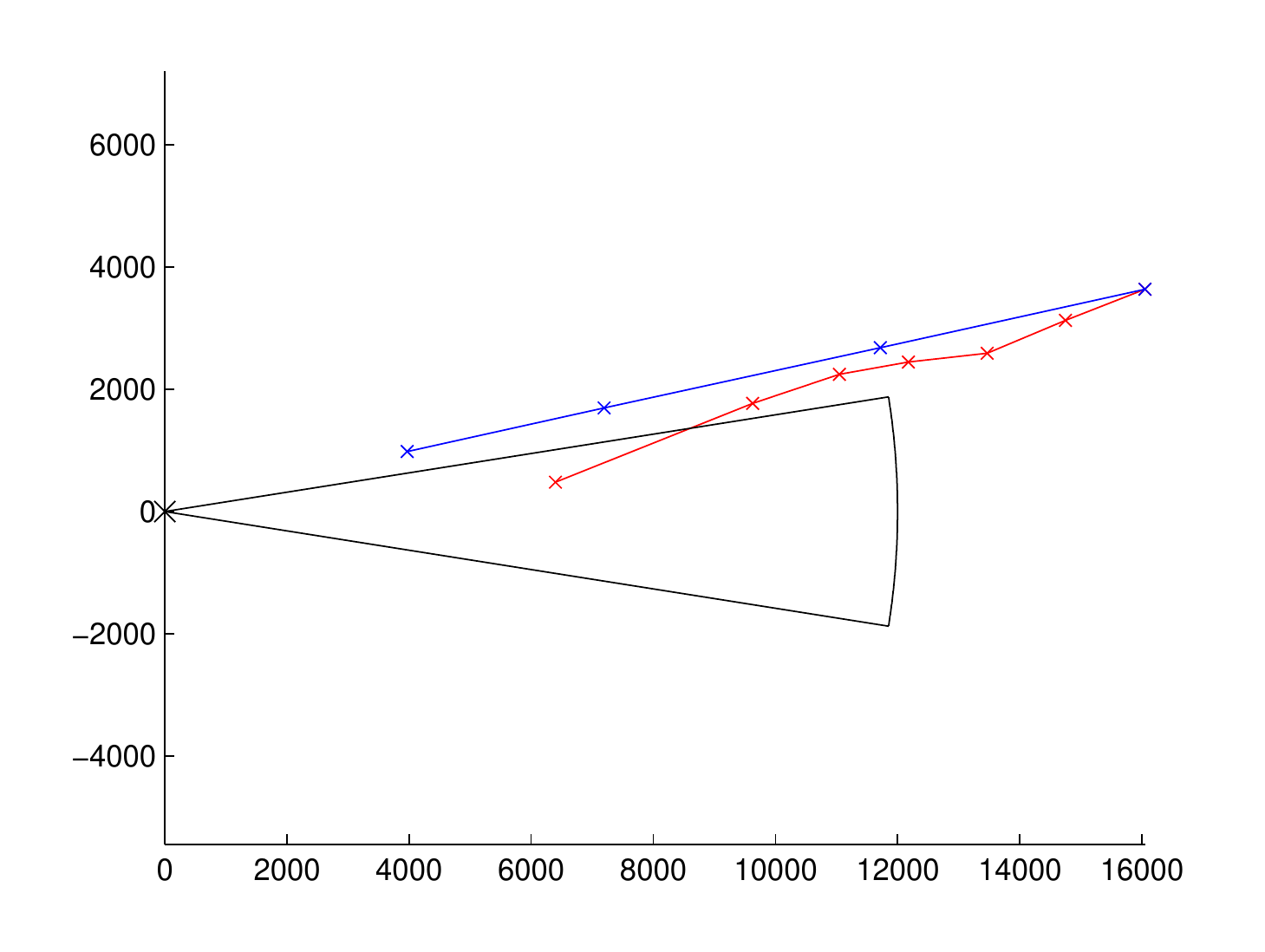}
\caption{Trajectory comparison of 3rd vehicle in 11.40am-12.40pm scenario Sept 14th (blue= real aircraft, red= simulated)}
\label{fig:outlier_3}
\end{center}
\end{figure}

Figure~\ref{fig:outlier_3} shows the single outlier case from the first high density simulation. Where in one of the shortest paths in the simulation 283\% more fuel has been used by the simulation than the fuel estimates. Similar to the case discussed in the low density simulations (Figure \ref{fig:outlier_1}) the simulated aircraft does a large airspeed acceleration the step before it enters the landing envelope. This serves to rush entering the envelope and get an improved cost for finishing sooner over a slower more fuel efficient approach. Since both fuel estimates are very low for this section of the path this rash acceleration shows up very poorly in the percentage fuel saving metric.

\begin{figure}[hptb]
\begin{center}
\def\picwidth{12cm}
\includegraphics[width=\picwidth]{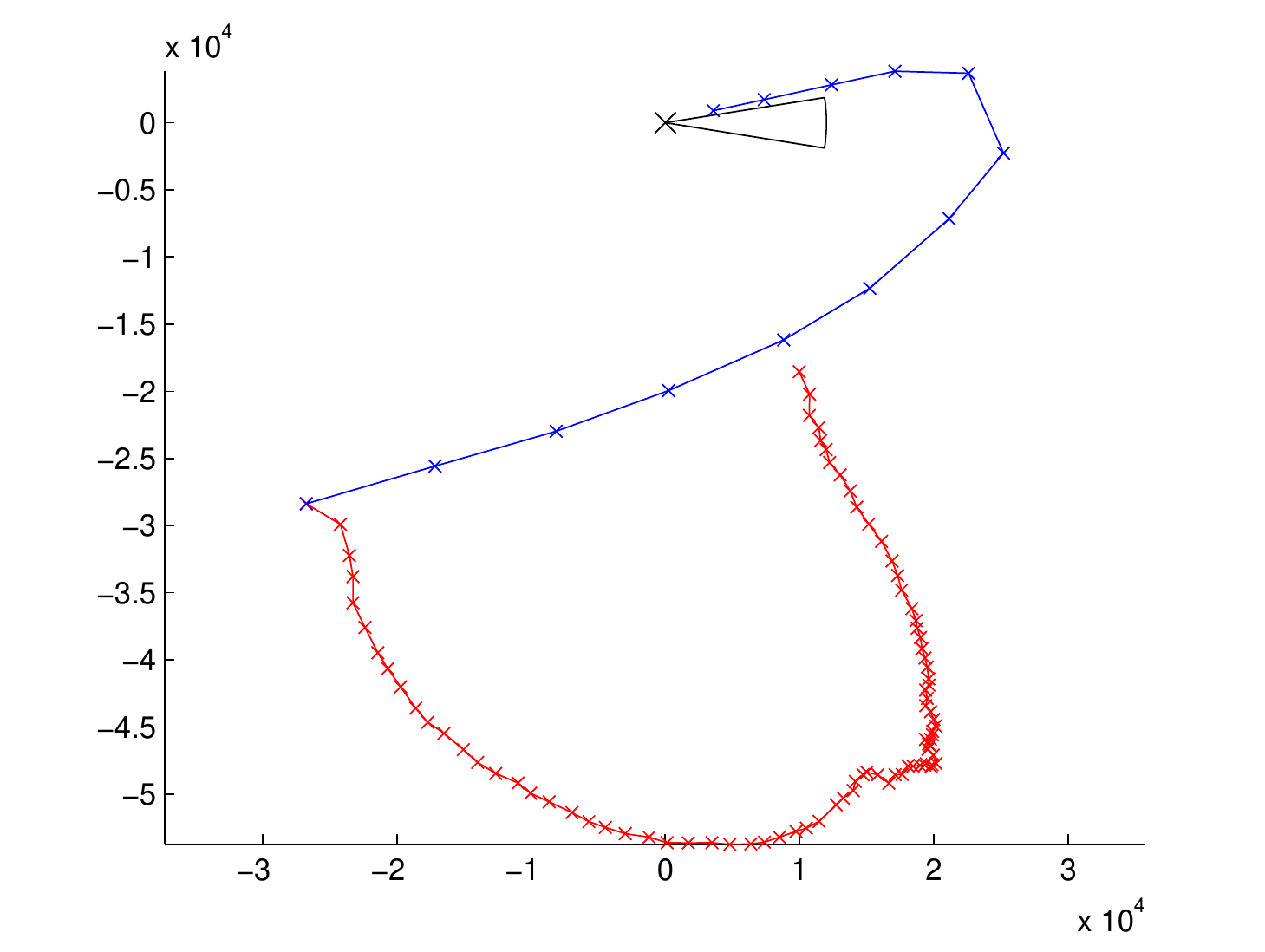}
\caption{Trajectory comparison of 18th vehicle in 4.40pm-5.40pm scenario Sept 14th (blue= real aircraft, red= simulated)}
\label{fig:outlier_18}
\end{center}
\end{figure}

Figure~\ref{fig:outlier_18} shows the first outlier case from the second high density simulation where an arrival aircraft has taken a very circuitous route towards the landing envelope. Initially pushed slightly off course by other aircraft also arriving at a similar point in the TMA the route taken is significantly further south than the standard arrival path. Additionally the aircraft has undesirable behaviour at the point the aircraft begins to turn back towards the airport and the landing envelope. Currently this path is recorded as having a 100\% increase in fuel use but since the aircraft has yet to land in the simulation window this would be far higher in practice. Altering the flow field to encourage less of a diversion to this arrival aircraft would improve matters. Lowering the cost function weighting to the flow field could possibly improve matters but only when used alongside an increase in horizon length.

\begin{figure}[hptb]
\begin{center}
\def\picwidth{12cm}
\includegraphics[width=\picwidth]{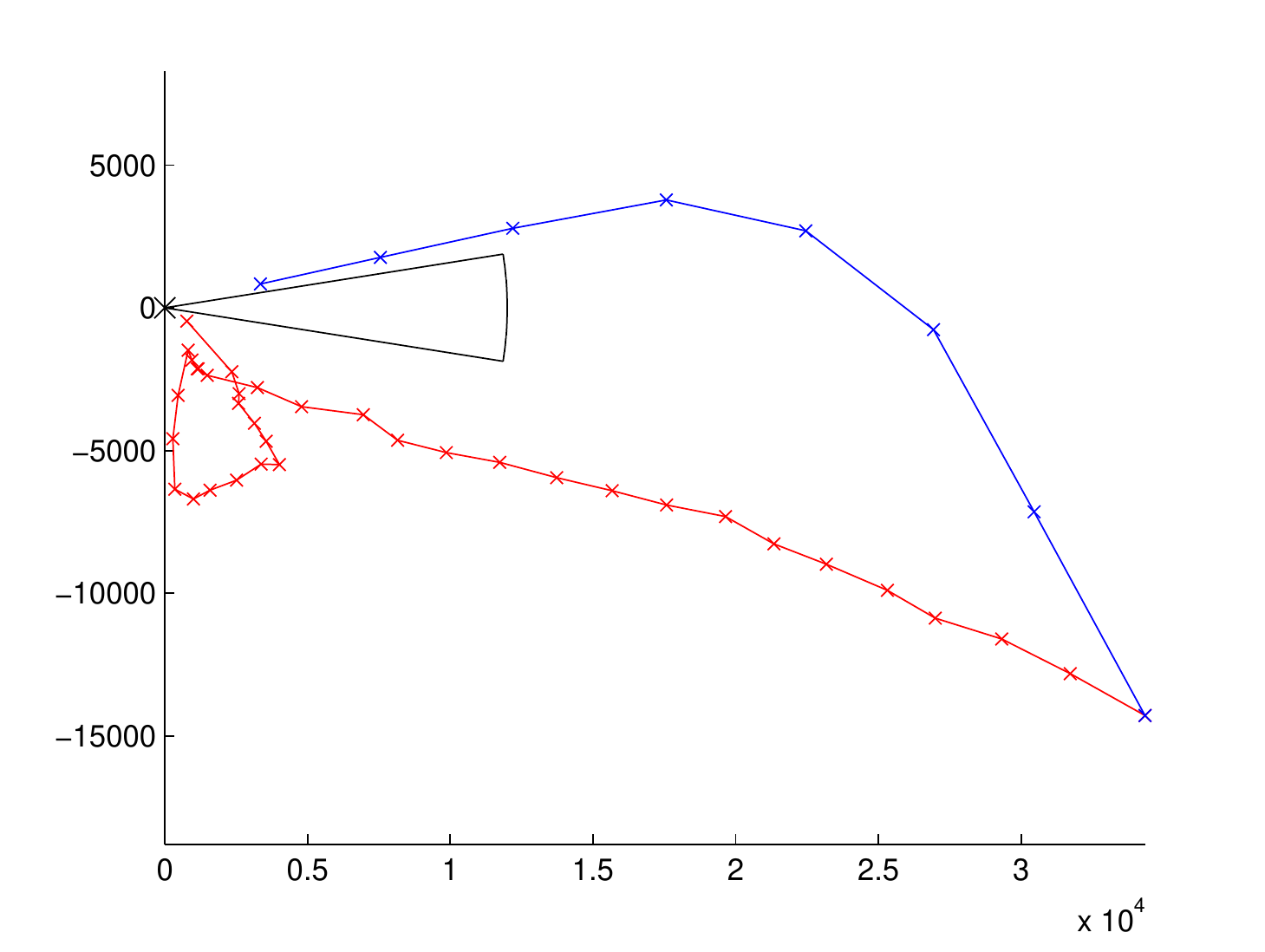}
\caption{Trajectory comparison of 22nd vehicle in 4.40pm-5.40pm scenario Sept 14th (blue= real aircraft, red= simulated)}
\label{fig:outlier_22}
\end{center}
\end{figure}

Figure~\ref{fig:outlier_22} shows the second outlier case from the second high density simulation where the 22nd vehicle of the simulation takes an aborted go-around to reach the landing envelope. This go-around appears to have been motivated by a collision avoidance manoeuvre and explains the 193\% increase in fuel used over the fuel estimates which follow a far simpler and smoother path.

\subsubsection{Overall Fuel and CO$_2$ Savings}
The summary tables for all four scenarios (both high and low density) do demonstrate an overall fuel saving summed across all aircraft with a magnitude of between 4.29\% to 53.93\%. In total across all four scenarios 14653.5 kg of fuel was saved. This leads to an overall fuel saving across every aircraft simulated of 29.7\%. On average during flight 3.15g of CO$_2$ is produced for each gram of aircraft fuel burnt \cite{CORINAIR} leading to a saving of 46158.52 kg of CO$_2$.

\section{Model Extension}
\subsection{Noise Pollution Reduction}
\begin{figure}[htbp]
\begin{center}
\def\picwidth{13cm}
\includegraphics[width=\picwidth]{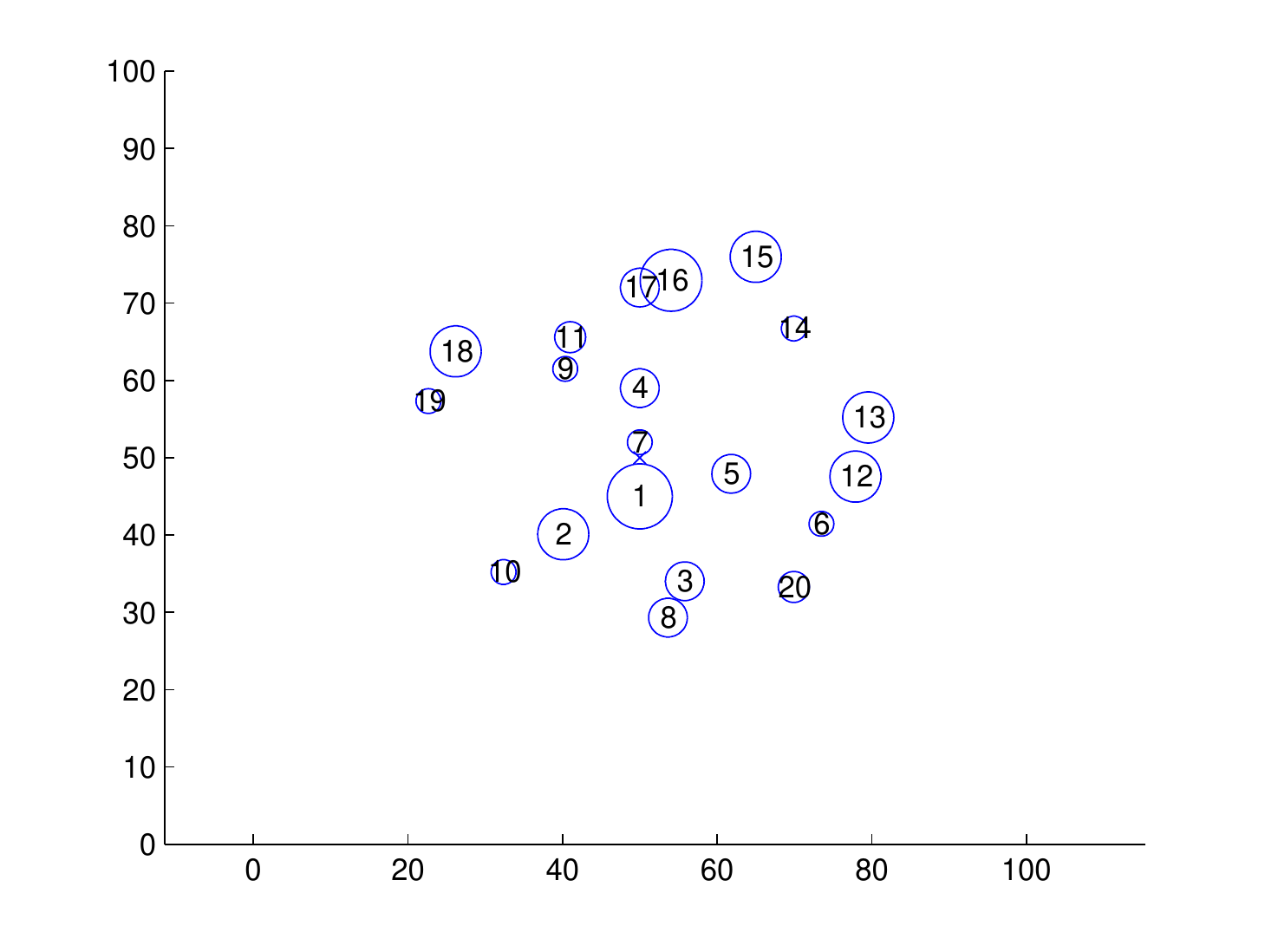}
\caption{Location of 20 population centres around Gatwick airport and their relative size. $1=$ Crawley, $2=$ Horsham, $3=$ Haywards Heath, $4=$ Reigate, $5=$ East Grinstead, $6=$ Crowborough, $7=$ Horley, $8=$ Burgess Hill, $9=$ Dorking, $10=$ Billinghurst, $11=$ Leatherhead, $12=$ Tunbridge Wells, $13=$ Tonbridge, $14=$ Sevenoaks, $15=$ Orpington, $16=$ Croydon, $17=$ Sutton, $18=$ Guildford, $19=$ Godalming, $20=$ Uckfield.}
\label{fig:pop_circles}
\end{center}
\end{figure}

\begin{figure}[bhtp]
\begin{center}
\def\picwidth{12cm}
\includegraphics[width=\picwidth]{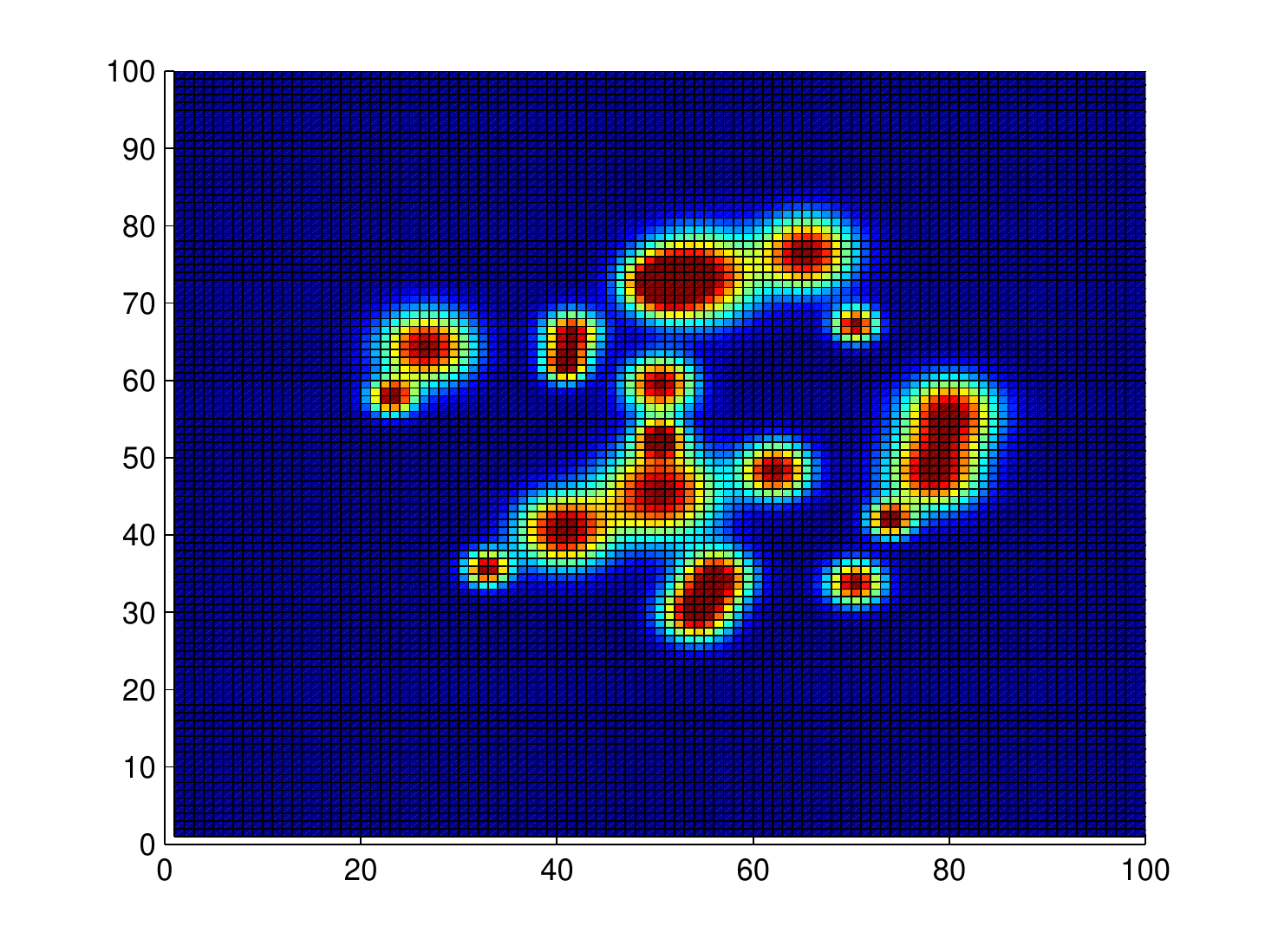}
\caption{Population density estimations in local area surrounding Gatwick airport}
\label{fig:pop_dense}
\end{center}
\end{figure}

For a noise pollution reduction cost a basic population density distribution was generated over the area surrounding Gatwick airport. Through using Google Maps satellite images to identify areas of population (assumed to be grey built up areas from the satellite map) which were equal to or greater than 1.5km in radius and then noting their co-ordinates relative to the airport 20 population centres were identified in the TMA area (shown in Figure \ref{fig:pop_circles}).
\begin{figure}[htbp]
\begin{center}
\def\picwidth{14cm}
\includegraphics[width=\picwidth,viewport=0.5in 8in 7in 11in,clip=true]{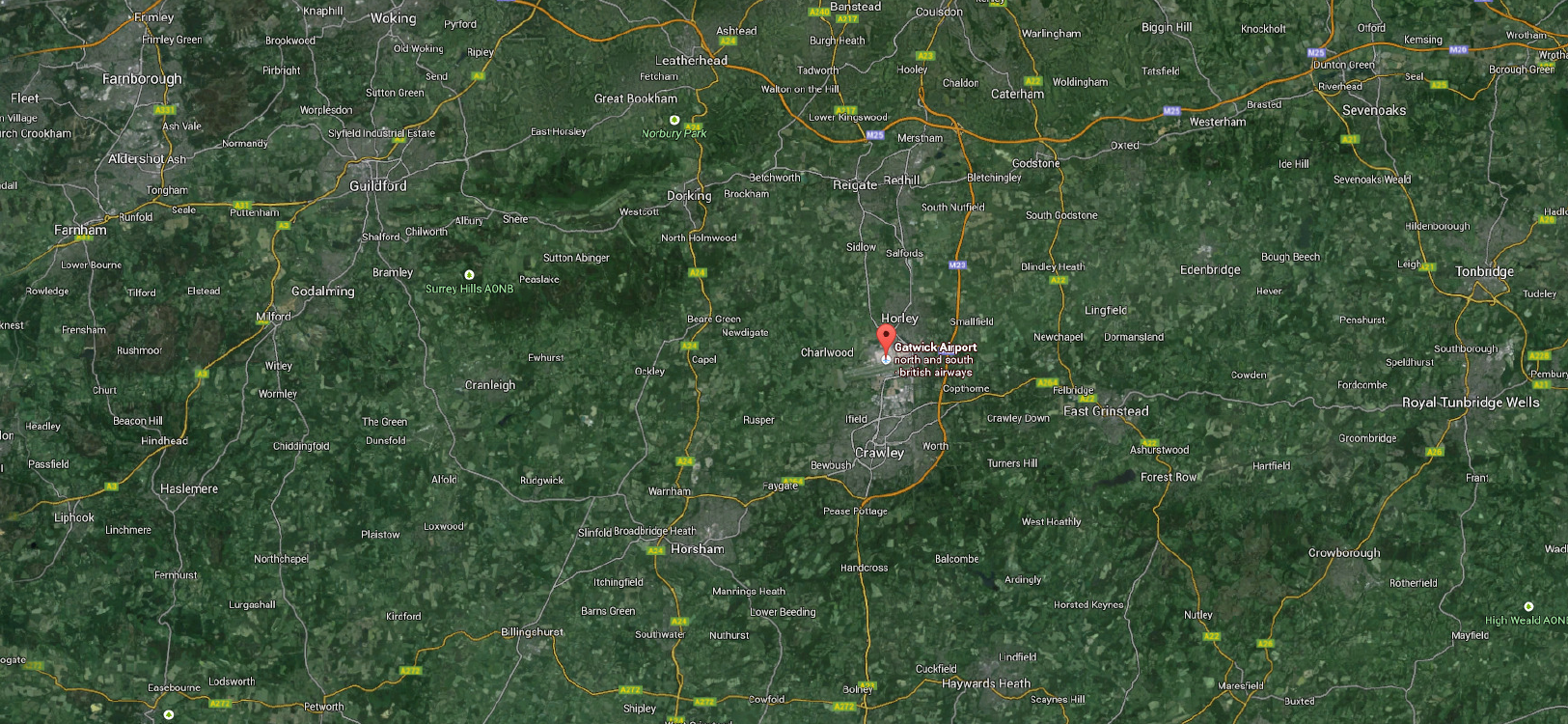}
\caption{Map of the local area around Gatwick airport (courtesy of GoogleMaps)}
\label{fig:gatwick_map}
\end{center}
\end{figure}
A grid of 1km spacing was then imposed across the entire area and the population density in each location was estimated using a Gaussian function:
\begin{equation}
\mathrm{popdense}(x,y)=\min \left(\sum_{i=1}^{20} {{e^{-{{b_i(x,y)}^2}\over{2{c_i}^2}}}\over{c_i\sqrt{2\pi}} }, 1\right).
\end{equation}
Where $b_i$ was taken as the distance from the grid point to the i-th population centre and $c_i$ was the radius of that i-th population centre. This results in a population density map shown in Figure \ref{fig:pop_dense} where brighter areas correlate with areas of higher population density. Ideally aircraft would not fly over populated areas, or if they did so, would fly at heights such that their noise pollution is negligible. So to augment the population density grid a cost on the aircraft based on their altitude is also included.
\begin{equation}
J_{\mathrm{noise}}=1-\max\left(\left(1-\left( 1- {{A_c-a}\over{A_c}}\right)^2\right),0\right)\mathrm{popdense}(x,y),
\end{equation}
where $A_c$ is effectively an altitude cutoff above which the noise of the aircraft begins to be negligible to the population on the ground. Experimentally this has been set at 4000m. The quadratic altitude element heavily penalises aircraft flying low to the ground. If aircraft fly above the altitude cutoff the quadratic element will go to 0 and thus $J_{\mathrm{noise}}$ will stay at 1 regardless of population density. This cost function could also be further augmented with terms relating the thrust of an aircraft's engines and the amount of noise produced. This cost function is then weighted based on the relative importance of noise pollution reduction vs other cost terms like following the flow field, following the nominal altitude and fuel saving.

This implementation for a noise reduction reward is an extension of a basic penalty method which could be trivially used to select no fly zones or minimum altitude constraints around an airport. Specific no fly zones could also be implemented through use of hard constraints. However softer constraints by penalising cost are easier for the SMC method to handle and, given enough time, will converge to an optimal solution.

\subsection{Noise Pollution Simulation Results}
To test the performance of the noise pollution cost function and its effects on aircraft trajectories 22 individual scenarios were run. Firstly with noise pollution reduction carrying a weighting of 10\% importance in the cost function and then with 20\% importance. The 22 scenarios chosen were single aircraft simulations similar to those used to test the original cost functions in Section \ref{sec:SAT}. With the exception of 2 missing scenarios where either an arrival aircraft flies directly into the path of take off aircraft or where a departure aircraft is asked to leave in the direction of the runway landing envelope. Single aircraft trajectories were chosen for this testing in order to remove behaviour caused by interaction between aircraft and limit the trajectory changes seen to those directly caused by the noise pollution element of the cost function. All simulations were done under the same wind conditions and with the same type of aircraft. Arrival aircraft were arranged to start at 30\degree intervals around the perimeter of the TMA and told to head to the landing envelope. Departure aircraft all start from the airport and are given a desired bearing to leave the TMA. Notably departure aircraft do not have to reach their desired bearing to leave the TMA and are considered to leave the TMA once they pass over 30km from the airport regardless of what bearing they are currently at.

\begin{figure}[bhtp]
\begin{center}
\def\picwidth{12cm}
\includegraphics[width=\picwidth]{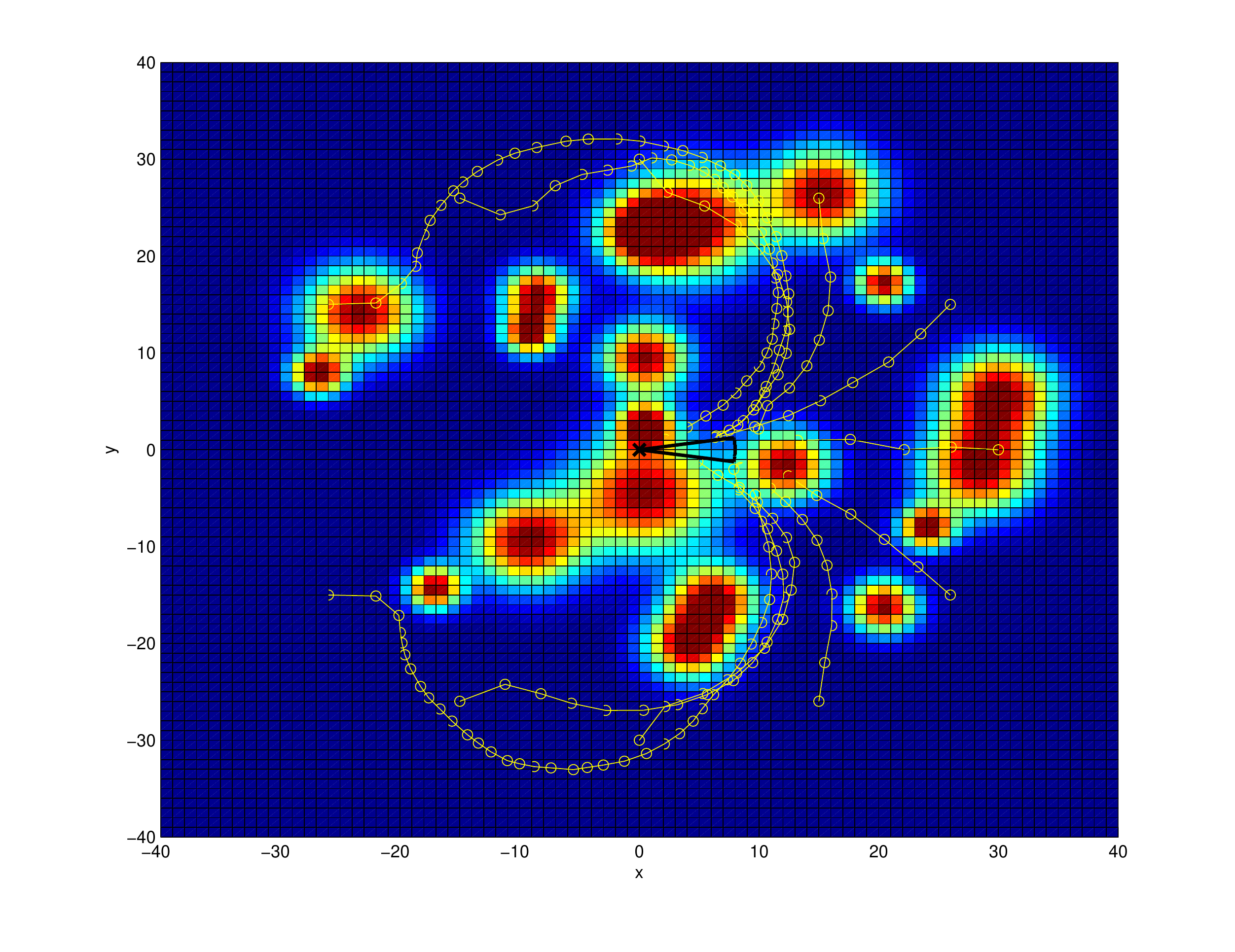}
\caption{Trajectories for 11 separate arrival simulations in presence of 10\% weighted noise pollution reduction cost}
\label{fig:10arrive}
\end{center}
\end{figure}
\begin{figure}[bhtp]
\begin{center}
\def\picwidth{12cm}
\includegraphics[width=\picwidth]{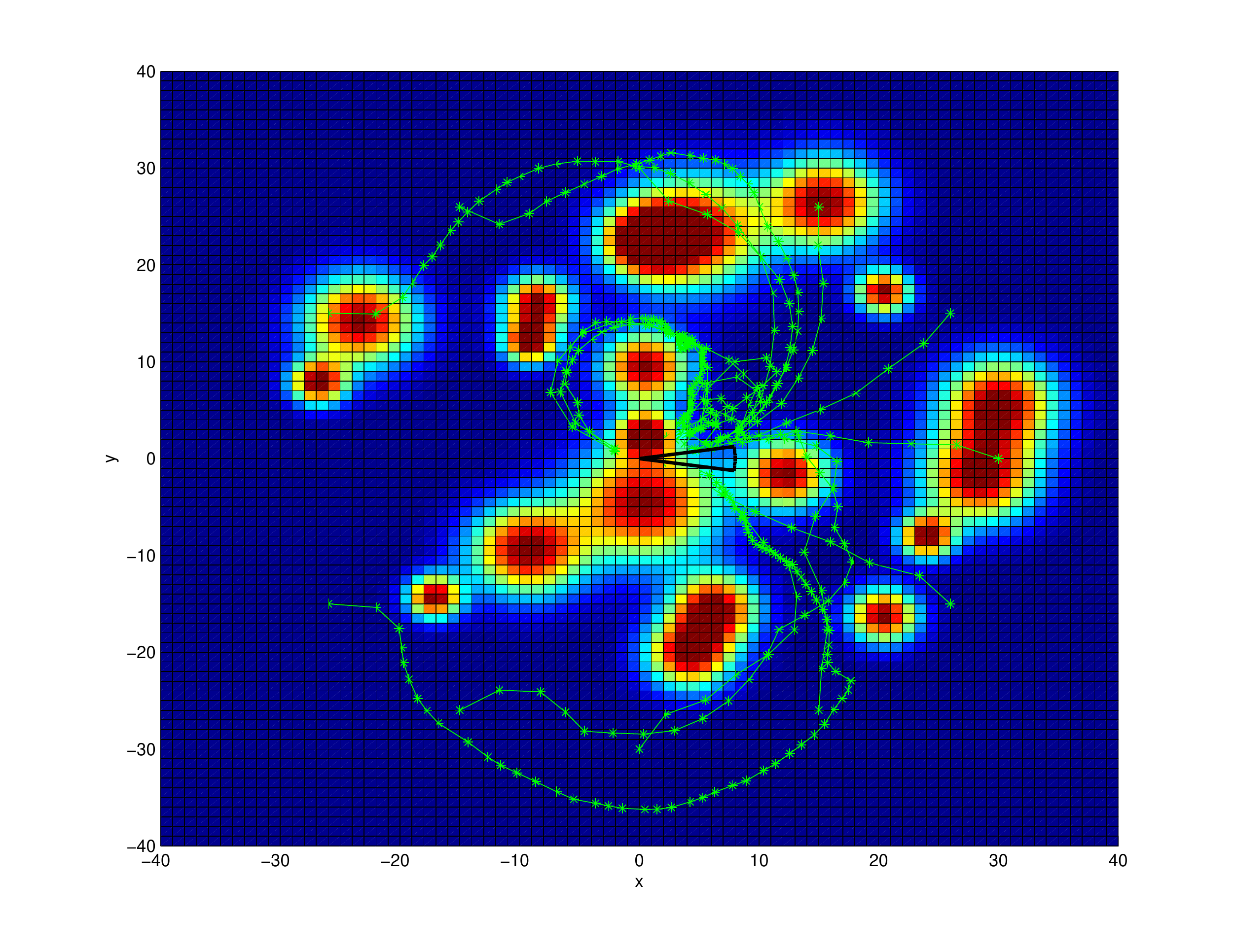}
\caption{Trajectories for 11 separate arrival simulations in presence of 20\% weighted noise pollution reduction cost}
\label{fig:20arrive}
\end{center}
\end{figure}

Figure \ref{fig:10arrive} shows the arrival aircraft simulations over the top of the population density map for 10\% and Figure \ref{fig:20arrive}20\% noise pollution cost waitings. Compared to the pure flow field with fuel reduction cost shown in Figure \ref{fig:arrive} (which operated on a shorter sampling time of 10 seconds instead of the 20 seconds used for the noise pollution simulation) there is a significant change in pathing of aircraft to the landing envelope. In the case of the 10\% performance weighting for noise pollution the paths are behave more like the original non augmented cost function than those for the 20\% importance. This is due to the 20\% importance weighting requiring a weighting reduction in many elements of the cost function which otherwise would ensure smooth aircraft pathing via the flow field and nominal altitude profiles. In the 20\% importance simulations at least 3 of the aircraft had to complete go-around manoeuvres to reach the landing envelope successfully. This behaviour was limited to aircraft entering from the southern side of the airport where they had to avoid a population zone at East Grinstead and could do so by either heading to the left or right of the region. In the case of the 10\% importance weighting aircraft were still prone to pass over East Grinstead in order to enter the landing envelope faster based on the importance of the other elements of the cost function. Arrival aircraft entering from the northern side of the TMA are seen to take a larger diversion than they would normally so as to avoid the Croydon/Sutton population centre. Since more outlying districts of London exist above Croydon and Sutton the population centres for this section of the map would need to be extended so as to encourage them to use the gap between Croydon and Leatherhead instead.

\begin{figure}[bhtp]
\begin{center}
\def\picwidth{12cm}
\includegraphics[width=\picwidth]{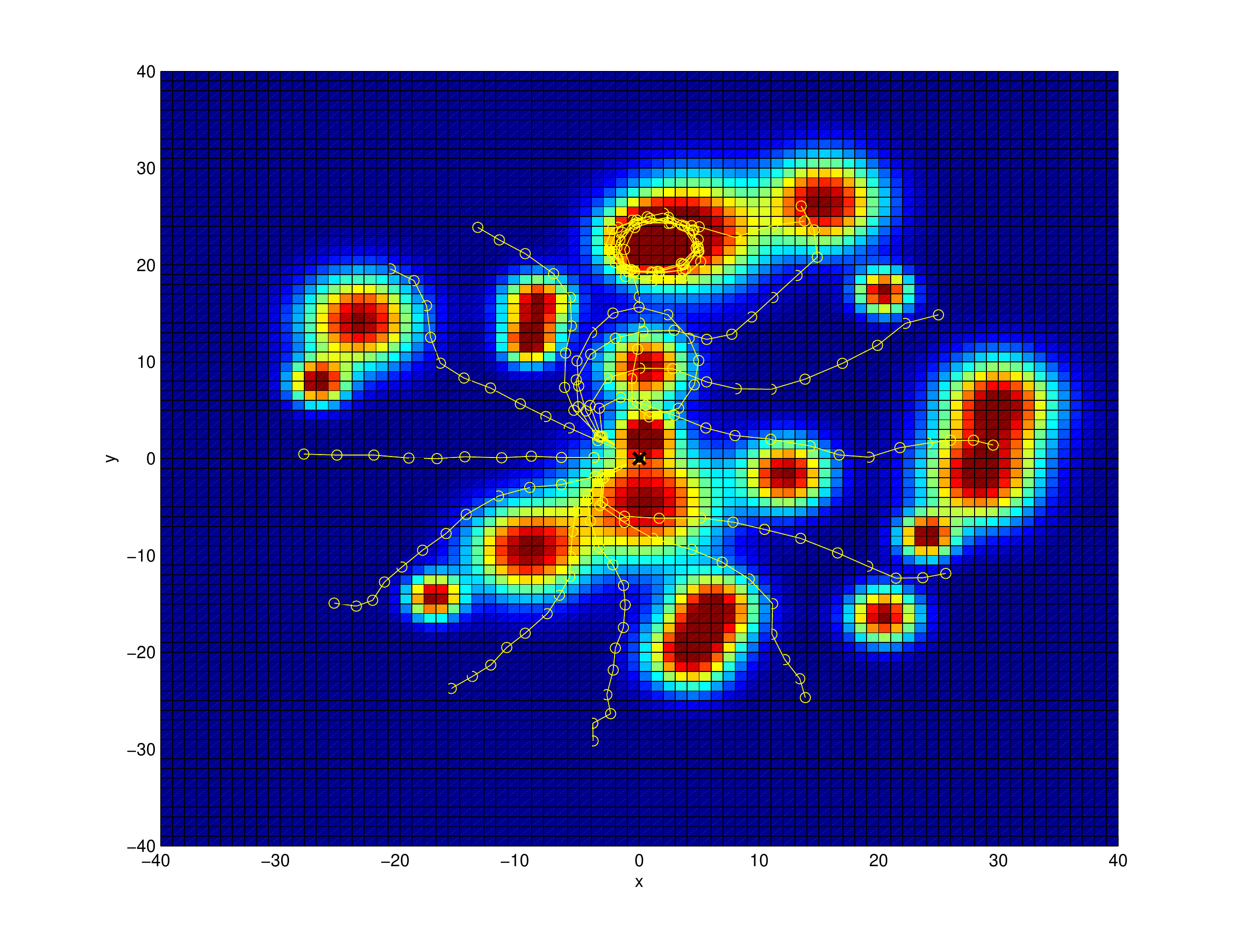}
\caption{Trajectories for 11 separate departure simulations in presence of 10\% weighted noise pollution reduction cost}
\label{fig:10depart}
\end{center}
\end{figure}
\begin{figure}[bhtp]
\begin{center}
\def\picwidth{12cm}
\includegraphics[width=\picwidth]{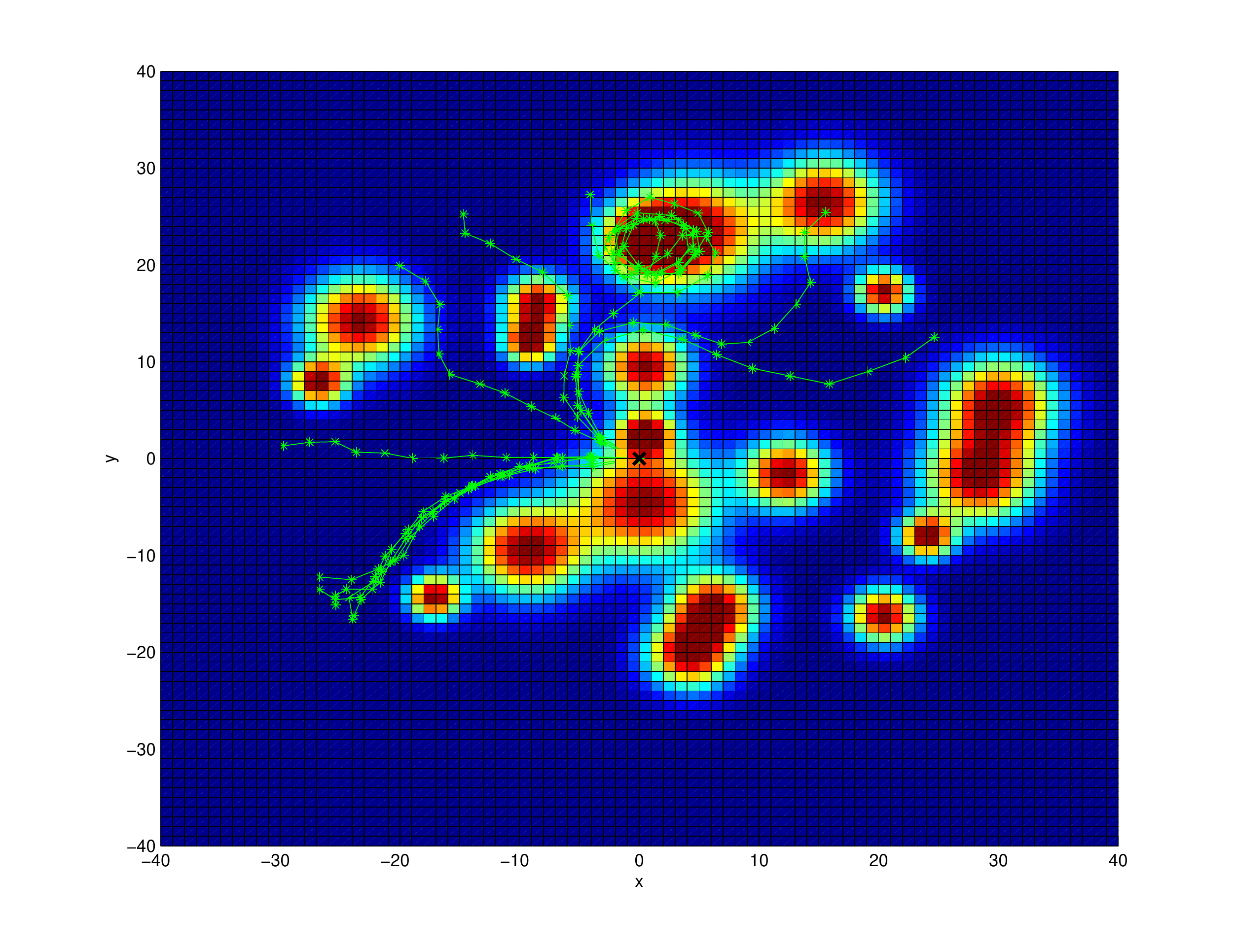}
\caption{Trajectories for 11 separate departure simulations in presence of 20\% weighted noise pollution reduction cost}
\label{fig:20depart}
\end{center}
\end{figure}
Figures \ref{fig:10depart} and \ref{fig:20depart} show the departure aircraft simulations. Here some immediate differences can be seen between both Figure \ref{fig:depart} and the 10\% and 20\% importance weightings. The 10\% importance weightings again roughly follow the pattern established in the original cost functions with the exception of the trajectory with the desired 0 \degree bearing to leave the TMA. The departure cost function contains elements which reward altitude gain, fuel minimisation, desired bearing matching and minimisation of noise pollution. No element of the function rewards gain of distance from the airport itself, instead relying on the general tendency of aircraft matching a bearing and gaining altitude to fly in a straight line along that bearing and thus leave the airport. The 0 \degree bearing simulation finds the degenerate case where the altitude gain at roughly the right bearing can be achieved by continuously circling and climbing. This behaviour is likely started by the encouragement of the noise pollution minimisation function reducing the penalty for noise pollution the higher the altitude that the aircraft passes over the population region (to 0 at the altitude cut off at 4km). In the case of the 10\% importance weighting this circling climb behaviour is continued until the maximum altitude is reached and the simulation terminated from constraint violation. This same behaviour is visible in the 20\% importance weighting although the aircraft does break the degenerate behaviour and leave the TMA before the maximum altitude constraint is reached. 

The 20\% importance weighting departures on the northern half display a preference to avoid flying over the Leatherhead population centre. Instead they take a diversion around the area rather than the tighter turn visible on the simulations with no noise pollution. Of significant interest is the affect that 20\% importance has had on the southern departures. In the 10\% importance simulations flying over Crawley's population centre seems to be considered still acceptable. Comparatively the lower population density gap between Crawley and Horsham is no longer considered valid to fly through for the 20\% importance. A secondary break in population hubs exists between Horsham and Billinghurst but by this point the aircraft are close enough to the TMA borders they gain a greater improvement in cost function by opting to leave the TMA early than turning to their desired bearing.

The noise pollution reduction cost is clearly shown to cause deviations in aircraft paths. Care is needed to maintain a balance between all the elements of the cost function such that degenerate behaviour like the uncontrolled climbing is avoided without clouding the importance of the noise pollution reduction. 

\section{Conclusions}
This report has presented the application of Sequential Monte Carlo (SMC) optimisation in the framework of Model Predictive Control (MPC) for the control of aircraft in the Terminal Manoeuvring Area (TMA). Through application of parallelisation on graphical processors (GPU) the slow method of SMC can be significantly accelerated. This allows for the solution of very complicated scenarios with both arrival and departure aircraft, in three dimensions, in the presence of a stochastic wind model and non-convex avoidance constraints.

The goal of the majority of the control optimisation done in this report is to reduce fuel usage and thus CO$_2$ production in and around the TMA. The benefits of this optimisation were demonstrated by comparison to real air traffic data around Gatwick airport over a 24 hour period of normal traffic patterns. Despite some aircraft trajectories looking far longer than real paths, significant fuel savings were achieved through improved usage of continuous descent and continuous ascent profiles. Some of this fuel saving comes from the freedom that the aircraft are given in our model to fly where they choose in the airspace. The use of a flow field to guide arrival aircraft, though very influential on our final aircraft trajectories, is far less strict than the specific corridors and beacons used in the actual airport. Our results therefore do rely on a level of flexibility and accuracy of aircraft positioning in the airspace. In total across all simulated scenarios 14653.5 kg of fuel was saved with comparison to the real data estimates. This leads to an overall fuel saving across every aircraft simulated of 29.7\% and a saving of 46158.52 kg of CO$_2$.

The model was also extended to include the effects of noise pollution. Firstly it was entirely feasible to incorporate this type of cost into the existing model. The minimal number of restrictions on the form of the cost for SMC (bounded maximisation) easily allows inclusions of complicated costs or objectives which may not have been suitable in other optimisation methods. The results from this extension also demonstrated the variability in behaviour based on the level of importance weighting noise pollution was given in the overall objective function and the importance of tuning these objective importance weightings thoroughly.

\section{Future Directions}
The work presented in this report has a wide scope for continuation. Some directions are inspired by challenging the smaller simplifying decisions taken along the way, others through the desire to reach a compromise between current working of ATC and the envisioned. This discussion of possible future directions of the work is by no means exhaustive.

\subsection{Algorithmic Improvements}
The underlying SMC algorithm was adapted to reduce the level of search space dimensionality through parallel solution of multiple agent problems. Aside from this adjustment and the minor alterations used to make the method run on a GPU kernel the base algorithm has been kept very simple. Options for algorithmic improvements could include standardizing the search patterns used by SMC. So instead of traversing the search space in a random fashion through the perturbations of controls  a hybrid of more ordered local searching could be combined with random elements in a similar fashion to simulated annealing. Work into how to best characterise the search space and how many particles are needed for a given complexity of problem would also directly benefit this work.

Currently all control variables optimised have been considered as continuous but bounded values. Although realistic from a control perspective this is not necessarily viable for real world air traffic control. Asking a pilot to set the thrust, bank and climb of an aircraft to very specific values (non round numbers) on such a rapid update cycle would be infeasible. One alternative for this would be to adapt the optimisation algorithm to operate on discrete controls. This would significantly diminish the size of the search space but would require thought on how the search space is traversed. How much discretisation of controls would also be an interesting topic to consider. The algorithm has the scope to be optimising between specific set piece moves of aircraft (such as climb 1000 feet and change heading by 10\degree) which would be very easy to relay in the current ATC environment but would impact on the potential of fuel savings due to the coarser nature of the controls. A cost benefit analysis of any discretisation would be vital in order to accurately weigh the alternatives.

Similarly in order to address the rapid update cycle required by the linearised dynamics equations vs the communications overload required to deliver these instructions to a pilot the MPC set-up could be altered. If the decision steps were placed further apart, with the addition of monitoring steps between these decisions to allow for smooth updates of the dynamics equations and constraint management, this would reduce the regularity of communication needed. However this would potentially place further conservatism on the planned trajectories for aircraft to remain within their constraints and away from avoidance regions with a reduced regularity of control. This could also allow the control horizon to be greatly extended which would give better longer term decision making for avoidance situations.

\subsection{Scenario Enhancements}
The scenarios considered in the report have had a number of limitations imposed upon them, most notably the use of a single runway which is outside of the control of the simulation itself. One clear extension would be to build in a runway management system to the optimisation either through basic runway `blocking' where an aircraft which has just used the runway blocks it from being used for a small section of time. This would need to be built in with the use of the landing envelope to predict how long from entering the landing envelope an aircraft blocks the runway. Additionally some form of priority system may need to be evolved so as to stop arrival aircraft always blocking a runway and delaying departure aircraft from taking off. This scenario enhancement would benefit from a longer control horizon or some estimate of time of arrival to plan runway movements accordingly.

Multiple runways are common in larger airports across the world and could be implemented within this work. In its simplest form some sort of decision on which runway any aircraft uses could be pre-allocated or if a runway management system is in place a more dynamic decision could be selected by the optimisation itself. This would require more investigation into how flow fields for arrival aircraft guide them to their desired runway along with emergency switching of flow fields for last minute runway changes. 

With respect to the noise pollution model extension a simple but potentially interesting adaptation to the noise pollution cost could be the variation of the population density field over time. Some airports operate periods of respite for population centres over the day, by forcing a harsh penalty for aircraft flying over a given population centre during its respite time. Modelling this would demonstrate an evolution over the day for air traffic paths based on how the underlying population density field is artificially shaped. 

\subsection{Computational Speed Ups}
Scope for additional computation speed up exists. The SMC algorithm was demonstrated as near real time for on-route traffic optimisation but the increase in complexity of scenarios for TMA optimisation necessitated an increase in particle numbers. This slowed the computational solve time. Some improvements can be made through data storage optimisation. Many values are currently stored as doubles where a float would be just as effective. How data is copied to the GPU could also be improved as this makes up a large overhead.

Some investigation into whether any benefit is gained from switching to larger or multiple GPUs along with how to optimise the spread of particles across these GPUs (similar to the analysis carried out in Section \ref{sec:compsave}) would be required.

Computational speed up could also result from a more algorithmic adjustment where the number of particles is allowed to vary throughout the SMC loops. The largest number of particles is required at the beginning of the optimisation to characterise as much of the search space as possible. As the method searches further and narrows down into areas of interest this number of active particles can be allowed to drop. This would require some changes in coding since we currently assume a fixed number of particles and thus threads for the kernels. It is however entirely possible that even with that code overhead varying particle numbers downwards through the optimisation would yield speed ups.

\section{Data access statement}

The data used for this study is archived in the University of Cambridge `DSpace@Cambridge' repository, under the title \emph{Air traffic in Gatwick TMA on 14 September 2013}. It can be accessed at\\
\url{http://www.repository.cam.ac.uk/handle/1810/248130}

\section{Acknowledgements}
The authors would like to acknowledge the funding from EPSRC (Engineering and Physical Sciences Research Council - UK) Grant No. EP/G066477/1.

Data for this study was downloaded from FlightRadar24.

The authors would also like to acknowledge their collaborators Thomas Chau and Wayne Luk from Imperial College, London for their provision of GPU facilities and expertise. Finally the authors would like to thank Thomas Kent and Arthur Richards from the University of Bristol for their expertise in accessing FlightRadar24 data.

\bibliographystyle{plainnat}
\bibliography{ajeae}

\end{document}